\baselineskip=12pt

\def\adots{\mathinner{\mskip1mu\raise1pt\hbox{.}\mskip2mu\raise4pt\hbox{.}\mskip2mu\raise7pt\vbox{\kern7pt\hbox{.}}\mskip1mu}}

\def\Hom{\mathop{\rm Hom}\nolimits}

\def\ker{\mathop{\rm ker}\nolimits}
\def\tr{\mathop{\rm tr}\nolimits}

\mathchardef\bfplus="062B
\mathchardef\bfminus="067B
\font\title=cmbx10 scaled\magstep5
\font\chapter=cmbx10 scaled\magstep4
\font\section=cmbx10 scaled\magstep2

\def\~#1{{\accent"7E #1}}
\def\sqr#1#2{{\vcenter{\hrule height.#2pt \hbox{\vrule width.#2pt height#1pt \kern#1pt \vrule width.#2pt}\hrule height.#2pt}}}
\def\square{\mathchoice\sqr63\sqr63\sqr{4.2}2\sqr{1.5}2}
\def\Square{\mathord{\hskip 1pt\square\hskip 1pt}}
\def\hk#1#2{{\vcenter{\hrule height0.0pt \hbox{\vrule width0.0pt \kern#1pt \vrule width.#2pt height#1pt}\hrule height.#2pt}}}
\def\hok{\mathchoice\hk63\hk63\hk{4.2}2\hk{1.5}2}
\def\hook{\mathbin{\hok}}
\def\operp{\mathrel{\bigcirc \hskip -11.5pt \perp}}

\def\mapright#1{\smash{\mathop{\longrightarrow}\limits^{#1}}}

\def\hookright#1{\smash{\mathop{\hookrightarrow}\limits^{#1}}}
\def\mapdown#1{\Big\downarrow\rlap{$\vcenter{\hbox{$\scriptstyle#1$}}$}}

\def\1{{\overline 1}}
\def\2{{\overline 2}}

\centerline{\chapter A Spinor Approach to Walker Geometry}
\vskip 24pt
\noindent {\section Abstract}\hfil\break
A four-dimensional Walker geometry is a four-dimensional manifold $M$ with a neutral metric $g$ and a parallel distribution of totally null two-planes. This distribution has a natural characterization as a projective spinor field subject to a certain constraint. Spinors therefore provide a natural tool for studying Walker geometry, which we exploit to draw together several themes in recent explicit studies of Walker geometry and in other work of Dunajski [11] and Pleba\~nski [30] in which Walker geometry is implicit. In addition to studying local Walker geometry, we address a global question raised by the use of spinors.
\vskip 24pt
\noindent Peter R Law. 4 Mack Place, Monroe, NY 10950, USA. prldb@member.ams.org; prldb@yahoo.com\hfil\break
Yasuo Matsushita. Section of Mathematics, School of Engineering, University of Shiga Prefecture, Hikone 522-8533, Japan. matsushita.y@e.usp.ac.jp
\vskip 24pt
\noindent 2000 MSC: 53B30, 53C27, 53C50\hfil\break
\vskip 24pt
\vfill\eject
\noindent{\section 1. Introduction}
\vskip 12pt
Our conventions and notation for the tensor and exterior algebras and curvature are stipulated in Appendix One. Let $(M,g)$ be a pseudo-Riemannian manifold of dimension $n$ and $\cal D$ a parallel distribution of $q$-planes on $M$, i.e., $\cal D$ is invariant under parallel translation. Let ${\cal D}_m$ be the plane at $m \in M$. Write ${\cal D}_m = N_m \oplus Q_m$ where $N_m := \ker\left(\xi_g\vert_{{\cal D}_m}\right)$ (where $\xi_g$ is defined just after (A1.2)) and $Q_m$ is a nondegenerate linear complement of $N_m$. If $N_m = \langle {\bf 0} \rangle_{\bf R}$, then ${\cal D}_m$ is itself nondegenerate and possesses a unique orthogonal complement ${\cal D}^\perp_m$, the latter forming a parallel $(n-q)$-distribution. The almost product structure $P$ corresponding to these complementary distributions is therefore orthogonal and parallel. Since $\nabla_g$ is torsion free, $P$ must be integrable and $(M,g,P)$ is locally decomposable pseudo-Riemannian, i.e., $M$ is locally product and each point $m$ has a neighbourhood $U$ with local coordinates $(u^{\bf A},x^{\bf B})$ which are simultaneously Frobenius for each distribution and with respect to which $(U,g)$ is a pseudo-Riemannian product; see, for example, [47].

Walker [41] studied the case when $N_m$ is nontrivial, showed one can choose local coordinates with respect to which the metric assumes a canonical form in [43], and treated the special case when ${\cal D}_m = N_m$ in [42]. When $(M,g)$ is $2n$-dimensional, $g$ is of neutral signature, and there exists a parallel distribution of totally null $n$-planes, we shall call this very natural geometry $(M,g,{\cal D})$ a {\sl Walker geometry}.

Walker geometry has proven of utility in several topics of independent interest: the holonomy Lie algebras of four-dimensional neutral metrics, [13]; isotropic K\"ahler metrics, [10] and [9]; the Osserman condition on the Jacobi operator, [2] and [8]; the Osserman condition on the conformal Jacobi operator, [4]; while Matsushita et al. [25] present a Walker geometry on an 8-torus which is a counterexample to the neutral analogue of the Goldberg conjecture, which asserts the integrability of the almost complex structure underlying an almost K\"ahler-Einstein Riemannian geometry on a compact manifold.

Of course, the structure of a Walker geometry is very natural in the context of neutral geometry, and has been investigated for its own interest: [22], [23] and [5]. In this paper, we show that four-dimensional Walker geometry has an intimate connexion with spinors and relate this fact to the local geometry of Walker four-manifolds. In the remainder of this introduction, we first make a few simple observations regarding arbitrary Walker geometry and then outline the subsequent sections of this paper. We first quote a result from [41].
\vskip 24pt
\noindent {\bf 1.1 Lemma}\hfil\break
{\it A distribution $\cal D$ on a pseudo-Riemannian manifold is parallel iff $\nabla_YX \in {\cal D}$ for all (local) sections $X$ of ${\cal D}$ and arbitrary (local) vector fields $Y$. Consequently, a parallel distribution is integrable and the integral surfaces are totally geodesic.}
\vskip 24pt
Now let $(M,g,{\cal D})$ be a Walker geometry of dimension $2n$. Walker [42] showed that one can find local coordinates $(x^1,\ldots,x^n,y^1,\ldots,y^n) =: (x^{\bf A},y^{\bf A})$, ${\bf A}=1,\ldots,n$, and hence an atlas of such, so that ${\cal D} = \langle\, \partial_{x^{\bf A}}:{\bf A}=1,\ldots,n\,\rangle_{\bf R}$ and with respect to which the metric takes the canonical form
$$\left(g_{\bf ab}\right) = \pmatrix{{\bf 0}_n&{\bf 1}_n\cr {\bf 1}_n&W\cr} \hskip 1in \left(g^{\bf ab}\right) = \pmatrix{-W&{\bf 1}_n\cr {\bf 1}_n&{\bf 0}_n\cr}\eqno(1.1)$$
where $W$ is an unspecified symmetric matrix (lower case concrete indices take values $1,\ldots,2n$ and upper case concrete indices $1,\ldots,n$). We call any coordinates, with respect to which the metric takes Walker's canonical form (1.1), {\sl Walker coordinates}. 

Consider two charts of Walker coordinates $(x^{\bf A},y^{\bf A})$ and $(u^{\bf A},v^{\bf A})$ with nontrivial intersection. Since ${\cal D} = \langle\, \partial_{x^{\bf A}}:{\bf A}=1,\ldots,n\,\rangle_{\bf R} = \langle\, \partial_{u^{\bf A}}:{\bf A}=1,\ldots,n\,\rangle_{\bf R}$, the Jacobian for the transformation between the two coordinate systems is of the form
$$J := \pmatrix{{\partial u^{\bf A} \over \partial x^{\bf B}}&{\partial u^{\bf A} \over \partial y^{\bf B}}\cr {\partial v^{\bf A} \over \partial x^{\bf B}}&{\partial v^{\bf A} \over \partial y^{\bf B}}\cr} = \pmatrix{B&C\cr {\bf 0}_n&D\cr}.$$
Since both coordinates systems are Walker, the metric has components of the form (1.1) with respect to each, i.e., on the intersection, $g = {^\tau\! J}g'J$, where $g$ and $g'$ are each of the form (1.1), from which one deduces that
$$J := \pmatrix{{\partial u^{\bf A} \over \partial x^{\bf B}}&{\partial u^{\bf A} \over \partial y^{\bf B}}\cr {\partial v^{\bf A} \over \partial x^{\bf B}}&{\partial v^{\bf A} \over \partial y^{\bf B}}\cr} = \pmatrix{B&C\cr {\bf 0}_n&{^\tau\! B}^{-1}\cr}\qquad\hbox{whence}\qquad\det(J) = 1.\eqno(1.2)$$
Note that
$$\pmatrix{{\partial x^{\bf A} \over \partial u^{\bf B}}&{\partial x^{\bf A} \over \partial v^{\bf B}}\cr {\partial y^{\bf A} \over \partial u^{\bf B}}&{\partial y^{\bf A} \over \partial v^{\bf B}}\cr} = J^{-1} = \pmatrix{B^{-1}&-B^{-1}C\,{^\tau\! B}\cr {\bf 0}_n&{^\tau\! B}\cr}.\eqno(1.3)$$
\vskip 24pt
\noindent {\bf 1.2 Lemma}\hfil\break
{\it The atlas $\{(x^{\bf A},y^{\bf A})\}$ of Walker coordinates  defines a canonical orientation for $M$, given by the orientation class $[\partial_{x^1},\ldots,\partial_{x^n},\partial_{y^1},\ldots,\partial_{y^n}]$ of the coordinate basis, and may be represented in the customary fashion by the equivalence class $[\partial_{x^1} \wedge \ldots \wedge \partial_{x^n}\wedge\partial_{y^1} \wedge\ldots\wedge\partial_{y^n}]$ in $\bigl(\Lambda^{2n}(T_mM)\bigr)/{\bf R}^+ \cong {\bf S}^0$ or the equivalence class $[dx^1 \wedge\ldots\wedge dx^n\wedge dy^1 \wedge \ldots\wedge dy^n]$ in $\left(\Lambda^{2n}\bigl((T_mM)_\bullet\bigr)\right)/{\bf R}^+ \cong {\bf S}^0$. Indeed, the globally defined $2n$-form
$$\Omega := dx^1 \wedge \ldots dx^n \wedge dy^1 \wedge \ldots dy^n\eqno(1.4)$$
is in fact the volume form for $(M,g)$ and the following $2n$-vector is also globally defined and equals the volume element $V$:
$$\partial_{x^1} \wedge \ldots \wedge \partial_{x^n} \wedge \partial_{y^1} \wedge \ldots \wedge \partial_{y^n}.$$
Thus, only orientable even-dimensional manifolds can admit a Walker geometry.}

{\it Proof.} These assertions follow immediately from (1.2--3) and
$$\left\vert\matrix{{\bf 0}_n&1_n\cr 1_n&W\cr}\right\vert = (-1)^n\left\vert\matrix{1_n&{\bf 0}_n\cr W&1_n\cr}\right\vert = (-1)^n.\eqno(1.5)$$
\vskip 24pt
\noindent {\bf 1.3 Observation}\hfil\break
Observe that $\xi_g:\partial_{x^{\bf A}} \mapsto dy^{\bf A}$ while $\xi_g:\partial_{y^{\bf A}} \mapsto dx^{\bf A} +\hbox{ terms in $dy$'s}$, whence one computes
$$\xi_g:V = \partial_{x^1} \wedge \ldots \wedge \partial_{x^n}\wedge\partial_{y^1} \wedge\ldots\wedge\partial_{y^n} \mapsto (-1)^ndx^1 \wedge\ldots\wedge dx^n\wedge dy^1 \wedge \ldots\wedge dy^n,$$
which concurs with (A1.3).

Now ${\cal D} = \langle \partial_{x^1},\ldots,\partial_{x^n} \rangle_{\bf R} = \ker\bigl(\hook (dy^1 \wedge \ldots\wedge dy^n)\bigr)$. Thus, $\cal D$ can be characterized by its image under the Pl\"ucker mapping, viz., $[\partial_{x^1} \wedge \ldots \wedge \partial_{x^n}] \in {\bf P}\bigl(\Lambda^n(TM)\bigr)$ or by $[dy^1 \wedge \ldots \wedge dy^n] \in {\bf P}\left(\Lambda^n\bigl((TM)_\bullet\bigr)\right)$.

One computes that $dy^1 \wedge \ldots \wedge dy^n$ is SD as an element of $\Lambda^n(TM_\bullet)$ with respect to the Hodge star operator defined by the canonical orientation and $g$. By (A1.6),
$$\xi^{-1}_g(dy^1 \wedge \ldots \wedge dy^n) = \xi^{-1}_g(* dy^1 \wedge \ldots \wedge dy^n) = (-1)^n* \xi^{-1}_g(dy^1 \wedge \ldots \wedge dy^n)$$
i.e.,
$$* \partial_{x^1} \wedge \ldots \wedge \partial_{x^n} = (-1)^n\partial_{x^1} \wedge \ldots \wedge \partial_{x^n}.\eqno(1.6)$$
Thus, while $dy^1 \wedge \ldots \wedge dy^n$ is always SD, $\partial_{x^1} \wedge \ldots \wedge \partial_{x^n}$ is SD/ASD according as $n$ is even/odd.
\vskip 24pt
\noindent {\bf 1.4 Observation}\hfil\break
From (1.2) and (1.3)
$$\partial_{x^1} \wedge \ldots \wedge \partial_{x^n} = \det(B)\partial_{u^1} \wedge \ldots \wedge \partial_{u^n} \hskip 1in dy^1 \wedge \ldots \wedge dy^n = \det({^\tau\! B})dv^1 \wedge \ldots \wedge dv^n$$
and $\det(B) > 0$ is the condition for the induced orientations on $\cal D$ to agree. Interchanging any $x^{\bf A}$ with any distinct $x^{\bf B}$ and $y^{\bf A}$ with $y^{\bf B}$ yields another set of Walker coordinates with the opposite orientation for $\cal D$ (but of course the same orientation for $M$). Thus, supposing $\cal D$ is indeed orientable, given a specified orientation for $\cal D$, one can choose a subatlas of Walker coordinates whose induced orientation on $\cal D$ agrees with the specified one.
\vskip 24pt
For four-dimensional Walker geometries, Walker coordinates will typically be denoted $(u,v,x,y)$ and $W$ in (1.1) is written
$$W = \pmatrix{a&c\cr c&b\cr}.\eqno(1.7)$$
In Appendix One, in addition to stipulating our conventions and some notation, we have collected together the expressions of standard objects with respect to Walker coordinates. We recast this information into spinorial form in \S 2, which focuses on the local geometry of four-dimensional Walker geometry. Appendix Two contains the essential background in spinors for four-dimensional neutral geometry we require. In \S 3, we impose a natural condition on four-dimensional Walker geometry and thereby refine a result of Dunajski [11]. In \S 4, we consider another natural restriction on Walker geometry; namely, the existence of a complementary parallel totally null distribution and demonstrate that some previously known formalisms arise naturally as special cases of Walker geometry. Finally, in \S 5, we address a global issue that arises naturally in the spinor approach to four-dimensional Walker geometry.
\vskip 24pt
\noindent {\section 2. Local Four-Dimensional Walker Geometry}
\vskip 12pt
Let $P$ be a totally null two-pane in the four-dimensional, pseudo-Euclidean linear space ${\bf R}^{2,2}$ of neutral signature. Under the identification (A2.1), any two linearly independent elements of $P$ can be written in the form $\kappa^A\lambda^{A'}$ and $\mu^A\nu^{A'}$. Orthogonality requires at least one of $\kappa^D\mu_D = 0$ and $\lambda^{D'}\nu_{D'} = 0$, i.e., either $\kappa^A \propto \mu^A$ or $\lambda^{A'} \propto \nu^{A'}$. In the latter case, $P$ takes the form $\{\,\eta^A\lambda^{A'}:\eta^A \in S\,\}$, in the former case $\{\,\kappa^A\eta^{A'}:\eta^{A'} \in S'\,\}$.

As is well known, the quadric Grassmannian $Q_2\left({\bf R}^{2,2}\right)$ of totally null two-planes in ${\bf R}^{2,2}$ is homeomorphic to {\bf O(2)}, and can be described as the planes of the form $\{\,\bigl(x,L(x)\bigr):x \in {\bf R}^2\,\}$, parametrized by $L \in {\bf O(2)}$. Under the identification (A2.1), each element, called an $\alpha$-plane, of the component $Q^+_2\left({\bf R}^{2,2}\right)\ \leftrightarrow {\bf SO(2)}$ takes the form $Z_{[\pi]} := \{\,\lambda^A\pi^{A'}:\lambda^A \in S\,\}$, for some $[\pi^{A'}] \in {\bf P}S' \cong {\bf S}^1$, while each element, called a $\beta$-plane, of the component $Q^-_2\left({\bf R}^{2,2}\right)\ \leftrightarrow\ {\bf ASO(2)}$ takes the form $W_{[\sigma]} := \{\,\sigma^A\lambda^{A'}:\lambda^{A'} \in S'\,\}$, for some $[\sigma^A] \in {\bf P}S \cong {\bf S}^1$. Under the Pl\"ucker embedding of $G_2\left({\bf R}^4\right)$ into ${\bf P}\left(\Lambda^2\left({\bf R}^4\right)\right) \cong {\bf RP}^5$, the element $Z_{[\pi]}$ of $Q^+_2\left({\bf R}^{2,2}\right)$ is mapped to the projective class of the SD bivector $\epsilon^{AB}\pi^{A'}\pi^{B'}$ while the element $W_{[\sigma]}$ of $Q^-_2\left({\bf R}^{2,2}\right)$ is mapped to the projective class of the ASD bivector $\sigma^A\sigma^B\epsilon^{A'B'}$, see (A2.2--3). 

A SD bivector $F^{ab}$ is simple iff null, $F^{ab}F_{ab} = 0$, which, with $F^{ab} = \epsilon^{AB}\psi^{A'B'}$, in turn is equivalent to $\psi^{A'B'}$ null, i.e., $\psi^{A'B'} = \pi^{A'}\pi^{B'}$, for some $\pi^{A'}$ (see [28], \S 3.5).

Thus, the Pl\"ucker embedding of $Q^\pm_2\left({\bf R}^{2,2}\right)$ is precisely the intersection of ${\bf P}\left(\Lambda^2_\pm\left({\bf R}^{2,2}\right)\right)$ with the projectivized subset of simple bivectors, which is a quadric surface in ${\bf RP}^5$. The standard basis of ${\bf R}^{2,2}$ yields, via (A1.17), a $\Psi$-ON basis for $\Lambda^2\left({\bf R}^{2,2}\right)$ which provides an explicit isomorphism $\Lambda^2\left({\bf R}^{2,2}\right) = \Lambda^2_+\left({\bf R}^{2,2}\right) \operp \Lambda^2_-\left({\bf R}^{2,2}\right) \cong {\bf R}^{1,2} \operp {\bf R}^{1,2}$ ($\operp$ denotes orthogonal direct sum) in which the space of simple SD/ASD bivectors is precisely the null cone in $\Lambda^2_\pm\left({\bf R}^{2,2}\right)$, whence the Pl\"ucker embedding of $Q^\pm_2\left({\bf R}^{2,2}\right)$ is identified with the space of generators of that null cone, i.e., with ${\bf S}^1$.

Thus, the image of an $\alpha$/$\beta$-plane under the Pl\"ucker mapping is the projective class of a SD/ASD bivector, whence they are also called SD/ASD planes. Moreover, if $P =\langle V_1,V_2 \rangle_{\bf R}$, with $\phi^i := \xi_g(V_i)$, then $\ker(\hook(\phi^1 \wedge \phi^2) = P$ and, by (A1.6), $\phi^1 \wedge \phi^2$ is SD/ASD according as $V_1 \wedge V_2$ is SD/ASD. Note that the employment of spinors to obtain this correspondence between the components of $Q_2\left({\bf R}^{2,2}\right)$ and SD/ASD simple bivectors is of course not necessary.

Now let $(M,g)$ be a connected, neutral four-dimensional manifold. Since any Walker manifold has a canonical orientation, we shall assume, without loss of generality, that $M$ is orientable. Let $\cal D$ be a totally null two-dimensional distribution on $(M,g)$. Thus, each ${\cal D}_m$ is either an $\alpha$-plane or $\beta$-plane in $T_mM$. In fact:
\vskip 24pt
\noindent {\bf 2.1 Lemma}\hfil\break
{\it $\cal D$ is a field of $\alpha$-planes or a field of $\beta$-planes.}

{\it Proof.} Since $M$ is orientable, there is a globally consistent notion of duality in $\Lambda^2(TM)$ and $\Lambda^2\bigl((TM)_\bullet\bigr)$. Let $U := \{\,m \in M:{\cal D}_m \hbox{ is SD}\,\}$. For $p \in U$, choose a neighbourhood $V$ of $p$ on which there are smooth vector fields $v$ and $w$ spanning ${\cal D}$. Then $\psi:= v \wedge w$ is SD at $p$ whence the continuous bivector $\psi + \ast\psi$ (which equals either $2\psi$ or zero) is nonzero at $p$ and thus near $p$. It follows that $U$ is open. The complement of $U$ is open by a similar argument.
\vskip 24pt
By 1.3, the distribution of a four-dimensional Walker geometry is SD with respect to the canonical orientation and, as we shall always adopt the canonical orientation, thus a distribution of $\alpha$-planes. (For any set of Walker coordinates $(u,v,x,y)$, (A1.18--20) confirms that $\partial_u \wedge \partial_v = s^+_1 - s^+_3$ is indeed SD with respect to the canonical orientation.)

Let, therefore, $\cal D$ be a distribution of $\alpha$-planes, i.e., every element of ${\cal D}_m$ is of the form $\eta^A\pi^{A'}$, with $[\pi^{A'}]$ fixed at $m$ and $\eta^A$ arbitrary. This statement only makes use of spinors locally. Supposing that $(M,g)$ is ${\bf SO^\bfplus}$-orientable, i.e., admits a reduction to ${\bf SO^\bfplus(2,2)}$, then one can construct a bundle $\cal B$ of ${\bf SO^\bfplus}$-oriented frames. Locally, one can construct a bundle of spin frames as a two-fold cover of the restriction of the bundle of ${\bf SO^\bfplus}$-oriented frames and, equally locally, associated bundles of spinors. The obstruction to gluing these local bundles together to obtain a two-fold covering of $\cal B$ and associated bundles over $M$ of spinors arises from the sign ambiguity in the two-fold covering of ${\bf SO^\bfplus(2,2)}$ by ${\bf Spin^\bfplus(2,2)}$ when nontrivial topology (specifically, the second Stiefel-Whitney class) of $M$ can obstruct a consistent choice of signs. But this problem does not arise when one employs the local lifts of transition functions for $\cal B$ to glue together local (trivial) bundles of projective classes of spinors. Thus, ${\bf P}S_M$ is well defined, provided $(M,g)$ is ${\bf SO^\bfplus}$-orientable, and unique (distinct spin structures arise from different ways of choosing signs, but the ambiguity at a point is always one of sign).

Hence, when $(M,g)$ is ${\bf SO^\bfplus}$-orientable, the distribution $\cal D$ is equivalent to a global section of ${\bf P}S'_M$. In the absence of this orientability condition, this characterization is purely local. The section $[\pi^{A'}]$, whether understood locally or globally, will be called the projective spinor field defining $\cal D$ (locally or globally).
\vskip 24pt
\noindent {\bf 2.2 Local Lifts of $[\pi^{A'}]$}\hfil\break
Let $U$ be a neighbourhood over which all bundles are trivial. Construct a trivial bundle of projective spinors on $U$ (either by restriction if ${\bf P}S'_M$ exists or otherwise by direct construction) and also construct trivial spinor bundles $S'_U$, with fibre $S'$ and $S_U$ with fibre $S$. Using the isomorphism (A2.1), one can construct an isomorphism $TM\vert_U \cong S_U \otimes S'_U$. Applying this isomorphism to a local smooth element of ${\cal D}\vert_U$ yields a spinor representation of the vector field in the form $\lambda^A\pi^{A'}$, where these spinors are defined up to scaling freedom $\pi^{A'} \mapsto f\pi^{A'}$ and $\lambda^A \mapsto f^{-1}\lambda^A$, for any nonvanishing smooth function $f$ on $U$. The spinor $\pi^{A'}$ projects onto $[\pi^{A'}]$. Any local spinor field which projects onto $[\pi^{A'}]$ will be called a {\sl a local scaled representative\/} (LSR) of $[\pi^{A'}]$.

Employing any such LSR $\pi^{A'}$ to describe $\cal D$, one easily checks that integrability of $\cal D$ is equivalent to
$$\pi_{A'}\pi^{B'}\nabla_b\pi^{A'} = 0,\eqno(2.1)$$
noting that if this equation holds for some LSR $\pi^{A'}$ then it holds for any. When integrable, the integral manifolds of $\cal D$ will be called $\alpha$-surfaces. Thus, when $\cal D$ is integrable, $M$ is foliated by $\alpha$-surfaces. 

By 1.1, a distribution $\cal D$ of $\alpha$-planes is parallel iff
$$\pi_{A'}\nabla_b\pi^{A'} = 0,\eqno(2.2)$$
which, like (2.1), if true for some LSR of $[\pi^{A'}]$ is true for any. (2.2) is equivalent to
$$\nabla_b\pi^{A'} = P_b\pi^{A'},\eqno(2.3)$$
where the one-form $P_b$ depends on the LSR $\pi^{A'}$ as follows:
$$\hbox{under}\qquad\pi^{A'} \mapsto f\pi^{A'}\hskip 1in P_b \mapsto P_b + f^{-1}\nabla_bf.\eqno(2.4)$$
Equation (2.3) expresses the fact that the LSR $\pi^{A'}$ is {\sl recurrent}, see Walker (1949). Indeed, Walker (1949), \S 5, characterized the condition for a distribution to be parallel in terms of recurrence, which, in the present circumstances, amounts to (2.3).
\vskip 24pt
The condition (2.1) in Lorentzian signature characterizes shear-free geodetic null congruences, see [29], \S 7.3. For complex spacetimes, the condition (2.1) is most usefully interpreted as describing distributions of totally null complex two-planes, see [29](7.3.18) and, for example, [31] and [3]. The geometrical interpretation of (2.1) in neutral geometry is thus a natural real analogue of that complex spacetime geometry. Walker geometry is therefore a specialization of the real neutral analogue of this complex geometry and those familiar with complex general relativity will recognize the parallels.
\vskip 24pt
\noindent{\bf 2.3 Walker's Canonical Form}\ [42]\hfil\break
Suppose that $\cal D$ is a parallel distribution of $\alpha$-planes with projective spinor field $[\pi^{A'}]$. Let $(p,q,x,y)$ be Frobenius coordinates for (the integrable) $\cal D$. Since $dx$ and $dy$ are zero when restricted to the distribution, one can write, for any LSR $\pi^{A'}$,
$$dx = \mu_A\pi_{A'} \hskip 1in dy = \nu_A\pi_{A'}\qquad\hbox{with}\qquad \nu^D\mu_D \not=0.$$
The vector fields $U^a := \mu^A\pi^{A'}$ and $V^a := \nu^A\pi^{A'}$ span $\cal D$. Noting that $\nabla_bV_c = \nabla_b\nabla_cy = \nabla_c\nabla_by = \nabla_cV_b$, then
$$\eqalign{U^b\nabla_bV^a &= g^{ac}U^b\nabla_bV_c\cr
&= g^{ac}U^b\nabla_cV_b\cr
&= g^{ac}U^b\pi_{B'}(\nabla_c\nu_B) + g^{ac}U^b\nu_B(\nabla_c\pi_{B'})\qquad\hbox{(the first summand of which is zero)}\cr
&= U^b\nu_B\pi_{B'}P^a = 0.\cr}$$
\vskip 6pt
Similarly, $V^b\nabla_bU^a = 0$ and it follows that $[U,V] = U^b\nabla_bV^a - V^b\nabla_bU^a = 0$. The pair of equations $Uf = 1$, $Vf = 0$ have trivial integrability conditions (see [42]), as do the equations $Ug = 0$, $Vg = 1$. Let $u$ and $v$ be solutions of these systems respectively. Let $B$ be the nonsingular matrix expressing $\{\partial_p,\partial_q\}$ in terms of $\{U,V\}$ (both are frames of $\cal D$). Then, by the definition of $u$ and $v$,
$${\partial(u,v,x,y) \over \partial(p,q,x,y)} = \pmatrix{B&C\cr {\bf 0}_2&1_2\cr}$$
is nonsingular, whence $(u,v,x,y)$ are legitimate local coordinates. Computing the metric $g^{ab}$: it vanishes on $\langle dx,dy \rangle_{\bf R}$; $g(du,dx) = du(U) = U(u) = 1$; $g(du,dy) = du(V) = V(u) = 0$; $g(dv,dx) = dv(U) = U(v) = 0$; and $g(dv,dy) = dv(V) = V(v) = 1$; which together give the form (1.1) of Walker's canonical form.

From Walker's canonical form for these coordinates, one observes that $U^a := g^{ab}(dx)_b = \partial_u$ and $V^a := g^{ab}(dy)_b = \partial_v$, i.e., the Walker coordinates are Frobenius coordinates for $\cal D$ satisfying
$$\partial_u = \mu^A\pi^{A'} \hskip .5in \partial_v = \nu^A\pi^{A'} \hskip .5in dx = \mu_A\pi_{A'} \hskip .5in dy = \nu_A\pi_{A'}\qquad \hbox{with }\nu^D\mu_D \not=0.\eqno(2.5)$$
\vskip 24pt
\noindent {\bf 2.4 Remarks}\hfil\break
From (2.5)
$$dx \wedge dy = (\nu^D\mu_D)\epsilon_{AB}\pi_{A'}\pi_{B'} =: (\nu^D\mu_D)\sigma_\pi \hskip 1in \partial_u \wedge \partial_v = (\nu^D\mu_D)\epsilon^{AB}\pi^{A'}\pi^{B'} =: (\nu^D\mu_D)\Sigma_\pi\eqno(2.6)$$
If one replaces $\pi^{A'}$ by $f\pi^{A'}$, then one must replace $\mu$ by $\mu/f$, $\nu$ by $\nu/f$, whence $\nu^D\mu_D$ by $(\nu^D\mu_D)/f^2$. Thus, a suitable scaling yields an LSR $\pi^{A'}$ so that $\nu^D\mu_D = \pm 1$, whence $dx \wedge dy = \pm \sigma_\pi$ and $\partial_u \wedge \partial_v = \pm\Sigma_\pi$. 

Thus, given any Walker coordinates $(u,v,x,y)$, there is an LSR $\pi^{A'}$ such that
$$dx \wedge dy = \pm\sigma_\pi, \hskip 1in \partial_u \wedge \partial_v = \pm\Sigma_\pi\eqno(2.7)$$
and this LSR is unique up to sign.

If $\cal D$ is orientable, then one can choose an atlas of Walker coordinates so that $\partial_u \wedge \partial_v$ defines a consistent orientation on $\cal D$, see 1.4. Thus, for an orientable distribution, one can choose an atlas of Walker coordinates so that (2.7) holds with constant sign for all charts of the atlas. We will have more to say about orientability of $\cal D$ in \S 5.

If $(u,v,x,y)$ are Walker coordinates, by (A1.7) so are $(v,u,y,x)$. Hence, by breaking the Walker symmetry (A1.7), it is always possible to choose Walker coordinates $(u,v,x,y)$ for which the plus sign occurs in (2.7). For such Walker coordinates, and the LSR $\pi^{A'}$ determined by (2.7), we will write 
$$\partial_u = \alpha^A\pi^{A'} \hskip .5in \partial_v = \beta^A\pi^{A'} \hskip .5in dx = \alpha_A\pi_{A'} \hskip .5in dy = \beta_A\pi_{A'},\eqno(2.8)$$
where $\{\alpha^A,\beta^A\}$ is an (unprimed) spin frame.
\vskip 24pt
\noindent {\bf 2.5 Proposition}\hfil\break
{\it The projective spinor $[\pi^{A'}]$ of a Walker geometry $(M,g,{\cal D},[\pi^{A'}])$ is a Weyl Principal Spinor {\rm (WPS, see [19])} of multiplicity at least two. Indeed, any LSR $\pi^{A'}$ satisfies
$$\tilde\Psi_{A'B'C'D'}\pi^{C'}\pi^{D'} + 2\Lambda\pi_{A'}\pi_{B'} = 0.\eqno(2.9)$$
Thus, $[\pi^{A'}]$ is a WPS of multiplicity at least three in scalar-flat Walker geometries and ${\cal W}^+ = 0\ \Rightarrow\ S = 0$. Moreover, $[\pi^{A'}]$ is also a principal spinor of the Ricci spinor $\Phi_{ABA'B'}$, whence the Einstein endomorphism $G^a{}_b = 2\Phi^{AA'}{}_{BB'}$ maps $\cal D$ to itself {\rm (see A1.8}.}

{\it Proof.} Working with some LSR $\pi^{A'}$ of $[\pi^{A'}]$,
$$0 = \nabla^b\left(\pi_{A'}\nabla_b\pi^{A'}\right) = \pi_{A'}\Square\,\pi^{A'} + \left(\nabla^b\pi_{A'}\right)\left(\nabla_b\pi^{A'}\right).$$
Since the second term vanishes by (2.3), then
$$0 = \pi_{A'}\Square\,\pi^{A'}.\eqno(2.10)$$
Furthermore,
$$0 = \nabla^B_{C'}\left(\pi^{A'}\nabla_b\pi_{A'}\right) = \pi^{A'}\nabla^B_{C'}\nabla_b\pi_{A'} + \left(\nabla^B_{C'}\pi^{A'}\right)\left(\nabla_b\pi_{A'}\right) = \pi^{A'}\nabla^B_{C'}\nabla_b\pi_{A'},$$
since the second term in the sum again vanishes by (2.3). Using (A2.8),
$$\nabla^B_{C'}\nabla_b = -\nabla_{BC'}\nabla^B_{B'} = -\left[\Square_{C'B'} + \nabla_{B[C'}\nabla^B_{B']}\right] = -\left[\Square_{C'B'} + {1 \over 2}\epsilon_{C'B'}\Square\right].$$
Thus, one obtains from the previous calculation
$$\eqalign{0 &= \pi^{A'}\left[\Square_{C'B'}\pi_{A'} + {1 \over 2}\epsilon_{C'B'}\Square\,\pi_{A'}\right]\cr
&= \pi^{A'}\Square_{C'B'}\pi_{A'} + {1 \over 2}\epsilon_{C'B'}\pi^{A'}\Square\,\pi_{A'}\cr
&= \pi^{A'}\left[\tilde\Psi_{C'B'A'D'}\pi^{D'} - \Lambda\pi_{C'}\epsilon_{B'A'} - \Lambda\epsilon_{C'A'}\pi_{B'}\right]\qquad\hbox{by (A2.9) and (2.10)}\cr
&= \tilde\Psi_{B'C'A'D'}\pi^{A'}\pi^{D'} + 2\Lambda\pi_{B'}\pi_{C'}.\cr}$$
Transvecting by $\pi^{C'}$ shows that $[\pi^{A'}]$ is a WPS of multiplicity at least two.

Similarly,
$$\eqalign{0 &= \nabla^{B'}_C(\pi^{A'}\nabla_b\pi_{A'})\cr
&= \pi^{A'}\nabla^{B'}_C\nabla_b\pi_{A'} + \left(\nabla^{B'}_C\pi^{A'}\right)\nabla_b\pi_{A'}\cr
&= -\pi^{A'}\nabla_{B'C}\nabla_B{}^{B'}\pi_{A'}\cr
&= - \pi^{A'}\left[\Square_{CB}\pi_{A'} + {1 \over 2}\epsilon_{CB}\Square\pi_{A'}\right]\cr
&= -\pi^{A'}\Square_{CB}\pi_{A'}\cr
&= -\Phi_{CBA'D'}\pi^{A'}\pi^{D'}.\cr}$$
(One can write
$$\Phi_{ABA'B'} = A_{AB}\pi_{A'}\pi_{B'} + B_{AB}\pi_{(A'}\xi_{B')} + C_{AB}\xi_{A'}\xi_{B'},$$
where $\xi^{D'}\pi_{D'} = 1$ and $A_{AB}$, $B_{AB}$, and $C_{AB}$ are symmetric. $\Phi_{ABA'B'}\pi^{A'}\pi^{B'} = 0$ entails $C_{AB} = 0$.)
\vskip 24pt
\noindent {\bf 2.6 Classification of the SD Weyl Curvature Endomorphism}\hfil\break
The classification of the Weyl curvature endomorphisms of four-dimensional neutral metrics according to their Jordan canonical form (JCF) was given in [17] and according to the algebraic structure of the corresponding Weyl spinors in [19]. By 2.5, the SD Weyl spinor of any four-dimensional Walker geometry is algebraically special. From (A1.27), the eigenvalues of the SD Weyl curvature endomorphism ${\cal W}^+$ of any four-dimensional Walker geometry are $-S/6$, $S/12$, and $S/12$ and D\'{\i}az-Ramos et al. [8] determined the Jordan canonical forms of ${\cal W}^+$. Here we refine their result. For any eigenvalue, let $m$($M$) denote the algebraic (geometric) multiplicity.  There are four cases; see (A1.23--24) for notation:\hfil\break

\noindent i) $S \not= 0$ but $S^2 + AS + 3B^2 = 0$; ${\cal W}^+$ can be diagonalized and has JCF
$$-{S \over 12}\pmatrix{2&0&0\cr 0&-1&0\cr 0&0&-1\cr},$$
$-S/6$ has $m=M=1$; $S/12$ has $m=M=2$. Since $[\pi^{A'}]$ is a real WPS of multiplicity at least two, from [19], p. 2106, this information indicates ${\cal W}^+$ must be of type $\{22\}$Ia, and the multiplicity of $[\pi^{A'}]$ is exactly two.\hfil\break
ii) $S \not= 0$ and $S^2 + AS + 3B^2 \not= 0$; ${\cal W}^+$ has JCF
$$\pmatrix{-{S \over 6}&0&0\cr 0&{S \over 12}&1\cr 0&0&{S \over 12}\cr},$$
$-S/6$ has $m=M=1$; $S/12$ has $M=1$, $m=2$. It follows from [19], p. 2106, that ${\cal W}^+$ is of type $\{211\}$II or $\{1\1 2\}$II; in each case the double WPS must be $[\pi^{A'}]$.\hfil\break
iii) $S=0$, $B=0$, $A \not= 0$; from (A1.25), the matrix representation of ${\cal W}^+$ is
$$^+{\bf W} = -{A \over 12}\pmatrix{1&0&1\cr 0&0&0\cr -1&0&-1\cr},\qquad\hbox{with JCF}\qquad\pmatrix{0&0&0\cr 0&0&1\cr 0&0&0\cr},$$
and 0 is the only eigenvalue, with $M=2$. From [19], p. 2106, ${\cal W}^+$ must be type $\{4\}$II, i.e., $\tilde\Psi_{A'B'C'D'}$ is null with four-fold WPS $[\pi^{A'}]$.\hfil\break
iv) $S=0$ but $B \not= 0$; ${\cal W}^+$ has JCF $J_3(0)$, i.e., 0 is the only eigenvalue and $M=1$. From [19], p. 2106, ${\cal W}^+$ is of type $\{31\}$III and $[\pi^{A'}]$ is a WPS of multiplicity three.
\vskip 24pt
\noindent {\bf 2.7 Remark}\hfil\break
Since ${\cal W}^+$ is algebraically special, it is never, in particular, of type Ib, see [17] and [19]. If, therefore, $(M,g)$ is compact Einstein, then $\chi(M) \leq 3\tau(M)/2$, see [17], [26]. If further, ${\cal W}^-$ is not of type Ib, then in fact $\chi(M) \leq -3\vert \tau(M)\vert/2 \leq 0$. This conclusion holds if $(M,g)$ is a compact Einstein Walker geometry which is, for example: SD (${\cal W}^-=0$); or algebraically special in ${\cal W}^-$; or K\"ahler with the orientation induced by the complex structure agreeing with the Walker canonical orientation, in which case ${\cal W}^-$ is of type Ia, see [15]. Examples of four-dimensional Walker geometries with ASD Weyl curvature of type Ib are presented in [4].
\vskip 24pt
In 2.6 we relied upon the computations recorded in Appendix One and the result of D\'{\i}az-Ramos et al. [8] which is also based on these computations. A systematic development of spinor analysis of Walker geometry would proceed by computing the neutral analogues of spin coefficients, [28], \S 4.5, and thence the spinor equivalents of the curvature. The notation for spin coefficients for neutral signature, however, requires modification of that employed in Lorentz signature; e.g., the priming operation of [28](4.5.17) is not appropriate for neutral signature. Spin coefficients for neutral signature will be presented elsewhere; here we shall follow expediency and further exploit the known results of Appendix One to deduce the spinor equivalents of curvature.

To this end, Walker's canonical form for the metric
$$ds^2 = 2dx(du + {a \over 2}dx + {c \over 2}dy) - 2dy(-dv - {c \over 2}dx - {b \over 2}dy),$$
suggests utilizing the following null tetrad:
$$\ell_a  := dx \hskip .5in \tilde m_a := dy \hskip .5in n_a := du + {a \over 2}dx + {c \over 2}dy \hskip .5in m_a := -(dv + {c \over 2}dx + {b \over 2}dy),$$
which we shall call the {\sl Walker null tetrad\/} associated to a set of Walker coordinates.

The null tetrad determines unique (up to an overall sign) spin frames. In particular, assuming we have chosen Walker coordinates $(u,v,x,y)$ and a LSR $\pi^{A'}$ so that (2.8) is satisfied, write $n^a = \nu^A\zeta^{A'}$, where $\zeta^{D'}\pi_{D'} \not = 0$. Then $n^a\tilde m_a = 0$ implies $\nu^D\beta_D = 0$. So $n^a = \beta^A\xi^{A'}$, with $\xi^{D'}\pi_{D'} \not= 0$. But then $n^a\ell_a = 1$ implies $\xi^{D'}\pi_{D'} = 1$. Writing $m^a = \gamma^A\kappa^{A'}$, then $m^a\tilde m_a = -1$ but $m^a\ell_a = 0$ implies $\gamma^D\alpha_D = 0$. Then, $m^an_a = 0$ implies $\kappa^{D'}\xi_{D'} = 0$, whence $m^a = \lambda\alpha^A\xi^{A'}$, for some $\lambda \in {\bf R}$. But then $m^a\tilde m_a = -1$ implies $\lambda = 1$. Thus, $\{\alpha^A,\beta^A\}$ and $\{\pi^{A'},\xi^{A'}\}$ are the unique (up to an overall sign) spin frames associated to the null tetrad, and
$$\vcenter{\openup2\jot \halign{$\hfil#$&&${}#\hfil$&\qquad$\hfil#$\cr
\ell_a &= dx = \alpha_A\pi_{A'} & \tilde m_a &= dy = \beta_A\pi_{A'}\cr
n_a &= \displaystyle du + {a \over 2}dx + {c \over 2}dy = \beta_A\xi_{A'} & m_a &= \displaystyle-(dv + {c \over 2}dx + {b \over 2}dy) = \alpha_A\xi_{A'}\cr
\ell^a &= \partial_u = \alpha^A\pi^{A'} & \tilde m^a &= \partial_v = \beta^A\pi^{A'}\cr
n^a &= \displaystyle-{a \over 2}\partial_u - {c \over 2}\partial_v + \partial_x = \beta^A\xi^{A'} & m^a &= \displaystyle{c \over 2}\partial_u + {b \over 2}\partial_v - \partial_y = \alpha^A\xi^{A'}\cr}}\eqno(2.11)$$
\vskip 6pt
Note that $[\partial_u,\partial_v,\partial_x,\partial_y] = [\ell^a,\tilde m^a,n^a,-m^a]$, the latter being a Witt frame which reduces to the standard Witt frame for ${\bf R}^{2,2} \cong {\bf R}^4_{\rm hb}$ when $a=b=c=0$ (${\bf R}^4_{\rm hb}$ being the standard {\sl hyperbolic\/} four-dimensional pseudo-Euclidean space as in [32]).

From the null tetrad, one constructs a $\Psi$-ON basis as follows:
$$\eqalignno{U^a &:= {1 \over \sqrt 2}(\ell^a + n^a) = {1 \over \sqrt 2}\left({2-a \over 2}\partial_u - {c \over 2}\partial_v + \partial_x\right)\cr
\noalign{\vskip 6pt}
V^a &:= {1 \over \sqrt 2}(\tilde m^a - m^a) = {1 \over \sqrt 2}\left(-{c \over 2}\partial_u + {2-b \over 2}\partial_v + \partial_y\right)\cr
\noalign{\vskip -6pt}
&&(2.12)\cr
\noalign{\vskip -6pt}
X^a &:= {1 \over \sqrt 2}(\ell^a - n^a) = {1 \over \sqrt 2}\left({2+a \over 2}\partial_u + {c \over 2}\partial_v - \partial_x\right)\cr
\noalign{\vskip 6pt}
Y^a &:= {1 \over \sqrt 2}(\tilde m^a + m^a) = {1 \over \sqrt 2}\left({c \over 2}\partial_u + {2+b \over 2}\partial_v - \partial_y\right).\cr}$$
\vskip 6pt
\noindent We note that $[U^a,V^a,X^a,Y^a] = [\partial_u,\partial_v,\partial_x,\partial_y]$, i.e., the $\Psi$-ON frame (2.12) has the canonical orientation; moreover it is well behaved under the Walker symmetry (A1.7): $U^a \leftrightarrow V^a$ and $X^a \leftrightarrow Y^a$.

From (2.12) and (A1.17), one obtains the following $\Psi$-ON bases for the spaces of SD and ASD bivectors:
$$\eqalignno{S^+_1 &:= {U^a \wedge V^b + X^a \wedge Y^b \over \sqrt 2}\cr 
&= {1 \over 2\sqrt 2}\left({4+ab-c^2 \over 2}\partial_u \wedge \partial_v + c\partial_u \wedge \partial_x - a\partial_u \wedge \partial_y + b\partial_v \wedge \partial_x - c\partial_v \wedge \partial_y + 2\partial_x \wedge \partial_y\right);\cr
\noalign{\vskip 6pt}
S^+_2 &:= {U^a \wedge X^b + V^a \wedge Y^b \over \sqrt 2} = -{1 \over \sqrt 2}(\partial_u \wedge \partial_x + \partial_v \wedge \partial_y);\cr
\noalign{\vskip 6pt}
S^+_3 &:= {U^a \wedge Y^b - V^a \wedge X^b \over \sqrt 2}&(2.13)\cr
&= {1 \over 2\sqrt 2}\left({4 - ab + c^2 \over 2}\partial_u \wedge \partial_v - c\partial_u \wedge \partial_x + a\partial_u \wedge \partial_y - b\partial_v \wedge \partial_x + c\partial_v \wedge \partial_y - 2\partial_x \wedge \partial_y\right);\cr
\noalign{\vskip 6pt}
S^-_1 &:= {U^a \wedge V^b - X^a \wedge Y^b \over \sqrt 2} = {1 \over \sqrt 2}\left(-{(a+b) \over 2}\partial_u \wedge \partial_v + \partial_u \wedge \partial_y - \partial_v \wedge \partial_x\right);\cr
\noalign{\vskip 6pt}
S^-_2 &:= {U^a \wedge X^b - V^a \wedge Y^b \over \sqrt 2} = {1 \over \sqrt 2}(c\partial_u \wedge \partial_v - \partial_u \wedge \partial_x + \partial_v \wedge \partial_y);\cr
\noalign{\vskip 6pt}
S^-_3 &:= {U^a \wedge Y^b + V^a \wedge X^b \over \sqrt 2} = {1 \over \sqrt 2}\left({b-a \over 2}\partial_u \wedge \partial_v - \partial_u \wedge \partial_y - \partial_v \wedge \partial_x\right).\cr}$$
Writing these (A)SD bivectors in terms of spinors:
$$\eqalignno{S^+_1 &= {\ell^a \wedge \tilde m^b - n^a \wedge m^b \over \sqrt 2} = {\epsilon^{AB}(\pi^{A'}\pi^{B'} + \xi^{A'}\xi^{B'}) \over \sqrt 2} =: \epsilon^{AB}\triangle_1^{A'B'};\cr 
S^-_1 &= {n^a \wedge \tilde m^b - \ell^a \wedge m^b  \over \sqrt 2} = -{(\alpha^A\alpha^B +\beta^A\beta^B)\epsilon^{A'B'} \over \sqrt 2} =: - \triangle_1^{AB}\epsilon^{A'B'};\cr
\noalign{\vskip 12pt}
S^\pm_2 &= {n^a \wedge \ell^b \pm \tilde m^a \wedge m^b \over \sqrt 2} = \cases{-\sqrt 2 \epsilon^{AB}\pi^{(A'}\xi^{B')} =: -\epsilon^{AB}\triangle_3^{A'B'};&\cr \cr -\sqrt 2\alpha^{(A}\beta^{B)}\epsilon^{A'B'} =: -\triangle_3^{AB}\epsilon^{A'B'};&\cr}&(2.14)\cr
\noalign{\vskip 12pt}
S^+_3 &= {\ell^a \wedge \tilde m^b + n^a \wedge m^b \over \sqrt 2} = {\epsilon^{AB}(\pi^{A'}\pi^{B'} - \xi^{A'}\xi^{B'}) \over \sqrt 2} =: \epsilon^{AB}\triangle_2^{A'B'};\cr
S^-_3 &= {n^a \wedge \tilde m^b + \ell^a \wedge m^b  \over \sqrt 2} = {(\alpha^A\alpha^B - \beta^A\beta^B)\epsilon^{A'B'} \over \sqrt 2} =: \triangle_2^{AB}\epsilon^{A'B'}.\cr}$$
Noting the conventions of Appendix One, and putting
$$^+Z^{ab}_i := \epsilon^{AB}\triangle^{A'B'}_i \hskip 1in ^-Z_i^{ab} := \triangle^{AB}\epsilon^{A'B'},$$
the SD Weyl curvature endomorphism ${\cal W}^+$ satisfies:
$$\eqalignno{^+Z_i \,{^+{\cal C}}^i{}_j &:= {\cal W}^+({^+Z}_j)\cr
&:= {1 \over 2}\epsilon^{AB}\epsilon_{CD}\tilde\Psi^{A'B'}{}_{C'D'}\,\epsilon^{CD}\triangle_j^{C'D'}\cr
&= \epsilon^{AB}\tilde\Psi^{A'B'}{}_{C'D'}\triangle_j^{C'D'}\cr
&=: \epsilon^{AB}\triangle_i^{A'B'}{\bf\tilde\Psi}^i{}_j&(2.15)\cr}$$
i.e., the matrix representation $^+{\cal C}$ of ${\cal W}^+$ (acting on $\Lambda^2_+$, the space of SD bivectors) with respect to the basis $\{{^+Z}_1,{^+Z}_2,{^+Z}_3\}$ coincides with the matrix representation ${\bf\tilde\Psi}$ of the endomorphism $\tilde\Psi^{A'B'}{}_{C'D'}$ (acting on the space $S' \odot S'$ of symmetric rank two primed spinors) with respect to the basis $\{\triangle^{A'B'}_1,\triangle^{A'B'}_2,\triangle^{A'B'}_3\}$. Similarly, the matrix representation $^-{\cal C}$ of the ASD Weyl curvature endomorphism ${\cal W}^-$ (acting on $\Lambda^2_-$) with respect to the basis $\{{^-Z}_1,{^-Z}_2,{^-Z}_3\}$ coincides with the matrix representation ${\bf\Psi}$ of the endomorphism $\Psi^{AB}{}_{CD}$ (acting on $S \odot S$) with respect to the basis $\{\triangle^{AB}_1,\triangle^{AB}_2,\triangle^{AB}_3\}$.

From the matrix representations of ${\cal W}^+$ and ${\cal W}^-$ given in (A1.21--25) with respect to the bases (A1.19--20) of $\Lambda^2_+$ and $\Lambda^2_-$, one can calculate $^+{\cal C}$ and $^-{\cal C}$. From the expression for ${\bf\Psi}$ given in [19], equation (20), and its analogue for ${\bf \tilde\Psi}$, one can then compute $\Psi_{ABCD}$ and $\tilde\Psi_{A'B'C'D'}$ themselves.

From (2.13) and (A1.19--20), one finds
$$\displaylines{S^+_1 = {1 \over 4\sqrt 2}\left((c^2+5)s^+_1 - 2cs^+_2 - (c^2+3)s^+_3\right)\hskip .75in S^+_3 = {1 \over 4\sqrt 2}\left((3-c^2)s^+_1 + 2cs^+_2 + (c-5^2)s^+_3\right)\cr
\hfill S^+_2 = {1 \over \sqrt 2}(cs^+_1 - s^+_2 - cs^+_3)\hfill\llap(2.16)\cr
S^-_1 = {s^-_1 \over \sqrt 2} \hskip 1in S^-_2 = -{s^-_2 \over \sqrt 2} \hskip 1in S^-_3 = -{s^-_3 \over \sqrt 2}.\cr}$$
Now, (2.14) and (2.16) entail:
$$\displaylines{\hskip 2in {^+Z_1} = {1 \over 4\sqrt2}\bigl((c^2+5)s^+_1 - 2cs^+_2 - (c^2+3)s^+_3\bigr)\hfill\cr
\hskip 2in{^+Z_2} = {1 \over 4\sqrt2}\bigl((3-c^2)s^+_1 + 2cs^+_2 + (c-5^2)s^+_3\bigr)\hfill\llap(2.17)\cr
\hskip 2in {^+Z_3} = {1 \over \sqrt2}(-cs^+_1 + s^+_2 + cs^+_3).\hfill\cr}$$
The matrix $J$ expressing the $^+Z_i$ in terms of the $s^+_i$ is therefore
$$J = {1 \over 4\sqrt2}\pmatrix{c^2+5&3-c^2&-4c\cr -2c&2c&4\cr -(3+c^2)&c^2-5&4c\cr}\qquad\hbox{whence}\qquad J^{-1} = {\sqrt2 \over 4}\pmatrix{c^2+5&2c&c^2+3\cr c^2-3&2c&c^2-5\cr 4c&4&4c\cr}.\eqno(2.18)$$
Hence, from (A1.25),
$$\eqalignno{{^+{\cal C}} &= J^{-1}\,{^+{\bf W}}\,J\cr
&= -{1 \over 48}\pmatrix{A - 6Bc - 3S(c^2+1)&-A+6Bc+S(3c^2-1)&6(B+Sc)\cr A - 6Bc - S(3c^2-1)&-A + 6Bc + S(3c^2-5)&6(B+Sc)\cr -6(B+Sc)&6(B+Sc)&8S\cr}.&(2.19)\cr}$$
Equating this expression to the tilde version of [19](20) one finds:
$$\tilde\Psi_0 = 0 = \tilde\Psi_1 \hskip .5in \tilde\Psi_2 = {S \over 12} \hskip .5in \tilde\Psi_3 = -{(B+Sc) \over 8} \hskip .5in \tilde\Psi_4 = {6Bc-A+S(3c^2-1) \over 24},\eqno(2.20)$$
whence
$$\tilde\Psi_{A'B'C'D'} = {S \over 2}\pi_{(A'}\pi_{B'}\xi_{C'}\xi_{D')} + {B+Sc \over 2}\pi_{(A'}\pi_{B'}\pi_{C'}\xi_{D')} + {6Bc-A+S(3c^2-1) \over 24}\pi_{A'}\pi_{B'}\pi_{C'}\pi_{D'}.\eqno(2.21)$$
From (2.20--21), one can obtain the results of 2.6 directly. Referring to [19](22--24), one computes $I = S^2/24$, $J = -S^3/288$, whence $I^3 = 6J^2$, and $0 = (\lambda + S/6)(\lambda - S/12)^2$ is the eigenvalue equation for $\tilde\Psi^{A'B'}{}_{C'D'}$. The expression (2.21) confirms 2.5: $[\pi_{A'}]$ is a real WPS of multiplicity at least two, and of multiplicity at least three when $S=0$. If $[\pi_{A'}]$ is of multiplicity exactly two, then $S \not= 0$. Given the eigenvalues just deduced, from the diagram in [19], p. 2106, if the geometric multiplicities coincide with the algebraic multiplicities, then $\tilde\Psi_{A'B'C'D'}$ is type $\{22\}$Ia (in particular, ${\cal W}^+$ is diagonalizable). It follows from [19], \S 5.6, that $\tilde\Psi_{A'B'C'D'} = 6(S/12)\pi_{(A'}\pi_{B'}\eta_{C'}\eta_{D')}$, for some spinor $\eta_{A'} = p\pi_{A'} + q\xi_{A'}$. Equating this expression with (2.21) yields the condition $S^2 + AS + 3B^2 = 0$ of 2.6(i). When $S^2 + AS + 3B^2 \not= 0$, $\tilde\Psi_{A'B'C'D'}$ cannot be of type $\{22\}$Ia; with $S \not= 0$ still, one sees from [19], p. 2106, that the only possible types with $[\pi_{A'}]$ a real WPS of multiplicity two, are types $\{211\}$II or $\{1\1 2\}$II, in which cases the geometric multiplicity of $S/12$ is one, rather than two (in particular, ${\cal W}^+$ is not diagonalizable). If now $S = 0$ but $B \not= 0$, then $[\pi_{A'}]$ is a WPS of multiplicity three, whence $\tilde\Psi_{A'B'C'D'}$ has type $\{31\}$III, and there is a single eigenvalue (namely zero) of algebraic multiplicity three and geometric multiplicity one. Finally, if $S = 0$, $B=0$, but $A \not= 0$, then $[\pi_{A'}]$ is a WPS of multiplicity four, whence $\tilde\Psi_{A'B'C'D'}$ is of type $\{4\}$II, and the zero eigenvalue now has geometric multiplicity two. Thus, the results of 2.6 are obtainable directly from (2.21) (bearing in mind that this expression could be computed directly by first computing spin coefficients without exploiting the results of Appendix One).

Turning now to the ASD Weyl curvature, one finds
$$^-Z_1 = -S^-_1 = -{s^-_1 \over \sqrt 2} \hskip .75in {^-Z_2} = S^-_3 = -{s^-_3 \over \sqrt 2} \hskip .75in {^-Z_3} = -S^-_2 = {s^-_2 \over \sqrt 2},\eqno(2.22)$$
whence the analogue of (2.18) is
$$K = {1 \over \sqrt 2}\pmatrix{-1&0&0\cr 0&0&1\cr 0&-1&0\cr}\qquad\hbox{and}\qquad K^{-1} = \sqrt 2\pmatrix{-1&0&0\cr 0&0&-1\cr 0&1&0\cr}.\eqno(2.23)$$
Hence, referring to (A1.22), one computes
$${^-{\cal C}} = K^{-1}\,{^-{\bf W}}\,K = {1 \over 12}\pmatrix{{\cal P}+3{\cal Q}&-3{\cal Y}&3({\cal T}+{\cal X})\cr 3{\cal Y}&{\cal P}-3{\cal Q}&3({\cal T}-{\cal X})\cr-3({\cal T}+{\cal X})&3({\cal T}-{\cal X})&-2{\cal P}\cr}.\eqno(2.24)$$
Equating this expression to [19](20) and, as in Appendix One, using the numerals 1, 2, 3 ,4 to denote $u$, $v$, $x$, $y$, yields:
$$\Psi_0 = {b_{11} \over 2} \hskip .5in \Psi_1 = {b_{12} - c_{11} \over 4} \hskip .5in \Psi_2 = {a_{11} + b_{22} - 4c_{12} \over 12} \hskip .5in \Psi_3 = {a_{12} - c_{22} \over 4} \hskip .5in \Psi_4 = {a_{22} \over 2},\eqno(2.25)$$
whence
$$\eqalignno{\Psi_{ABCD} &= {b_{11} \over 2}\beta_A\beta_B\beta_C\beta_D - (b_{12}-c_{11})\alpha_{(A}\beta_B\beta_C\beta_{D)} + {a_{11} + b_{22}-4c_{12} \over 2}\alpha_{(A}\alpha_B\beta_C\beta_{D)}\cr
&\qquad - (a_{12}-c_{22})\alpha_{(A}\alpha_B\alpha_C\beta_{D)} + {a_{22} \over 2}\alpha_A\alpha_B\alpha_C\alpha_D.&(2.26)\cr}$$
It is clear from the dyad components in (2.25) that, generically, there will be no relation between the $I$ and $J$ of [19](24), which therefore impose no constraints on $\Psi_{ABCD}$ for the general Walker metric.

We note that the $\Psi$-ON frame (A1.18) determines, as in (2.12), a null tetrad $\{L^a,N^a,M^a,\tilde M^a\}$, which one computes to be
$$\vcenter{\openup2\jot \halign{$\hfil#$&&${}#\hfil$&\qquad$\hfil#$\cr
L^a &= \displaystyle{1 \over \sqrt 2}(-a\partial_u + 2\partial_x) = \beta^A\left(\sqrt 2\xi^{A'} + {c \over \sqrt 2}\pi^{A'}\right) & N^a &= \displaystyle{1 \over \sqrt 2}\ell^a = {1 \over \sqrt 2}\alpha^A\pi^{A'}\cr
\tilde M^a &= \displaystyle{1 \over \sqrt 2}\left(-2c\partial_u - b\partial_v + 2\partial_y\right) = -\alpha^A\left(\sqrt 2 \xi^{A'} + {c \over \sqrt 2}\pi^{A'}\right) & M^a &= \displaystyle-{1 \over \sqrt 2}\tilde m^a = -{1 \over \sqrt 2}\beta^A\pi^{A'}.\cr}}\eqno(2.27)$$
The associated spin frames are (up to an overall sign) $\{\beta^A,-\alpha^A\}$ and $\{\sqrt 2\xi^{A'} + (c/\sqrt 2)\pi^{A'},-\pi^{A'}/\sqrt 2\}$. One can therefore obtain the components of the Weyl spinors with respect to these spin frames either in the manner followed above or simply by re-expressing (2.21) and (2.26) in terms of the relevant spin frames. Putting
$$o^A := \beta^A \hskip .5in \iota^A := -\alpha^A \hskip .5in o^{A'} := \sqrt 2\xi^{A'} + {c \over \sqrt 2}\pi^{A'} \hskip .5in \iota^{A'} := -{\pi^{A'} \over \sqrt 2},\eqno(2.28)$$
one finds
$$\tilde\Psi_{A'B'C'D'} = -{(A+S) \over 6}\iota_{A'}\iota_{B'}\iota_{C'}\iota_{D'} - Bo_{(A'}\iota_{B'}\iota_{C'}\iota_{D')} + {S \over 2}o_{(A'}o_{B'}\iota_{C'}\iota_{D')},\eqno(2.29)$$
the analysis of which is similar to that of (2.21), and
$$\eqalignno{\Psi_{ABCD} &= {a_{22} \over 2}\iota_A\iota_B\iota_C\iota_D - (c_{22} - a_{12})o_{(A}\iota_B\iota_C\iota_{D)} + {a_{11} + b_{22}-4c_{12} \over 2}o_{(A}o_B\iota_C\iota_{D)}\cr
&\qquad - (c_{11} - b_{12})o_{(A}o_Bo_C\iota_{D)} + {b_{11} \over 2}o_Ao_Bo_Co_D.&(2.30)\cr}$$
Returning to the Walker null tetrad, consider the Ricci spinor:
$$\Phi_{ABA'B'} = {1 \over 2}\left(R_{ab} - {S \over 4}g_{ab}\right) = {1 \over 2}E_{ab}.$$
From A1.7--8, one can compute the components of $E_{ab}$ with respect to the null tetrad (2.11) and thence the dyad components of $\Phi_{ABA'B'}$. (Note that by 2.5 $\Phi_{ABA'B'}\pi^{A'}\pi^{B'} = 0$, i.e., $\Phi_{AB0'0'} = 0$.) One finds:
$$\vcenter{\openup2\jot \halign{$\hfil#$&&${}#\hfil$&\qquad$\hfil#$\cr
\Phi_{00} &= \displaystyle\Phi_{000'0'} = {1 \over 2}\ell_aE^a{}_b\ell^b = 0; & \Phi_{01} &= \displaystyle\Phi_{000'1'} = {1 \over 2}\ell_aE^a{}_bm^b = -{\mu \over 2};\cr
\Phi_{02} &= \displaystyle\Phi_{001'1'} = {1 \over 2}m_aE^a{}_bm^b = {\Upsilon \over 2}; & \Phi_{10} &= \displaystyle\Phi_{010'0'} = {1 \over 2}\ell_aE^a{}_b\tilde m^b = 0;\cr
\Phi_{11} &= \displaystyle\Phi_{010'1'} = {1 \over 2}\ell_aE^a{}_bn^b = {\theta \over 2}; & \Phi_{20} &= \displaystyle\Phi_{110'0'} = {1 \over 2}\tilde m_aE^a{}_b\tilde m^b = 0;\cr
\Phi_{21} &= \displaystyle\Phi_{110'1'} + {1 \over 2}\tilde m_aE^a{}_bn^b = {\nu \over 2}; & \Phi_{22} &= \displaystyle\Phi_{111'1'} = {1 \over 2}n_aE^a{}_bn^b = {\zeta \over 2};\cr}}\eqno(2.31)$$
$$\Phi_{12} = \Phi_{011'1'} = {1 \over 2}m_aE^a{}_bn^b = -{1 \over 4}(2\eta + 2c\theta + b\nu - a\mu).$$
It follows, see also 2.5, that
$$\Phi_{ABA'B'} = A_{AB}\pi_{A'}\pi_{B'} + B_{AB}\pi_{(A'}\xi_{B')},\eqno(2.32)$$
with
$$A_{AB} = {1 \over 2}\left(\Upsilon\beta_A\beta_B + (2\eta + 2c\theta + b\nu - a\mu)\alpha_{(A}\beta_{B)} + \zeta\alpha_A\alpha_B\right)\qquad
B_{AB} = \mu\beta_A\beta_B + 2\theta\alpha_{(A}\beta_{B)} - \nu\alpha_A\alpha_B.\eqno(2.33)$$
With respect to the basis induced by $\{\alpha^A,\beta^B\}$ for $S \otimes S$, $A_{AB}$ and $B_{AB}$ have components
$$A_{\bf AB}  = {1 \over 2}\pmatrix{\Upsilon&-(2c\theta + 2\eta + b\nu - a\mu)/2\cr -(2c\theta + 2\eta + b\nu - a\mu)/2&\zeta\cr} \hskip .75in B_{\bf AB} =\pmatrix{\mu&-\theta\cr -\theta&-\nu\cr}.$$
\vskip 24pt
\noindent {\section 3. Walker Geometry with Parallel LSRs}
\vskip 12pt
Dunajski [11] considered a four manifold $M$ with a neutral metric $g$ and a global parallel spinor $\pi^{A'}$ and called such a {\sl null-K\"ahler\/} four-manifold. Since a global parallel spinor $\pi^{A'}$ satisfies (2.2), $[\pi^{A'}]$ defines a parallel distribution of $\alpha$-planes and hence a Walker geometry. In fact, the considerations in [11] are essentially local, and his main result is naturally subsumed as a feature of Walker geometry, as we explain in this section. 

It is natural to study Walker geometries $(M,g,{\cal D},[\pi^{A'}])$ which admit parallel LSRs. We first note the following restrictions imposed on the curvature by the presence of a parallel spinor ([11]). Note that these restrictions are local, rather than global, in nature, in the sense that they hold on the domain of any parallel spinor.
\vskip 24pt
\noindent {\bf 3.1 Proposition}\hfil\break
{\it If $(M,g)$ is a neutral four-manifold and $\pi^{A'}$ a parallel spinor on some domain (with nonempty interior) in $M$, then on that domain: $g$ is Ricci-scalar flat; $[\pi_{A'}]$ is a principal spinor of the Ricci spinor of multiplicity two, i.e., $\Phi_{ABA'B'} = F_{AB}\pi_{A'}\pi_{B'}$ for some symmetric $F_{AB}$; and $[\pi_{A'}]$ is a WPS of multiplicity four, i.e., if nonzero, $\tilde\Psi_{A'B'C'D'}$ is of type $\{4\}${\rm II} and case {\rm (iii)} of {\rm 2.6} pertains.}

{\it Proof.} Work locally on the interior of the domain of $\pi^{A'}$; the curvature conditions then extend to the closure of that domain by continuity. The proof follows the argument of 2.5 but beginning with $\nabla_b\pi_{A'} = 0$ rather than (2.2), which eliminates a factor of $\pi^{A'}$ from the computations and therefore increases the multiplicity of $[\pi^{A'}]$ as a principal spinor of both $\Phi_{ABA'B'}$ and $\tilde\Psi_{A'B'C'D'}$ by one. From (A2.9), $\Square_{A'B'}\pi^{A'} = -3\Lambda\pi_{A'}$, whence $\pi^{A'}$ parallel entails $\Lambda = 0$, which fact entails $[\pi^{A'}]$ is a WPS of multiplicity four.

Thus, $\tilde\Psi_{A'B'C'D'} = c\pi_{A'}\pi_{B'}\pi_{C'}\pi_{D'}$. Substituting the expressions obtained for the Weyl and Ricci curvature spinors into (A2.10) yields, since $\pi_{A'}$ is parallel,
$$\pi_{A'}\pi_{B'}\pi_{C'}\pi_{D'}\nabla^{A'}_Bc = \pi_{(C'}\pi_{D'}\nabla_{B')}^AF_{AB} \hskip 1.25in \pi_{A'}\pi_{B'}\nabla^{CA'}F_{CD} = 0.$$
In fact, transvecting the first, by say $\pi^{B'}$, yields the second. Transvecting the first by $\eta^{B'}\eta^{C'}\eta^{D'}$, where $\eta^{D'}\pi_{D'} = 1$, yields $\pi_{A'}\nabla^{A'}_Bc = \eta^{B'}\nabla^A_{B'}F_{AB}$. Transvecting the second equation above by $\eta^{B'}$ yields \hfil\break
$\pi_{A'}\nabla^{A'C}F_{CD} = 0$. It follows that
$$\pi_{A'}\nabla^{A'}_Bc = 0\ \Leftrightarrow\ \nabla^C_{A'}F_{CD} = 0,\eqno(3.1)$$
\vskip 24pt
\noindent {\bf 3.2 Proposition}\hfil\break
{\it A Walker geometry $(M,g,{\cal D},[\pi^{A'}])$ admits, on a neighbourhood of each point, a parallel LSR $\pi^{A'}$ iff the curvature conditions of {\rm 3.1} pertain on $M$.}

{\it Proof.} The necessity follows immediately from 3.1. For sufficiency, choose any LSR $\pi^{A'}$ on a neighbourhood of some point and consider $f\pi^{A'}$, where $f$ is a smooth positive function.

Observe that $0 = \nabla_b(f\pi^{A'}) = (\nabla_bf)\pi^{A'} + f\nabla_b(\pi^{A'})$ is equivalent to $(\nabla_b f)\pi^{A'} = -f\nabla_b\pi^{A'} = -fP_b\pi^{A'}$, by (2.3), i.e., equivalent to $\nabla_b\bigl(\ln(f)\bigr) = -P_b$. One can therefore find an $f$ to rescale $\pi^{A'}$ to be parallel iff $-P_b$ is a gradient. By the Poincar\'e lemma, this is so iff $P_b$ is closed, at least locally, i.e., iff $\nabla_{[c}P_{b]} = 0$. Now
$$\eqalign{\nabla_c\nabla_b\pi^{A'} &= (\nabla_cP_b)\pi^{A'} + P_b\nabla_c\pi^{A'}\cr
&= \left(\nabla_cP_b + P_bP_c\right)\pi^{A'}.\cr}$$
Hence, $\nabla_{[c}\nabla_{b]}\pi^{A'} = \nabla_{[c}P_{b]}\pi^{A'}$, i.e., $P_b$ is closed iff $\nabla_{[c}\nabla_{b]}\pi^{A'} = 0$. This condition is not true in general, so one cannot always find an LSR of $[\pi^{A'}]$ that is parallel. From (A2.7) and (A2.9), however, the curvature conditions of 3.1 do entail $\nabla_{[c}\nabla_{b]}\pi^{A'} = 0$.
\vskip 24pt
\noindent {\bf 3.3 Lemma}\hfil\break
{\it For any Walker geometry $(M,g,{\cal D},[\pi^{A'}])$, suppose given a parallel LSR $\pi^{A'}$ of $[\pi^{A'}]$ on a neighbourhood of some point $m \in M$. Then one can choose Walker coordinates $(u,v,x,y)$ on a (smaller) neighbourhood of $m$ for which {\rm (2.8)} is valid with this parallel LSR.}

{\it Proof.} Put $L_a := \mu_A\pi_{A'}$ and $\tilde M_a := \nu_A\pi_{A'}$, where $\{\mu^A,\nu^A\}$ is a spin frame, i.e., $\nu^D\mu_D = 1$. Then $L_a \wedge \tilde M_b = \epsilon_{AB}\pi_{A'}\pi_{B'} =: \sigma_\pi$ is parallel and therefore a closed two-form. Thus, locally, $\sigma_\pi = d\phi$, for some one-form $\phi$. Now, $\sigma_\pi$ is of rank one, i.e., $\sigma_\pi \wedge \sigma_\pi = 0$. Since $\phi$ is defined only up to the addition of an exact one-form, one can exploit this freedom to ensure that $\phi \wedge d\phi$ is nonvanishing. It follows from Darboux's theorem (e.g., Theorem 6.2 in [38]) that there exist local coordinates $(x,p,q,y)$ such that $\phi = xdy + dp$, whence $\sigma_\pi = d\phi = dx \wedge dy$.

Now ${\cal D} = \langle L^a,\tilde M^a \rangle_{\bf R}$ and is the kernel of the endomorphism of the tangent space: $v \mapsto v \hook \sigma_\pi$. In particular,
$$0 = L^a \hook \sigma_\pi = L^a \hook (dx \wedge dy) = dx(L)dy - dy(L)dx.$$
Since $dx$ and $dy$ are linearly independent one-forms,then $dx(L) = dy(L) = 0$, whence $L \in \langle \partial/\partial p,\partial/\partial q \rangle_{\bf R}$. Together with a similar computation for $\tilde M^a$, one deduces that 
$$\bigl(\ker(dx) \cap \ker(dy)\bigr) = \langle \partial/\partial p,\partial/\partial q \rangle_{\bf R} = \langle L^a,\tilde M^a \rangle_{\bf R} = {\cal D}.$$
Thus, $(p,q,x,y)$ are Frobenius coordinates for $\cal D$. It follows, as in 2.3, that $dx = \alpha_A\pi_{A'}$, $dy = \beta_A\pi_{A'}$, for some spinors $\alpha_A$ and $\beta_A$ satisfying $\beta^D\alpha_D \not= 0$. But $\sigma_\pi = dx \wedge dy = (\beta^D\alpha_D)\epsilon_{AB}\pi_{A'}\pi_{B'}$, whence $\beta^D\alpha_D = 1$. Putting $U^a = \alpha^A\pi^{A'}$ and $V^a = \beta^A\pi^{A'}$, Walker's argument 2.3 constructs functions $u$ and $v$ so that $(u,v,x,y)$ are Frobenius coordinates, with $\partial_u = \alpha^A\pi^{A'}$ and $\partial_v = \beta^A\pi^{A'}$, which is the desired result. 
\vskip 24pt
\noindent {\bf 3.4 Lemma}\hfil\break
{\it Given any Walker coordinates $(u,v,x,y)$, the LSR $\pi^{A'}$ defined by {\rm (2.7)} is parallel iff $dx \wedge dy$, equivalently, $\partial_u \wedge \partial_v$, is parallel.}

{\it Proof.} All that needs to be verified here is that $\pi^{A'}\pi^{B'}$ parallel entails $\pi^{A'}$ is parallel. Supposing $\pi^{A'}\pi^{B'}$ is parallel, then by (2.3) $0 = \nabla_c(\pi^{A'}\pi^{B'}) = 2P_c\pi^{A'}\pi^{B'}$, which entails that $P_c$ is zero, i.e., $\pi^{A'}$ is indeed parallel.

In a slightly different formulation, for any LSR $\pi^{A'}$, by (2.6) $dx \wedge dy = (\nu^D\mu_D)\sigma_\pi =: \pm f\sigma_\pi$, with $f$ positive. Hence, $\nabla_b(dx \wedge dy) = \pm(\nabla_bf + 2fP_b)\Sigma_\pi$. Thus, $dx \wedge dy$ is parallel iff $-P_b = \nabla_b(\ln\sqrt{f})$. From the proof of 3.2, one notes that this last equation is precisely the condition for $\sqrt{f}\pi^{A'}$ to be parallel; and $dx \wedge dy = \pm\Sigma_{\sqrt{f}\pi}$.
\vskip 24pt
\noindent {\bf 3.5 Lemma}\hfil\break
{\it The LSR $\pi^{A'}$ determined up to sign by the condition {\rm (2.7)} is parallel iff}
$$a_1 + c_2 = 0 \hskip 1.5in b_2 + c_1 = 0.\eqno(3.2)$$

{\it Proof.} From (A1.8), a direct computation yields
$$\displaylines{\nabla_1(\partial_u \wedge \partial_v) = 0 \hskip 1in \nabla_2(\partial_u \wedge \partial_v) = 0\cr
\cr
\nabla_3(\partial_u \wedge \partial_v) = {a_1 + c_2 \over 2}\partial_u \wedge \partial_v \hskip 1in
\nabla_4(\partial_u \wedge \partial_v) = {b_2 + c_1 \over 2}\partial_u \wedge \partial_v.\cr}$$
\vskip 24pt
\noindent {\bf 3.6 Remarks}\hfil\break
The conditions (3.2) imply
$$a_{12} + c_{22} = 0 = b_{12} + c_{11} \hskip 1in a_{11} = b_{22} = -c_{12},\eqno(3.3)$$
which entail simplifications of the curvature. Obviously one expects to recover the results of 3.1. Indeed, from A1.6 one sees immediately that $S = 0$ and from A1.5 that $R_{\bf ij} = 0$ unless {\bf i}, {\bf j} are 3 or 4. More geometrically, in A1.8 $\theta = \mu = \nu = 0$, whence the Einstein endomorphism (here equal to the Ricci endomorphism since $S = 0$) maps the tangent space to $\cal D$ with kernel $\cal D$. Since the contraction of both $\partial_u$ and $\partial_v$ (which have linearly independent unprimed spinor parts) on the Ricci tensor is zero, then (since $S = 0$) one deduces that the Ricci spinor $\Phi_{ABA'B'}$ is null of the form $F_{AB}\pi_{A'}\pi_{B'}$.

By (3.2), the $B$ of (A1.24) vanishes. Since $S=0$, assuming ${\cal W}^+ \not= 0$, it follows that case (iii) of 2.6 pertains, i.e., $\pi^{A'}$ is WPS of multiplicity four and $\tilde\Psi_{A'B'C'D'}$ is null. Substituting $S= B= 0$ in (2.21) yields
$$\tilde\Psi_{A'B'C'D'} = -{A \over 24}\pi_{A'}\pi_{B'}\pi_{C'}\pi_{D'}.\eqno(3.4)$$ 

(3.3) entails a simplification in the formulae (2.25--26) for the ASD Weyl curvature:
$$\Psi_0 = {b_{11} \over 2} \hskip .3in \Psi_1 = {b_{12} \over 2} = -{c_{11} \over 2} \hskip .3in \Psi_2 = {a_{11} \over 2} = {b_{22} \over 2} = -{c_{12} \over 2} \hskip .3in \Psi_3 = {a_{12} \over 2} = -{c_{22} \over 2} \hskip .3in \Psi_4 = {a_{22} \over 2}.\eqno(3.5)$$
Hence, in the present circumstances, the ASD Weyl curvature is zero iff $a$, $b$, and $c$ are affine functions of $u$ and $v$, whose coefficients are functions of $x$ and $y$. Imposing (3.2) then yields
$$a = Eu + Fv + G \hskip .5in b = Mu + Nv + P \hskip .5in c = -Nu - Ev + T,\eqno(3.6)$$
i.e., for a Walker geometry $(M,g,{\cal D},[\pi^{A'}])$, if some LSR $\pi^{A'}$ is parallel and $(u,v,x,y)$ are Walker coordinates for which (2.7) holds, then the ASD curvature vanishes iff $a$, $b$, and $c$ are of the form as in (3.6), leaving 7 arbitrary functions of $x$ and $y$.
\vskip 24pt
\noindent {\bf 3.7 Observation}\hfil\break
The equations (3.2) can be simply solved as follows. Take $c$ to be any smooth function of $(u,v,x,y)$ such that $c_v$ is integrable wrt $u$ and $c_u$ is integrable wrt $v$; then one directly solves for $a$ and $b$.

Presume in fact that $c$ is at least $C^4$. Assuming (3.2) holds, then $a$ is an antiderivative of $-c_v$ wrt $u$, and $b$ is an antiderivative of $-c_u$ wrt $v$. Let $g :=: \int c\,du$ denote a specific antiderivative of $c$ wrt $u$ and $h :=: \int c\,dv$ a specific antiderivative of $c$ with respect to $v$. Since $g_u = h_v$ there exists a function $k :=:\int\int c\,dudv$ such that $k_u = h$ and $k_v = g$. Put
$$\vartheta := {1 \over 2}\int\int c\, dudv + m(u,x,y)+ n(v,x,y),\eqno(3.7)$$
where $m$ and $n$ are arbitrary (suitably) smooth functions. Then, $\vartheta_{uv} = c/2$, $2\vartheta_{vv} = \int c_v\,du + 2n_{vv} = -\int a_u\,du + 2n_{vv}$, and $2\vartheta_{uu} = \int c_u\,dv + 2m_{uu} = -\int b_v + 2m_{uu}$, where the remaining integrals are specific antiderivatives of their integrands. Use the freedom in the choice of $m$ and $n$ to ensure that in fact
$$W = \pmatrix{a&c\cr c&b\cr} = 2\pmatrix{-\vartheta_{vv}&\vartheta_{uv}\cr \vartheta_{uv}&-\vartheta_{uu}\cr} = 2\pmatrix{-\vartheta_{22}&\vartheta_{12}\cr \vartheta_{12}&-\vartheta_{11}\cr}.\eqno(3.8)$$
By (3.4), the vanishing of $\tilde\Psi_{A'B'C'D'}$ is equivalent to $A=0$. It is therefore natural to re-express $A$ in terms of $\vartheta$. With the advantage of hindsight obtained from [11] and the results of \S 4, we wrote (A1.23) in the form we require here.
Observing that under (3.2), in (A1.23) each of the first and second pairs of terms in the final line cancel, while the last three terms together equal $-S$, which is zero under (3.2), one may ignore the final line.

Using A1.1, for an arbitrary function $f$:
$$\Square f = g^{ab}\nabla_a\nabla_bf = g^{\bf ij}\partial_{\bf i}\partial_{\bf j}f - 2\left(\Gamma^{\bf k}_{13} + \Gamma^{\bf k}_{24}\right)\partial_{\bf k}f = g^{\bf ij}\partial_{\bf i}\partial_{\bf j}f - (a_1 + c_2)\partial_1f - (b_2 + c_1)\partial_2f.\eqno(3.9)$$
Thus, when (3.2) holds,
$$\Square f = g^{\bf ij}\partial_{\bf i}\partial_{\bf j}f = -af_{11} -2cf_{12} - bf_{22} + 2f_{13} + 2f_{24}.\eqno(3.10)$$
Using (3.2) to rewrite certain terms, it is now routine to verify that: the first line of (A1.23) equals $12\Square(\vartheta_{13})$; the second line equals $12\Square(\vartheta_{24})$, the third and fourth lines together equal $3\Square(ab) = 12 \Square(\vartheta_{11}\vartheta_{22})$ and the fifth line equals $-3\Square(c^2) = -12\Square\bigl((\vartheta_{12})^2\bigr)$. Hence,
$$A = 12\Square\bigl(\vartheta_{13} + \vartheta_{24} + \vartheta_{11}\vartheta_{22} - (\vartheta_{12})^2\bigr).\eqno(3.11)$$
Thus, when (3.2) holds, one may write the metric in the form (1.1) with $W$ determined by a function $\vartheta$ as in (3.8), the curvature conditions of 3.1 hold and the vanishing of $\tilde \Psi_{A'B'C'D'}$ is equivalent to the vanishing of (3.11). This result refines that of Dunajski [11], who studied ASD, four-dimensional, neutral metrics admitting a parallel spinor and also wrote the condition for the vanishing of the SD Weyl curvature in the form:
$$\vartheta_{13} + \vartheta_{24} + \vartheta_{11}\vartheta_{22} - (\vartheta_{12})^2 =: P\hskip 1in \Square P = 0.\eqno(3.12)$$
The homogeneous form of this second order PDE is Pleba\~nski's [30] {\sl second heavenly equation\/} characterizing, locally, the metrics of (in the real case) neutral four-manifolds whose only nontrivial curvature is the ASD Weyl curvature. We elucidate this connexion in \S 4. Of course, any Walker metric with $W$ of the form (3.8) trivially satisfies (3.2).

Now consider the form that the ASD Weyl curvature (2.26) and (3.5) takes under (3.2). Substituting (3.8) into (2.25) or (3.5) yields:
$$\Psi_0 = -\vartheta_{1111} \hskip .5in \Psi_1 = -\vartheta_{1112} \hskip .5in \Psi_2 = -\vartheta_{1122} \hskip .5in \Psi_3 = -\vartheta_{1222} \hskip .5in \Psi_4 = -\vartheta_{2222}.\eqno(3.13)$$
With
$$\delta_A := \alpha_A\partial_v - \beta_A\partial_u = \alpha_A\tilde m^b\nabla_b - \beta_A\ell^b\nabla_b = \pi^{A'}\nabla_{AA'},\eqno(3.14)$$
one finds that
$$\eqalignno{\Psi_{ABCD} &= -\vartheta_{1111}\beta_A\beta_B\beta_C\beta_D + 4\vartheta_{1112}\alpha_{(A}\beta_B\beta_C\beta_{D)} - 6\vartheta_{1122}\alpha_{(A}\alpha_B\beta_C\beta_{D)}\cr
&\qquad + 4\vartheta_{1222}\alpha_{(A}\alpha_B\alpha_C\beta_{D)} - \vartheta_{2222}\alpha_A\alpha_B\alpha_C\alpha_D&(3.15)\cr
& = -\delta_A\delta_B\delta_C\delta_D\vartheta.\cr}$$
(3.15) is the form of the ASD Weyl spinor obtained by Pleba\~nski [30] under the assumption that all the other curvature vanishes, and by Dunajski [11] assuming 3.1 and the vanishing of the SD Weyl curvature. We have found this form is valid wherever a Walker metric has a parallel LSR, i.e., under the constraints imposed on the curvature in 3.1; note, in particular, that the vanishing of $\tilde\Psi_{A'B'C'D'}$ is not required. In short, Walker geometry provides the natural context for this generalization of Dunajski's result.

In accordance with 3.1, under (3.2) $\theta = \mu = \nu = 0$ in A1.8, whence $B_{AB} = 0$ in (2.33). Hence,
$$\Phi_{ABA'B'} = A_{AB}\pi_{A'}\pi_{B'}\qquad\hbox{where}\qquad A_{AB} = {1 \over 2}\left(\Upsilon\beta_A\beta_B + 2\eta\alpha_{(A}\beta_{B)} + \zeta\alpha_A\alpha_B\right).\eqno(3.16)$$
Under (3.2), one further finds that $\Upsilon/2 = P_{11}$, $\eta/2 = -P_{12}$, and $\zeta/2 = P_{22}$, whence one can express (3.16) in the form
$$\Phi_{ABA'B'} = \left(\delta_A\delta_B P\right)\pi_{A'}\pi_{B'}.\eqno(3.17)$$
\noindent {\bf 3.8 Theorem}\hfil\break
{\it In summary, a Walker geometry $(M,g,[\pi^{A'}])$ with a parallel LSR on some coordinate domain, is determined on that domain by a single function $\vartheta$ of Walker coordinates with: the metric given by {\rm (3.8)}; $\tilde\Psi_{A'B'C'D'}$ of type $\{4\}$ with WPS $[\pi^{A'}]$ or vanishes when (3.12) holds; $\Psi_{ABCD} = -\delta_D\delta_C\delta_B\delta_A\vartheta$ and thus vanishes iff $\vartheta$ is a cubic polynomial in $u$ and $v$ of the form
$$-12\vartheta = Mu^3 + 3Euv^2 + 3Nu^2v + Fv^3 + 3Gv^2 -6Tuv + 3Pu^2 + K_1u + K_2v + K_3,$$
consistent with {\rm (3.6)}, where the coefficients are arbitrary functions of $(x,y)$ and $K_1$, $K_2$ and $K_3$ express the residual freedom in $\vartheta$ not constrained by {\rm (3.8)}; $S=0$; the Ricci spinor is given by {\rm (3.17)}, whence the Walker geometry is Einstein iff $P$ is affine in $u$ and $v$.}
\vskip 24pt
\noindent{\section 4. Complementary Distributions}
\vskip 12pt
We begin this section by considering an orientable four-dimensional neutral manifold $(M,g)$ which admits a pair of complementary distributions of totally null two-planes. By 2.1, each distribution is a distribution of either $\alpha$-planes or of $\beta$-planes. Since at any point, any $\alpha$-plane intersects any $\beta$-plane in a one-dimensional subspace, a complementary pair of totally null distributions must both be distributions of $\alpha$-planes or distributions of $\beta$-planes. If one of the distributions is parallel, then one has a Walker geometry, and with respect to the canonical orientation the two distributions are each distributions of $\alpha$-planes.

Hence, we consider two complementary distributions of $\alpha$-planes on $(M,g)$. Denote these distributions by ${\cal D}^\pi$ and ${\cal D}^\chi$, where the distributions have associated projective spinor fields $[\pi^{A'}]$ and $[\chi^{A'}]$ respectively. The complementarity is equivalent to $\chi^{D'}\pi_{D'} \not= 0$ for any LSRs.

We now suppose both distributions are integrable. In a suitable neighbourhood of any point $m \in M$, one can choose Frobenius coordinates $(p,q,x,y)$ for ${\cal D}^\pi$ with ${\cal D}^\pi = \langle \partial_p,\partial_q \rangle_{\bf R}$ and Frobenius coordinates $(w,z,r,s)$ for ${\cal D}^\chi$ with ${\cal D}^\chi = \langle \partial_w,\partial_z \rangle_{\bf R}$. Since $dx$ and $dy$ separately annihilate ${\cal D}^\pi$, for any LSR $\pi_{A'}$ of $[\pi_{A'}]$ they are of the form $\gamma_A\pi_{A'}$, for suitable $\gamma_A$. One can rescale $\pi_{A'}$, and interchange the roles of $x$ and $y$ if necessary, so that one may suppose the coordinates $(p,q,x,y)$ and LSR $\pi_{A'}$ are chosen so that
$$dx = \alpha_A\pi_{A'} \hskip .75in dy = \beta_A\pi_{A'} \hskip .75in \beta^D\alpha_D = 1,\qquad\hbox{whence}\qquad dx \wedge dy = \sigma_\pi.\eqno(4.1)$$
Similarly, one may suppose the coordinates $(w,z,r,s)$ and an LSR $\chi_{A'}$ are chosen so that
$$dr = \mu_A\chi_{A'} \hskip .75in ds = \nu_A\chi_{A'} \hskip .75in \nu^D\mu_D = 1,\qquad\hbox{whence}\qquad dr \wedge ds = \sigma_\chi,\eqno(4.2)$$
where $\sigma_\chi = \epsilon_{AB}\chi_{A'}\chi_{B'}$. Since the distributions are complementary, the functions $(r,s,x,y)$ constitute local coordinates which are simultaneously Frobenius with respect to both distributions, specifically
$${\cal D}^\pi = \langle \partial_r,\partial_s \rangle_{\bf R} \hskip 1.25in {\cal D}^\chi = \langle \partial_x,\partial_y \rangle_{\bf R},\eqno(4.3)$$
(see for example [16], p. 182) and (4.1) and (4.2) are of course still valid, being coordinate-independent statements. With respect to the coordinates $(r,s,x,y)$, the metric must take the form
$$(g_{\bf ab}) = \pmatrix{{\bf 0}_2&D\cr {^\tau\! D}&{\bf 0}_2\cr},\qquad\hbox{whence}\qquad \left(g^{\bf ab}\right) = \pmatrix{{\bf 0}_2&{^\tau\! D}^{-1}\cr D^{-1}&{\bf 0}_2\cr}\eqno(4.4)$$
for some $(2 \times 2)$-matrix $D$. Note that $\det(g_{\bf ab}) = \bigl(\det(D)\bigr)^2$. Writing the metric as
$$2dx(D_{11}dr + D_{21}ds) - 2dy(-D_{12}dr - D_{22}ds),\eqno(4.5)$$
one can extract the null tetrad:
$$\vcenter{\openup2\jot \halign{$\hfil#$&&${}#\hfil$&\qquad$\hfil#$\cr
\ell_a &:= dx \qquad& n_a &:= D_{11}dr + D_{21}ds\cr
\tilde m_a &:= dy \qquad& m_a &:= -(D_{12}dr + D_{22}ds).\cr}}\eqno(4.6)$$
One computes
$$\vcenter{\openup2\jot \halign{$\hfil#$&&${}#\hfil$&\qquad$\hfil#$\cr
\ell^a &= (D^{-1})_{11}\partial_r + (D^{-1})_{12}\partial_s & n^a &= \partial_x\cr
\tilde m^a &= (D^{-1})_{21}\partial_r + (D^{-1})_{22}\partial_v & m^a &= -\partial_y.\cr}}\eqno(4.7)$$
From (4.1--2) and (4.4), one computes
$$\vcenter{\openup2\jot \halign{$\hfil#$&&${}#\hfil$&\qquad$\hfil#$\cr
(D^{-1})_{11} &:= g^{ab}(dr,dx) = (\mu^D\alpha_D)(\chi^{D'}\pi_{D'}) & (D^{-1})_{21} &:= g^{ab}(dr,dy) = (\mu^D\beta_D)(\chi^{D'}\pi_{D'});\cr
(D^{-1})_{12} &:= g^{ab}(ds,dx) = (\nu^D\alpha_D)(\chi^{D'}\pi_{D'}) & (D^{-1})_{22} &:= g^{ab}(ds,dy) = (\nu^D\beta_D)(\chi^{D'}\pi_{D'})\cr}}\eqno(4.8)$$
whence
$$\det(D) = \bigl(\det(D^{-1})\bigr)^{-1} = \left[(\mu^D\alpha_D)(\nu^E\beta_E) - (\mu^D\beta_D)(\nu^E\alpha_E)\right]^{-1}(\chi^{D'}\pi_{D'})^{-2}.$$
But, $\{\alpha^A,\beta^A\}$ and $\{\mu^A,\nu^A\}$ are both spin frames, and therefore related by an element of {\bf SL(2;R)}; indeed,
$$\mu^A = -(\mu^D\beta_D)\alpha^A + (\mu^D\alpha_D)\beta^A \hskip 1in \nu^A = -(\nu^D\beta_D)\alpha^A + (\nu^D\alpha_D)\beta^A,\eqno(4.9)$$
whence $1 = -(\mu^D\beta_D)(\nu^E\alpha_E) + (\nu^D\beta_D)(\mu^E\alpha_E)$, and therefore
$$\det(D) = (\chi^{D'}\pi_{D'})^{-2} > 0.\eqno(4.10)$$
Putting $\xi^{A'} := \chi^{A'}/(\chi^{D'}\pi_{D'})$ so that $\xi^{D'}\pi_{D'} = 1$, then, from (4.6),
$$n_a = (\chi^{D'}\pi_{D'})(D_{11}\mu_A + D_{21}\nu_A)\xi_{A'}\hskip .75in m_a = -(\chi \cdot \pi)(D_{12}\mu_A + D_{22}\nu_A)\xi_{A'}$$
With
$$\kappa^A := -(\chi^{D'}\pi_{D'})(D_{12}\mu_A + D_{22}\nu_A) \hskip 1.25in \lambda^A = (\chi^{D'}\pi_{D'})(D_{11}\mu_A + D_{21}\nu_A),$$
the matrix 
$${\bf F} = (\chi\cdot\pi)\pmatrix{-D_{12}&D_{11}\cr -D_{22}&D_{21}\cr}$$
has determinant $(\chi^{D'}\pi_{D'})^2\det(D) = 1$, by (4.10). Thus, {\bf F} is an element of {\bf SL(2;R)} and $\{\kappa^A,\lambda^A\}$ is a spin frame. From (4.8--10), one checks that {\bf F} is inverse to the element of {\bf SL(2;R)} in (4.9) expressing $\{\mu^A,\nu^A\}$ in terms of $\{\alpha^A,\beta^A\}$, whence $\kappa^A = \alpha^A$ and $\lambda^A = \beta^A$. 

In summary, the null tetrad (4.6) is
$$\vcenter{\openup2\jot \halign{$\hfil#$&&${}#\hfil$&\qquad$\hfil#$\cr
\ell_a &= \alpha_A\pi_{A'} = dx & \tilde m_a &= \beta_A\pi_{A'} = dy\cr
m_a &= \alpha_A\xi_{A'} = -(D_{12}dr + D_{22}ds) & n_a &= \beta_A\xi_{A'} = D_{11}dr + D_{21}ds.\cr}}\eqno(4.11)$$
The complementary distributions are equivalent to an almost product structure $P$, an endomorphism of $TM$ satisfying $P^2 = 1$. Here, $P = 1_2 \oplus -1_2$ with respect to the decomposition $TM = {\cal D}^\pi \oplus {\cal D}^\chi$. It is straightforward to check that $P$ is in fact an anti-orthogonal automorphism of $(T_pM,g)$. By analogy with the K\"ahler form, define $\omega := g(P\ ,\ )$, which is indeed a two-form. In fact, one has:
$$P^2 = 1 \hskip .75in g_{ab}P^a{}_cP^b{}_d = -g_{cd} \hskip .75in \omega_{ab} := g_{cb}P^c{}_a = -\omega_{ba},\eqno(4.12)$$
any two of which entails the third. 

The integrability of the two complementary distributions is equivalent to the integrability of the almost product structure (i.e., the vanishing of the Nijenhuis tensor of $P$) and equivalent to the integrability of the induced ${\bf GL(2;R)} \times {\bf GL(2;R)}$-structure on the bundle of frames. It is well known that an integrable almost product structure (in any of these equivalent forms) is in turn equivalent to a locally product structure on $M$ (see, for example, [47], Ch. XI); indeed, an atlas of coordinates of the form $(r,s,x,y)$ as constructed above constitute an atlas for the locally product structure in the present circumstances; the Jacobian of the transformations between such coordinates systems must be of the form $\left({M \atop {\bf 0}}{{\bf 0} \atop N}\right)$ with $M$, $N \in {\bf GL(2;R)}$. The following result is also well known, e.g., [47], Ch. X.
\vskip 24pt
\noindent {\bf 4.1 Lemma}\hfil\break
{\it The almost product structure is parallel iff the two complementary distributions are each parallel.}
\vskip 24pt
So, now suppose that ${\cal D}^\pi$ and ${\cal D}^\chi$ are each parallel; in particular, $(M,g)$ is a Walker geometry, in two ways. Now $P$ is obviously parallel iff the associated two-form $\omega$ is parallel, whence $\omega$ is a symplectic form naturally associated with the complementary parallel totally null distributions. One computes
$$\omega = D_{11}dr \wedge dx + D_{12} dr \wedge dy + D_{21} ds \wedge dx + D_{22} ds \wedge dy.\eqno(4.13)$$
In the more usual theory of (almost) product spaces in which the almost product structure is an orthogonal automorphism of the tangent space with respect to a metric which restricts to nondegenerate metrics on the complementary distributions, the condition $\nabla P = 0$ is equivalent to $(M,g)$ being locally a Riemannian product ([47], Ch. X). In the present circumstances, we should therefore expect that $\nabla P = 0$ is equivalent to some condition on the metric.

To investigate, it will prove convenient to write the (locally product) coordinates $(r,s,x,y)$ derived above in the form $(x^{\bf a}) = (x^{\bf A},x^{\bf A'})$, where upper case indices range over 1 and 2 here (as in \S 1, context should prevent confusion of these concrete indices with concrete spinor indices). Then,
$$\displaylines{\nabla_{\bf a}\partial_r = \Gamma^{\bf c}_{{\bf a}1}\partial_{\bf c} \hskip 1in \nabla_{\bf a}\partial_s = \Gamma^{\bf c}_{{\bf a}2}\partial_{\bf c};\cr
\nabla_{\bf a}\partial_x = \Gamma^{\bf c}_{{\bf a}3}\partial_{\bf c} \hskip 1in \nabla_{\bf a}\partial_y = \Gamma^{\bf c}_{{\bf a}4}\partial_{\bf c};\cr}$$
from which it follows, from 1.1 and (4.3), that
$${\cal D}^\pi\hbox{ is parallel iff } \Gamma^{\bf C'}_{\bf aB} = 0\qquad\hbox{and}\qquad {\cal D}^\chi\hbox{ is parallel iff } \Gamma^{\bf C}_{\bf aB'} = 0\eqno(4.14)$$
(see also [47], Ch. X). From (4.4), one directly computes that:
$$\eqalign{\Gamma^{\bf C'}_{\bf aB} &= {1 \over 2}g^{\bf C'm}(g_{{\bf am},{\bf B}} + g_{{\bf Bm},{\bf a}} - g_{{\bf aB},{\bf m}})\cr
&= {1 \over 2}g^{\bf C'M}(g_{{\bf aM},{\bf B}} + g_{{\bf BM},{\bf a}} - g_{{\bf aB},{\bf M}})\cr
&= {1 \over 2}g^{\bf C'M}(g_{{\bf A'M},{\bf B}} - g_{{\bf A'B},{\bf M}})\cr
&= {1 \over 2}(D^{-1})^{\bf C'M}\bigl(({^\tau\! D})_{{\bf A'M},{\bf B}} - ({^\tau\! D})_{{\bf A'B},{\bf M}}\bigr)\cr} \hskip .75in \eqalign{\Gamma^{\bf C}_{\bf aB'} &= {1 \over 2}g^{\bf Cm}(g_{{\bf am},{\bf B'}} + g_{{\bf B'm},{\bf a}} - g_{{\bf aB'},{\bf m}})\cr
&= {1 \over 2}g^{\bf CM'}(g_{{\bf aM'},{\bf B'}} + g_{{\bf B'M'},{\bf a}} - g_{{\bf aB'},{\bf M'}})\cr
&= {1 \over 2}g^{\bf CM'}(g_{{\bf AM'},{\bf B'}} - g_{{\bf AB'},{\bf M'}})\cr
&= {1 \over 2}({^\tau\! D}^{-1})^{\bf CM'}\bigl(D_{{\bf AM'},{\bf B'}} - D_{{\bf AB'},{\bf M'}}\bigr)\cr}$$
\vskip 6pt
\noindent whence ${\cal D}^\pi$ is parallel iff $D_{{\bf MA'},{\bf B}} = D_{{\bf BA'},{\bf M}}$ and ${\cal D}^\chi$ is parallel iff $D_{{\bf AM'},{\bf B'}} = D_{{\bf AB'},{\bf M'}}$, i.e., 
$${\cal D}^\pi\hbox{ is parallel iff } {\partial D_{\bf MA'} \over \partial x^{\bf B}} = {\partial D_{\bf BA'} \over \partial x^{\bf M}}\qquad\hbox{and}\qquad{\cal D}^\chi\hbox{ is parallel iff } {\partial D_{\bf AM'} \over \partial x^{\bf B'}} = {\partial D_{\bf AB'} \over \partial x^{\bf M'}}.\eqno(4.15)$$
(The analogues of (4.15) for nondegenerate complementary distributions may be found in [47], Ch. X.) The nontrivial components of the second condition of (4.15) (${\bf B'} \not= {\bf M'}$) are the integrability conditions for functions $\phi_{\bf A}$ satisfying $\left(\partial\phi_{\bf A}/\partial x^{\bf B'}\right) = D_{\bf AB'}$. Substituting this expression into the first condition of (4.15), the nontrivial components (${\bf B} \not= {\bf M}$) are now the integrability conditions for a function $\Omega$ (context should preclude confusion with the Walker volume form; the reason for the choice of this symbol will become apparent shortly) satisfying $\Omega_r = \partial\phi_0/\partial x^{\bf A'}$ and $\Omega_s = \partial\phi_1/\partial x^{\bf A'}$, whence
$$D_{\bf AB'} = {\partial\phi_{\bf A} \over \partial x^{\bf B'}} = {\partial^2\Omega \over \partial x^{\bf A}\partial x^{\bf B'}}.\eqno(4.16)$$
Thus, the metric in locally product coordinates $(x^{\bf A},x^{\bf A'})$, see (4.4), is determined by a single function according to (4.16).

Note that spinors do not, in fact, play a necessary role in this derivation (so far in this section, spinors have played an essentially descriptive role) nor does the four dimensionality; stripping away all reference to spinors and the null tetrad (4.6), one obtains an argument valid for a neutral manifold of arbitrary dimension $2n$ with parallel complementary totally null distributions (now {\bf A} and ${\bf A'}$ each range over $1,\ldots,n$ and the null tetrad is replaced by a Witt frame):
\vskip 24pt
\noindent {\bf 4.2 Theorem}\hfil\break
{\it Let $(M,g)$ be a neutral manifold of dimension $2n$ admitting complementary totally null distributions. This structure is equivalent to an almost product structure $P$ which is an anti-orthogonal automorphism of each $(T_pM,g)$. Integrability of $P$ (vanishing of its Nijenhuis tensor) is equivalent to integrability of the two distributions, and coordinates $(x^{\bf A},x^{\bf A'})$ which are simultaneously Frobenius for both distributions provide a locally product structure for $(M,g)$. When $P$ is parallel, equivalently the two distributions are parallel, with respect to the local product coordinates the metric takes the form
$$(g_{\bf ab}) = \pmatrix{{\bf 0}_n&D\cr {^\tau\! D}&{\bf 0}_n\cr}$$
with $D$ as in {\rm (4.16)}, and the associated two-form $\omega$ is a symplectic form with coordinate expression 
$$\omega = {\partial^2\Omega \over \partial x^{\bf A}\partial x^{\bf B'}}dx^{\bf A} \wedge dx^{\bf B'}.$$
Note that the canonical orientation induced by the one distribution is $(-1)^n$ times that of the other.}
\vskip 24pt
We will refer to this geometry as {\sl double Walker\/} geometry. It is already known, however, under other names. As an instance of locally product geometry, there is an obvious strong analogy with K\"ahler geometry. Indeed, an almost product structure $P$ on a $2n$-dimensional manifold whose eigenspaces are both of dimension $n$ is also known as a {\sl paracomplex structure}. When, in addition, $M$ carries a (necessarily) neutral metric $g$ compatible with $P$ in the sense that $P$ is an anti-orthogonal automorphism of each $(T_pM,g)$, then the triple $(M,g,P)$ is called almost paraHermitian, paraHermitian when $P$ is integrable, and paraK\"ahler when $\nabla P = 0$ (that a parallel almost product structure $P$ is indeed integrable is again analogous to the case for almost complex structures, deduced, in particular, by Walker [45--46]; see also [47]). The results above indicate this analogy is far reaching. 

Indeed, there is an extensive body of literature on paracomplex geometry, see the reviews [6], [7]. (Almost) paraHermitian geometry has also been called biLagrangian geometry because the eigenspaces of $P$ are totally null with respect to $g$, see [14]. These various names reflect different emphases upon $g$, $P$, and $\omega$; and perhaps because the topic is not widely known. Indeed, the result 4.2 has been presented independently on several occasions, e.g., in [39] and  [1], though it dates back to the earliest work in paracomplex geometry. Indeed, Rashevskij [33] studied the properties of a metric of the form (4.16) on a locally product $2n$-dimensional manifold and then Rozenfeld [34] explicitly drew the parallel with K\"ahler geometry. Paracomplex geometry has multiple independent origins, [20--21] being notable; a short, but informative, history of the origins of the subject may be found in [7].

We suggest that (almost) paraHermitian geometry is most naturally construed as a special kind of neutral geometry. ParaK\"ahler geometry is then a special kind of Walker geometry which we have called double Walker. Our intention here is not to simply make matters worse by adding yet further terminology. Both paraK\"ahler and double Walker are useful terms which stress different aspects of this interesting geometry. But it is the underlying neutral geometry which is, in our opinion, the natural setting. It may be argued, however, that the root of these geometries is paracomplex geometry, which involves no metric. We would suggest that the significance of (almost) paracomplex geometry within (almost) product geometries stems from its algebraic origins in the {\sl paracomplex\/} (also called {\sl Lorentz\/}) numbers. Though one can treat this algebra purely algebraically, it is most naturally regarded as a `normed' ($\Psi$-Euclidean) algebra, i.e., as an algebraic structure on ${\bf R}^{1,1}$, and thus naturally a feature of neutral geometry (it is argued in [18] that ${\bf R}^{1,1}$ is best understood as a neutral geometry rather than two-dimensional Lorentzian geometry due to the role of anti-isometries).

All this geometry also has a natural description as $G$-structures on the frame bundle $F(M)$ of $M$. Paracomplex geometry is a ${\bf GL(n;R)} \times {\bf GL(n;R)}$-structure. Let $P(M)$ denote the bundle of almost-product frames for the given reduction. Adding the compatible neutral metric $g$ allows one to construct from any almost-product frame, an almost-product frame that is also a Witt frame with respect to $g$ and thus a further reduction to $W(M)$, say. The group of this bundle of admissible frames is the intersection of ${\bf GL(n;R)} \times {\bf GL(n;R)}$ with {\bf O(n,n)}, but where the latter must be expressed in the form appropriate to Witt bases rather than $\Psi$-ON bases. The result is that the symmetry group of $W(M)$ takes the form
$$\left\{\,\pmatrix{B&{\bf 0}_n\cr {\bf 0}_n&{^\tau\! B}^{-1}\cr}:B \in {\bf GL(n;R)}\,\right\},$$
and thus is isomorphic to {\bf GL(n;R)}. Thus, almost paraHermitian geometry and its refinements may be viewed as certain {\bf GL(n;R)}-structures on $2n$-dimensional manifolds.

Returning now to the four-dimensional case of our primary concern, in which a double Walker geometry has an unambiguous canonical orientation, we encounter another co-incidence. Pleba\~nski [30] derived two canonical coordinate forms for a complex space-time (i.e., a four-dimensional complex manifold carrying a holomorphic Riemannian metric) whose only nontrivial curvature is the ASD Weyl curvature. Pleba\~nski's results carry over to the real category to apply to neutral metrics in four dimensions. With respect to the coordinates $(r,s,x,y)$, (4.16) takes the form
$$D = \pmatrix{\Omega_{rx}&\Omega_{ry}\cr \Omega_{sx}&\Omega_{sy}\cr}$$
which, given our choice of ordering of the coordinates, is of Pleba\~nski's {\sl first heavenly form\/}. Pleba\~nski's [30] first result is that, when $\det(D) = 1$ (Pleba\~nski's first heavenly equation), this form is a canonical form for a neutral metric in four dimensions whose only nontrivial curvature is the ASD Weyl curvature. Such metrics are thereby given, locally, by a single function $\Omega$ satisfying the first heavenly equation. What we have derived is the first heavenly form, without the constraint of the first heavenly equation, as a canonical form for a neutral metric with parallel complementary totally null distributions, and we recognize that this form is nothing but the paraK\"ahler form of the metric for a paraK\"ahler geometry. We note that (4.11) is a slight variation on Pleba\~nski's heavenly null tetrad.

To appreciate the significance of the first heavenly equation we state some definitions and facts well known to those familiar with complex general relativity and twistor theory, adapted to the context of neutral signature. But first we note in passing that Pleba\~nski's first heavenly equation actually appeared many years earlier in [35]; see also [12].
\vskip 24pt
\noindent {\bf 4.3 Definitions \& Facts}\hfil\break
Let $(M,g)$ be a four-dimensional neutral manifold. One can apply the Hodge star operator to the fully covariant Riemann tensor by applying it to either of the pair of (abstract) indices in which the Riemann tensor is skew to obtain two closely related `duals', see [28], \S 4.6. It then turns out that the Riemann tensor is SD/ASD with respect to one of these notions of duality iff it is also SD/ASD with respect to the other, and that when SD/ASD, the Ricci curvature vanishes and the Riemann tensor equals the SD/ASD Weyl tensor. One says $(M,g)$ is {\sl half-flat\/} when the Riemann tensor is SD/ASD; specifically, {\sl right flat\/} when the Riemann tensor is ASD, ie., the only nontrivial curvature is the ASD Weyl curvature, and {\sl left flat\/} when the Riemann tensor is SD, i.e., the only nontrivial curvature is the SD Weyl curvature. $(M,g)$ itself is said to be SD/ASD according as the Weyl curvature is SD/ASD, a weaker condition.

In $(M,g)$, one can construct, locally, parallel primed/unprimed spin frames iff $(M,g)$ is right/left flat. The forward implication follows by applying 3.1 to the two elements of the spin frame; the converse can be proved by a straightforward adaptation of the classical result that one can construct, locally, parallel ON frame fields iff the Riemann tensor vanishes (e.g., [37], pp. 261--263). In more modern terms, $(M,g)$ is right/left flat iff the induced connexion on the bundle of primed/unprimed spinors is flat, see. e.g., [16], \S II.9, or [37], pp. 402--403.
\vskip 24pt
Now suppose that ${\cal D}^\pi$ and ${\cal D}^\chi$ admit parallel local scaled representatives $\pi^{A'}$ and $\chi^{A'}$, at least on some common domain. On that domain, $\chi^{D'}\pi_{D'}$ is also constant, whence by a constant scaling of $\chi$ one can suppose, with out loss of generality, that $\chi^{D'}\pi_{D'} = 1$. For these choices of $\pi^{A'}$ and $\chi^{A'}$, one can choose, by 3.3, coordinates $(p,q,x,y)$ and $(w,z,r,s)$ as in (4.1--2) and proceed with the previous construction up to (4.16). Of course, $\xi^{A'} = \chi^{A'}$. By (4.10), $\det(D) = 1$, i.e., Pleba\~nski's first heavenly equation now holds. Moreover, the results of 3.1 must hold with respect to both $\pi^{A'}$ and $\xi^{A'}$. But these facts entail that $\tilde\Psi_{A'B'C'D'}$, $\Phi_{ABA'B'}$, and $S$ must vanish, whence $(M,g)$ is right flat.

What if we assume $\det(D) = 1$, equivalently $\xi^{A'} = \chi^{A'}$, in the construction leading up to (4.16) but not assume that parallel local scaled representatives are available? Pleba\~nski's [30] result assures one that $(M,g)$ is indeed right flat and that $\pi^{A'}$ and $\xi^{A'}$, characterized through the heavenly tetrad (4.11), are indeed parallel. 

Rather than repeat this analysis, however, we return to the assumption that ${\cal D}^\pi$ and ${\cal D}^\xi$, where from now on we employ the LSR $\xi^{A'}$ of $[\chi^{A'}]$ defined just after (4.10), are parallel complementary distributions without any presumption that either has parallel LSRs, i.e., to the double Walker (paraK\"ahler) case. The next result follows immediately from 2.5 and demonstrates the utility of viewing the geometry as Walker geometry.
\vskip 24pt
\noindent {\bf 4.4 Corollary}\hfil\break
{\it Let $(M,g,{\cal D}^\pi,[\pi^{A'}],{\cal D}^\xi,[\xi^{A'}])$ be a four-dimensional double Walker (paraK\"ahler) geometry. Then $[\pi^{A'}]$ and $[\xi^{A'}]$ are each WPSs of multiplicity at least two, whence $\tilde\Psi_{A'B'C'D'}\propto \pi_{(A'}\pi_{B'}\xi_{C'}\xi_{D')}$, i.e., $\tilde\Psi_{A'B'C'D'}$ is of type $\{22\}${\rm Ia}, or zero; $S = 0$ iff $\tilde\Psi_{A'B'C'D'}$ is zero; and the Ricci spinor is of the form $\Phi_{ABA'B'} = B_{AB}\pi_{(A'}\xi_{B')}$. This geometry is therefore Einstein iff $B_{AB} = 0$.}
\vskip 24pt
We now note that, when $\pi^{A'}$ can be chosen parallel, (3.8) substituted into (1.1) yields Pleba\~nski's second heavenly form for the metric of a right-flat space, though without the constraint of Pleba\~nski's second heavenly equation; rather (3.11) pertains. Hence, for a double Walker geometry, we follow Pleba\~nski [30] and define functions
$$u := \Omega_x \hskip 1.25in v := \Omega_y.\eqno(4.17)$$
Considering $(u,v,x,y)$ as functions of $(r,s,x,y)$, the Jacobian is
$$\pmatrix{\Omega_{xr}&\Omega_{xs}&\Omega_{xx}&\Omega_{xy}\cr \Omega_{yr}&\Omega_{ys}&\Omega_{yx}&\Omega_{yy}\cr 0&0&1&0\cr 0&0&0&1\cr} = \pmatrix{{^\tau\! D}&K\cr{\bf 0}_2&1_2\cr},\eqno(4.18)$$
which therefore has nonzero determinant. Thus, one may employ $(u,v,x,y)$ as local coordinates. Now (4.1) is still valid. One computes
$$\displaylines{du = \Omega_{xr}dr + \Omega_{xs}ds + \Omega_{xx}dx + \Omega_{xy}dy\cr
\noalign{\vskip -6pt}
\hfill\llap(4.19)\cr
\noalign{\vskip -6pt}
dv = \Omega_{yr}dr + \Omega_{ys}ds + \Omega_{yx}dx + \Omega_{yy}dy,\cr}$$
whence one can rewrite (4.11) as:
$$\displaylines{\ell_a = dx \hskip 1.25in \tilde m_a = dy\cr
\hfill n_a = D_{11}dr + D_{21}ds = \Omega_{rx}dr + \Omega_{sx}ds = du - \Omega_{xx}dx - \Omega_{xy}dy\hfill\llap(4.20)\cr
m_a = -D_{12}dr - D_{22}ds = -\Omega_{ry}dr - \Omega_{sy}ds = -dv + \Omega_{yx}dx + \Omega_{yy}dy,\cr}$$
and
$$du = n_a + \Omega_{xx}\ell_a + \Omega_{xy}\tilde m_a \hskip 1in dv = -m_a + \Omega_{xy}\ell_a + \Omega_{yy}\tilde m_a.\eqno(4.21)$$
One can now compute $g^{ab}$ with respect to the coordinates $(u,v,x,y)$:
$$\displaylines{g^{ab}(du,du) = 2\Omega_{xx} \hskip .75in g^{ab}(du,dv) = 2\Omega_{xy} \hskip .75in g^{ab}(dv,dv) = 2\Omega_{yy}\cr
g^{ab}(du,dx) = 1 \hskip .5in g^{ab}(du,dy) = 0 \hskip .5in g^{ab}(dv,dx) = 0 \hskip .5in g^{ab}(dv,dy) = 1\cr
g^{ab}(dx,dx) = g^{ab}(dx,dy) = g^{ab}(dy,dy) = 0.\cr}$$
Hence, one obtains, with respect to $(u,v,x,y)$, Walker's canonical form (1.1) for the metric
with
$$W = -2\pmatrix{\Omega_{xx}&\Omega_{xy}\cr \Omega_{xy}&\Omega_{yy}\cr}.\eqno(4.22)$$
Thus, $(u,v,x,y)$ are Walker coordinates for $(M,g,{\cal D}^\pi,[\pi^{A'}])$ and the heavenly tetrad (4.11), equivalently (4.20), is exactly the Walker null tetrad (2.11) for these Walker coordinates. Hence, (2.21), (2.26), and (2.32--33) give the curvature spinors with respect to the spin frames associated to the heavenly tetrad but with 4.4 in force, so in fact
$$\displaylines{\tilde\Psi_{A'B'C'D'} = {S \over 2}\pi_{(A'}\pi_{B'}\xi_{C'}\xi_{D')}\cr
\hfill\llap(4.23)\cr
\Phi_{ABA'B'} = B_{AB}\pi_{(A'}\xi_{B')}\cr}$$
with $B_{AB}$ as in (2.33). Moreover, the following equations must hold:
$$B + Sc = 0 \hskip .75in 0 = 6Bc - A + S(3c^2 - 1) = 3Bc - A - S \hskip .75in A_{AB} = 0,\eqno(4.24)$$
with $A_{AB}$ as in (2.33). Note that the first two equations entail $S^2 + AS + 3B^2 = 0$ in accord with 2.6(i). In particular, every four-dimensional paraK\"ahler metric must locally be of this form. Note also, from (4.10) and (4.18) that $[\partial_r,\partial_s,\partial_x,\partial_v] = [\partial_u,\partial_v,\partial_x,\partial_v]$, i.e., the locally product coordinates have the canonical orientation of the double Walker geometry.

Suppose now that $[\pi^{A'}]$ has a parallel LSR. One can employ 3.3 to choose the original set of coordinates $(p,q,x,y)$ satisfying (4.1) with that parallel LSR $\pi^{A'}$. Then one can introduce $\vartheta$ as in (3.7) and $W$ takes the form (3.8), where $\vartheta$ is subject to (3.12). The ASD Weyl curvature is given by (3.15). By 3.1, $S=0$ whence $\tilde\Psi_{A'B'C'D'} = 0$ by (4.23) (equivalently, $[\pi^{A'}]$ is a WPS of multiplicity four but $\xi^{A'}$ is also a WPS of multiplicity two; impossible). It also follows from 3.1, and 2.5 applied to $\xi^{A'}$, that $\Phi_{ABA'B'} = 0$ (equivalently, combine (3.16) and (4.23)). From (3.17), $\delta_A\delta_Bf = 0$. Thus, if only ${\cal D}^\pi$ has a parallel LSR, then a double Walker geometry $(M,g,{\cal D}^\pi,[\pi^{A'}],{\cal D}^\xi,[\xi^{A'}])$ is right flat.

It then follows from 3.2 that $\xi^{A'}$ can be rescaled so as to be parallel, say $f\xi^{A'}$ is parallel. Then $f = f\xi^{D'}\pi_{D'}$ is constant, whence $\xi^{A'}$ itself must be parallel. Although the description of the double Walker geometry relative to Walker coordinates is asymmetrical in the roles played by ${\cal D}^\pi$ and ${\cal D}^\xi$, the geometry itself is not and one can of course reverse the roles. Hence:
\vskip 24pt
\noindent {\bf 4.5 Theorem}\hfil\break
{\it Given a four-dimensional double Walker (paraK\"ahler) geometry, if either distribution admits a parallel LSR, the geometry is right flat on that domain and the other distribution also admits a parallel LSR on that domain.}
\vskip 24pt
Consequently, suppose now that both $\pi^{A'}$ and $\xi^{A'}$ in (4.11 \& 20) are parallel. Then $m^a \wedge n^b = \epsilon ^{AB}\xi^{A'}\xi^{B'}$ is parallel. One computes
$$m^a \wedge n^b = {ab - c^2 \over 4}\partial_u \wedge \partial_v + {c \over 2}\partial_u \wedge \partial_x - {a \over 2}\partial_u \wedge \partial_y + {b \over 2}\partial_v \wedge \partial_x - {c \over 2}\partial_v \wedge \partial_y + \partial_x \wedge \partial_y.\eqno(4.25)$$
Using (A1.8), one computes
$$\nabla_1(m^a \wedge n^b) = 0 = \nabla_2(m^a \wedge n^b).\eqno(4.26)$$
$$\eqalignno{\nabla_3(m^a \wedge n^b) &= (a_u + c_v)\left[{(c^2 - ab) \over 8}\partial_u \wedge \partial_v - {b \over 4}\partial_v \wedge \partial_x + {a \over 4}\partial_u \wedge \partial_y - {1 \over 2}\partial_x \wedge \partial_y\right]&(4.27)\cr
&\qquad + {2c_x - 2a_y - (c^2)_v + ba_v - ac_u \over 4}\partial_u \wedge \partial_x + {2c_x - 2a_y + ba_v - ac_u + 2ca_u \over 4}\partial_v \wedge \partial_y\cr}$$
$$\eqalignno{\nabla_4(m^a \wedge n^b) &= (b_v + c_u)\left[{(c^2 - ab) \over 8}\partial_u \wedge \partial_v - {b \over 4}\partial_v \wedge \partial_x + {a \over 4}\partial_u \wedge \partial_y - {1 \over 2}\partial_x \wedge \partial_y\right]&(4.28)\cr
&\qquad + {2b_x - 2c_y - 2cb_v + bc_v- ab_u  \over 4}\partial_u \wedge \partial_x + {2b_x - 2c_y + (c^2)_u + bc_v - ab_u \over 4}\partial_v \wedge \partial_y\cr}$$
Hence, under (3.2), and with $P$ as in (3.12), (4.27--28) become
$$\nabla_3(m^a \wedge n^b) = P_v(\partial_u \wedge \partial_x + \partial_v \wedge \partial_y) \hskip 1in \nabla_4(m^a \wedge n^b) = -P_u(\partial_u \wedge \partial_x + \partial_v \wedge \partial_y)$$
from which one deduces that $\xi^{A'}$ is parallel iff $P_u = P_v = 0$, i.e., $P$ depends only on $x$ and $y$. Following Pleba\~nski [30] again, let $F$ be an antiderivative of $P$ with respect to $x$: $F_x := P$. Define $\Theta := \vartheta - uF$. Then, $\Theta_{uu} = \vartheta_{uu}$, $\Theta_{uv} = \vartheta_{uv}$, $\Theta_{vv} = \vartheta_{vv}$, i.e., (3.8) holds with $\Theta$ replacing $\vartheta$, as do (3.13 \& 15). Of course, $\delta_A\delta_B P = 0$. Finally,
$$\Theta_{ux} + \Theta_{vy} - (\Theta_{uv})^2 + \Theta_{uu}\Theta_{vv} = \vartheta_{ux} -P + \vartheta_{vy} - (\vartheta_{uv})^2 + \vartheta_{uu}\vartheta_{vv} = 0,\eqno(4.29)$$
which is Pleba\~nski's second heavenly equation.

Thus, Pleba\~nski's second heavenly form is the special case of Walker's canonical form for a double Walker (paraK\"ahler) geometry with parallel LSRs for one, hence both, distributions. Of course, every right-flat neutral four-fold is, locally, such a double Walker geometry in many ways.
\vskip 24pt
\noindent {\section 5 Global Lifts of $[\pi^{A'}]$}
\vskip 12pt
Our considerations so far have been essentially local in nature. In this section we consider a four-dimensional Walker geometry $(M,g,{\cal D},[\pi^{A'}])$ for which the frame bundle admits a reduction with structure group ${\bf SO^\bfplus(2,2)}$. As noted just prior to 2.2, with respect to this reduction the bundle of projective spinors ${\bf P}S'_M$ is well defined globally, so the spinor field $[\pi^{A'}]$ corresponding to the Walker distribution $\cal D$ is a global section of ${\bf P}S'_M$. The question we address here is whether this section can be lifted to a section of the bundle $S'_M$? Obviously, we must assume the bundle $S'_M$ exists which, given the ${\bf SO^\bfplus(2,2)}$-reduction, requires that the second Stiefel-Whitney class of $M$ vanishes. We will also assume $M$ is connected and paracompact. 

Before proceeding, we say a few words about our assumptions. One can reformulate the existence of an ${\bf SO^\bfplus(2,2)}$-reduction in various ways, most obviously as the existence of a distribution of oriented two-planes, for which, when $M$ is compact, there are well known necessary and sufficient topological conditions. Matsushita has studied these issues and we direct the reader to [24] for a recent review. In general, the assumption of a Walker geometry admitting an ${\bf SO^\bfplus(2,2)}$-reduction cannot have a purely topological characterization as evidenced by the explicit examples in the literature cited in the Introduction of such Walker geometries on ${\bf R}^4$. It would, however, be of interest to determine in the compact case any topological conditions imposed by the existence of the Walker geometry beyond those equivalent to an ${\bf O(2,2)}$- or ${\bf SO^\bfplus(2,2)}$-reduction.

Returning to the question of a global lifting of $[\pi^{A'}]$, we first note that if $\pi^{A'}$ is such a global lift, then $\epsilon^{AB}\pi^{A'}\pi^{B'}$ is globally defined and defines an orientation for $\cal D$, see 2.4. Thus, orientability of $\cal D$ is a necessary condition. We elaborate on this point at the end of the section. By {\sl an open covering\/} of the Walker geometry, we shall mean an open covering ${\cal U} = \{\,U_i:i \in {\cal I}\,\}$ of $M$ such that each $U_i$ carries an LSR $\pi^{A'}_i$ of $[\pi^{A'}]$. We can frame our question as: given an open covering of the Walker geometry, can one scale the LSRs so that they agree on nontrivial intersections $U_{ij} := U_i \cap U_j$, $i$, $j \in {\cal I}$?
\vskip 24pt
\noindent {\bf 5.1 Lemma}\hfil\break
{\it With notation and assumptions as in the previous paragraphs, the obstruction to a global lifting of $[\pi^{A'}]$ to $S'_M$ is a specific element $\pi$ in $H^1(M,{\cal C}^*)$, where ${\cal C}^*$ is the sheaf of germs of nowhere-vanishing smooth functions on $M$. Thus, whether a global lift of $[\pi^{A'}]$ exists or not is determined by the topology of $M$ and the sheaf ${\cal C}^*$.}

{\it Proof.} For any open covering of the Walker geometry, on a nontrivial intersection $U_{ij}$ one has $\pi^{A'}_i = f_{ij}\pi^{A'}_j$, where $f_{ij}$ is a nowhere-vanishing, smooth function on $U_{ij}$. One easily checks that $U_{ij} \mapsto f_{ij}$ defines a \v Cech 1-cocycle with coefficients in ${\cal C}^*$ on $\cal U$. This assignment is well defined under restriction to refining coverings, i.e., on passing to a refinement, the induced LSRs define the 1-cocycle obtained by restricting the original 1-cocycle to the refinement.

It is routine to confirm that given two open coverings of the Walker geometry, with open coverings $\cal U$ and $\cal V$ of $M$ respectively, the 1-cocycles induced on a common refinement of $\cal U$ and $\cal V$ are cohomologous. Hence, the possible open coverings of the Walker geometry define a certain element $\pi$ of $H^1(M,{\cal C}^*)$.

Since $H^1({\cal U},{\cal S}) \to H^1(M,{\cal S})$ is injective for any sheaf $\cal S$, if $\pi$ is zero then for any open covering of the Walker geometry, $f_{ij}$ is a coboundary, $f_{ij} = (\delta h)_{ij} = h_jh_i^{-1}$. Then $h_i\pi^{A'}_i = h_j\pi^{A'}_j$, which thereby defines a global lift of $[\pi^{A'}]$.

Conversely if $\pi^{A'}$ is a global lift of $[\pi^{A'}]$ then one can construct open coverings of the Walker geometry with $f_{ij} = 1$, whence the 1-cocycle is a coboundary: $f_{ij} = (\delta h)_{ij} = h_j/h_i$, where $h_i \equiv 1$, for all $i$, i.e., $\pi=0$.

Thus, $\pi \in H^1(M,{\cal C}^*)$ is the obstruction to a global lifting of $[\pi^{A'}]$.
\vskip 24pt
If a Walker geometry $(M,g,{\cal D},[\pi^{A'}])$ satisfies the curvature conditions of 3.1, then 3.2 says one can construct an open covering of the Walker geometry with parallel LSRs.  Call such a {\sl parallel open covering\/} of the Walker geometry. One can repeat the argument of 5.1, the difference being that the 1-cocycle takes coefficients in a different sheaf.
\vskip 24pt
\noindent {\bf 5.2 Lemma}\hfil\break
{\it If a Walker geometry admits parallel open coverings, the obstruction to constructing a global parallel lifting of $[\pi^{A'}]$ is a certain cohomology class $\rho$ in $H^1(M,{\bf R}^*)$, where ${\bf R}^*$ is the constant (multiplicative) sheaf of nonzero real numbers (the sheaf of germs of locally constant ${\bf R}^*$-valued functions).}
\vskip 24pt
To explore these obstructions, consider the commutative diagram of short exact sequences:
\vskip 12pt
$$
\matrix{{\bf 0}&\mapright{}&{\bf S}^0&\hookright{j}&{\bf R}^*&\mapright{\psi}&{\bf R}^+&\mapright{}&{\bf 0}\cr
&&\mapdown{1}&&\mapdown{p}&&\mapdown{q}&&\cr
{\bf 0}&\mapright{}&{\bf S}^0&\hookright{i}&{\cal C}^*&\mapright{\phi}&{\cal C}^+&\mapright{}&{\bf 0}\cr}\eqno(5.1)$$
where ${\bf S}^0$ is the constant multiplicative sheaf with fibres isomorphic as groups to ${\bf Z}_2$, ${\bf R}^+$ is the constant (multiplicative) sheaf of positive real numbers, ${\cal C}^+$ is the sheaf of germs of positive {\bf R}-valued smooth functions, both $\phi$ and $\psi$ are $f \mapsto \vert f \vert$ for the appropriate domain, and the remaining mappings are the obvious inclusions. 

If $\cal C$ denotes the sheaf germs of smooth functions, $\exp:{\cal C} \to {\cal C}^+$ is a sheaf isomorphism, whence these two sheaves have the same cohomology. Since $\cal C$ is a fine sheaf, $H^p(M,{\cal C}^+) = 0$ for $p \geq 1$. Also $\exp:{\bf R} \to {\bf R}^+$ is a sheaf isomorphism, whence $H^p(M,{\bf R}^+) \cong H^p(M,{\bf R}) \cong H^p_{\rm DR}(M)$, where the last is the de Rham cohomology of $M$. One obtains from the long exact cohomology sequences the commutative diagram:
\vskip 12pt
$$
\matrix{\scriptstyle{\bf 0}&\scriptstyle{\mapright{}}&\scriptstyle{H^0(M,{\bf S}^0)}&\scriptstyle{\hookright{j_*}}&\scriptstyle{H^0(M,{\bf R}^*)}&\scriptstyle{\mapright{\psi_*}}&\scriptstyle{H^0(M,{\bf R}^+)}&\scriptstyle{\mapright{\delta_*}}&\scriptstyle{H^1(M,{\bf S}^0)}&\scriptstyle{\mapright{j_*}}&\scriptstyle{H^1(M,{\bf R}^*)}&\scriptstyle{\mapright{\psi_*}}&\scriptstyle{H^1(M,{\bf R}^+)}&\scriptstyle{\mapright{\delta_*}}\cr
&&\scriptstyle{\mapdown{1}}&&\scriptstyle{\mapdown{p_*}}&&\scriptstyle{\mapdown{q_*}}&&\scriptstyle{\mapdown{1}}&&\scriptstyle{\mapdown{p_*}}&&\scriptstyle{\mapdown{q_*}}&\cr
\scriptstyle{{\bf 0}}&\scriptstyle{\mapright{}}&\scriptstyle{H^0(M,{\bf S}^0)}&\scriptstyle{\hookright{i_*}}&\scriptstyle{H^0(M,{\cal C}^*)}&\scriptstyle{\mapright{\phi_*}}&\scriptstyle{H^0(M,{\cal C}^+)}&\scriptstyle{\mapright{\delta_*}}&\scriptstyle{H^1(M,{\bf S}^0)}&\scriptstyle{\mapright{i_*}}&\scriptstyle{H^1(M,{\cal C}^*)}&\scriptstyle{\mapright{\phi_*}}&\scriptstyle{H^1(M,{\cal C}^+)={\bf 0}}&\scriptstyle{\mapright{\delta_*}}\cr}\eqno(5.2)$$
\vskip 12pt
Beyond the portion shown, the bottom row is just
$$
\matrix{0&\mapright{\delta_*}&H^p(M,{\bf S}^0)&\mapright{i_*}&H^p(M,{\cal C}^*)&\mapright{\phi_*}&0\cr}$$
for each $p > 1$, whence $i_*$ is clearly an isomorphism for $p >1$.

Since ${\cal C}^+ \subset {\cal C}^*$ and ${\bf R}^+ \subset {\bf R}^*$, a cohomology class $c$ of $H^p(M,{\cal C}^+)$ also defines an element $d$ of $H^p(M,{\cal C}^*)$ which is mapped by $\phi_*$ to $c$, i.e., $\phi_*$ is always onto. Similarly, $\psi_*$ is always onto. It follows by exactness, equivalently by their very definition, that all the connecting homomorphisms $\delta_*$ are trivial mappings. 

Since $\phi_*:H^1(M,{\cal C}^*) \to H^1(M,{\cal C}^+)$ is a surjection onto a trivial space and $\delta_*$ a trivial mapping, it follows that $i_*:H^1(M,{\bf S}^0) \to H^1(M,{\cal C}^*)$ is an isomorphism. Hence, for $p > 0$,
$$H^p(M,{\cal C}^*) \cong H^p(M,{\bf S}^0) \cong H^p(M,{\bf Z}_2).\eqno(5.3)$$

Since the $\delta_*$'s are trivial, one can isolate the following commutative diagram, with exact rows, from (5.2):
\vskip 12pt
$$
\matrix{{\bf 0}&\mapright{}&H^1(M,{\bf S}^0)&\mapright{j_*}&H^1(M,{\bf R}^*)&\mapright{\psi_*}&H^1(M,{\bf R}^+)&\mapright{}&{\bf 0}\cr
&&\mapdown{1}&&\mapdown{p_*}&&\mapdown{q_*}&&\cr
{\bf 0}&\mapright{}&H^1(M,{\bf S}^0)&\mapright{i_*}&H^1(M,{\cal C}^*)&\mapright{\phi_*}&H^1(M,{\cal C}^+) = {\bf 0}&\mapright{}&{\bf 0}\cr}\eqno(5.4)$$
\vskip 12pt
Thus
$$H^1_{\rm DR}(M) \cong H^1(M,{\bf R}) \cong H^1(M,{\bf R}^+) \cong {H^1(M,{\bf R}^*) \over H^1(M,{\bf S}^0)}.\eqno(5.5)$$
\vskip 12pt
\noindent {\bf 5.3 Theorem}\hfil\break
{\it Suppose $M$ is connected and paracompact. Let $(M,g,{\cal D},[\pi^{A'}])$ be a four-dimensional Walker geometry. Suppose $M$ admits an ${\bf SO^\bfplus(2,2)}$-reduction of the frame bundle, whence $[\pi^{A'}]$ is a global section of the projective spinor bundle ${\bf P}S'_M$. Suppose further that $w_2(M) = 0$, whence $(M,g)$ admits spinor structures. The obstruction to a global lifting of $[\pi^{A'}]$ to a section of $S'_M$ is an element of}
$$H^1(M,{\cal C}^*) \cong H^1(M,{\bf Z}_2) \cong \Hom\left(H_1(M,{\bf Z}),{\bf Z}_2\right) \cong \Hom\left(\pi_1(M),{\bf Z}_2\right).\eqno(5.6)$$
The nontrivial elements of $\Hom\left(\pi_1(M),{\bf Z}_2\right)$ are in bijective correspondence with the subgroups of $\pi_1(M)$ of index two, which are, in turn, in bijective correspondence with the isomorphism classes of two-fold covering spaces of $M$.

If the curvature conditions of 3.1 pertain, the open coverings of the Walker geometry define $\pi \in H^1(M,{\cal C}^*)$ and the parallel open coverings define $\rho \in H^1(M,{\bf R}^*)$; the latter group is isomorphic to $H^1(M,{\bf Z}_2)$ precisely when $H^1_{\rm DR}(M)$ vanishes. From (5.4), $\pi = i_*(\beta) = (p_* \circ j_*)(\beta) = p_*(\rho)$, for a unique $\beta \in H^1(M,{\bf S}^0)$, whence $\rho - j_*(\beta) \in \ker p_*$. If $\pi = 0$, then $\beta = 0$ and $\rho \in \ker p_*$ but is not necessarily zero, i.e., $[\pi^{A'}]$ may have a global lifting but not a parallel global lifting, even when the conditions of 3.1 pertain. If $H^1_{\rm DR}(M) = 0$, however, then $j_*$, whence $p_*$, are isomorphisms, $\rho = j_*(\beta)$, $p_*(\rho) = \pi$, and $\pi = 0$ iff $\beta = 0$ iff $\rho = 0$.

{\it Remarks.} The isomorphisms in (5.6) and thereafter are standard interpretations of $H^1(M,{\bf Z}_2)$ following, in the first instance, from the Universal Coefficient Theorem.

When the curvature conditions of 3.1 pertain, one can construct a parallel open covering using a {\sl simple\/} covering in the sense of [16], pp. 167--168. The cohomology of $M$ is then given by the \v Cech cohomology with respect to this covering. Let $\{f_{ij}\}$ be the 1-cocycle obtained from the parallel open covering, which has coefficients in ${\bf R}^* \subset {\cal C}^*$, i.e., it represents both $\rho$ and $\pi$. Now $\vert f_{ij}\vert$ defines a 1-cocycle with coefficients in ${\cal C}^+$. As $H^1(M,{\cal C}^+) = 0$, then $\vert f_{ij}\vert = h_j/h_i$, for some 0-cochain $\{h_i\}$ with coefficients in ${\cal C}^+$. But $\{f_{ij}/\vert f_{ij} \vert\}$ is also a 1-cocycle with coefficients in ${\bf S}^0 \subset {\cal C}^*$, whence $\{f_{ij}\}$ is cohomologous to $\{f_{ij}/\vert f_{ij} \vert\} = \{f_{ij}h_j/h_i\}$ as 1-cocycles with coefficients in ${\cal C}^*$ and $\{f_{ij}/\vert f_{ij} \vert\}$ represents $\beta$. If $\pi = 0$, then $f_{ij} = f_j/f_i$, for some 0-cochain $\{f_j\}$ with coefficients in ${\cal C}^*$. Then $\vert f_{ij}\vert = \vert f_j \vert/\vert f_i \vert$ and $f_{ij}/\vert f_{ij} \vert = (f_j/\vert f_j\vert)(f_i/\vert f_i \vert)^{-1}$, i.e., $\beta$ is trivial. But $\{f_{ij}\}$ may not be a coboundary with coefficients in ${\bf R}^*$. Now $\vert f_{ij} \vert$ also defines a 1-cocycle with coefficients in ${\bf R}^+$, so if $H^1_{\rm DR}(M) \cong H^1(M,{\bf R}^+) = 0$, then $\vert f_{ij} \vert = g_j/g_i$, where $\{g_i\}$ is a 0-cochain with coefficients in ${\bf R}^+$. Now $\{f_{ij}\}$ is cohomologous to $\{f_{ij}/\vert f_{ij} \vert\} = \{f_{ij}g_j/g_i\}$ as 1-cocycles with coefficients in ${\bf R}^*$, i.e., $j_*(\beta) = \rho$. Now $\pi = 0$ iff $\beta = 0$ iff $\rho = 0$.

Finally, if $\{f_{ij}\}$ is the 1-cocycle defined by an open covering of the Walker geometry, so $\pi_i^{A'} = f_{ij}\pi^{A'}_j$, then $\Sigma_i := \epsilon^{AB}\pi^{A'}_i\pi^{B'}_i  = (f_{ij})^2\epsilon^{AB}\pi^{A'}_j\pi^{B'}_j =: (f_{ij})^2\Sigma_j$. Now $(f_{ij})^2$ is a smooth, positive function, which indicates $\cal D$ must be orientable. Indeed, $(f_{ij})^2/(\vert f_{ij}\vert)^2$ is a one-cocycle with coefficients in ${\bf S}^0$ which everywhere takes the trivial value and thus belongs to the trivial cohomology class. In other words, $\{(f_{ij}/\vert f_{ij}\vert)^2\}$ must represent the first Stiefel-Whitney class of $\cal D$ viewed as a bundle. We thus see that $\cal D$ must be orientable merely because we have assumed the existence of a global spinor bundle, since this fact allows one to compare LSRs on overlaps and deduce that the $\Sigma_i$s are positive multiples of each other. On the other hand, the cohomology class $\pi = \{f_{ij}/\vert f_{ij}\vert\}$ need not be trivial, so orientability of $\cal D$ is merely a necessary condition of the context in which the question of global lifts arises.
\vskip 24pt
\noindent {\section Acknowledgments}
\vskip 12pt
We thank Mike Eastwood for a useful tip in connexion with \S 5, Pedro Gadea for kindly and expeditiously providing us with copies of his papers on paracomplex geometry, and an anonymous referee for remarks that improved the paper.
\vskip 24pt
\noindent {\section Appendix One: Local Geometry with Respect to}\hfil\break
{\section Walker Coordinates}
\vskip 12pt
In this appendix we record, for ease of reference in this and other papers, local coordinate expressions for standard geometric objects. We also specify our choice of conventions and in this sense standardize the coordinate expressions.

We employ the abstract index notation of [28]; italic indices will be `abstract', i.e., serve to denote the tensor space to which the object they are attached to belongs, while bold Roman indices are `concrete', i.e., they take numerical values and typically serve to label components of geometric objects with respect to some basis or the elements of a basis themselves. We also employ the standard summation convention for concrete indices when convenient; while repeated abstract indices in a formula, one as superscript and one as subscript, indicates the standard pairing between a linear space and its dual.

For exterior algebras, of the two conventions commonly found in the literature, we employ those of [40] (and many other texts such as [36]). This choice simplifies the expression for volume forms and elements. We shall, however, retain the definition of symmetrization and skew symmetrization, denoted by round and square brackets respectively, of abstract and concrete indices employed in [28]. Thus, if $v_a$ and $w_b$ are two vectors, then
$$v^a \wedge w^b = v^a \otimes w^b - w^a \otimes v^b = 2v^{[a}w^{b]}.\eqno({\rm A}1.1)$$
Given a real linear space ${\cal V}$ equipped with a scalar product $g$ (of arbitrary signature), the induced scalar product on the exterior algebra $\Lambda^p({\cal V})$ of $p$-vectors is given by
$$(U,W) := {1 \over p!}U^{a_1\ldots a_p}W^{b_1 \ldots b_p}g_{a_1b_1}\ldots g_{a_pb_p},\eqno({\rm A}1.2)$$
where $U^{a_1 \ldots a_p}$ and $W^{b_1 \ldots b_p}$ denote the multivectors $U$ and $W$ as tensors. The scalar product $g_{ab}$ induces an identification of ${\cal V}$ with its linear dual ${\cal V}_\bullet$ by $v \mapsto g(v,\ )$, which we shall call the correlation $\xi_g:{\cal V} \to {\cal V}_\bullet$ with inverse $\xi_g^{-1}:{\cal V}_\bullet \to {\cal V}$. These identifications extend to arbitrary tensor spaces and are conveniently represented by index lowering and raising via $g_{ab}$ and $g^{ab}$ using abstract indices, where $g^{ab}$ is the scalar product induced on ${\cal V}_\bullet$ in the standard way. The induced scalar product on the exterior algebra $\Lambda^p({\cal V}_\bullet)$ is given by a formula analogous to (A1.2) with forms replacing multivectors and $g^{ab}$ replacing $g_{ab}$. The definition (A1.2) ensures that if $\{v_1,\ldots,v_n\}$ is a pseudo-orthonormal ($\Psi$-ON) basis  for ${\cal V}$ then the multivectors $v_{i_1} \wedge \ldots \wedge v_{i_p}$, $i_1 < \cdots < i_p$, form a $\Psi$-ON basis for $\Lambda^p({\cal V})$.

Now suppose $g$ is of signature $(r,s)$, $r+s = n$, that ${\cal V}$ is oriented, and that $\{v_1,\ldots,v_n\}$ is an oriented $\Psi$-ON basis with dual basis $\{\phi^1,\ldots,\phi^n\}$. Putting $\nu^j := \xi_g(v_j)$, then $\nu^j = \epsilon_j\phi^j$, where $\epsilon_j = \pm$ according as $j \leq r$ or $j > r$. The orientation class $[v_1,\ldots,v_n]$ can also be represented by the equivalence class of $v_1 \wedge \ldots \wedge v_n$ in $\Lambda^n({\cal V})/{\bf R}^+$ or, equivalently, by the equivalence class of $\phi^1 \wedge \ldots \wedge \phi^n$ in $\Lambda^n({\cal V}_\bullet)/{\bf R}^+$ (as those $n$-forms which are positive when evaluated on oriented bases). The $n$-fold wedge product of the elements of any oriented $\Psi$-ON basis yield one and the same $n$-vector in $\Lambda^n({\cal V})$. This element we call the {\sl volume element\/} of $({\cal V},g)$ and represent as a tensor by $V^{a_1 \ldots a_n}$. Similarly, the dual bases of oriented $\Psi$-ON bases all generate the same $n$-form $\phi^1 \wedge \ldots \wedge \phi^n$, which we call the {\sl volume form\/} of $({\cal V},g)$ and denote by $e_{a_1 \ldots e_n}$ as a tensor. Note that
$$V_{a_1 \ldots a_n} := g_{a_1b_1}\ldots g_{a_nb_n}V^{b_1\ldots b_n} = \xi_g(v_1 \wedge \ldots \wedge v_n) = \nu^1 \wedge \ldots \wedge \nu^n = (-1)^se_{a_1 \ldots a_n}.\eqno({\rm A}1.3)$$
The Hodge star operator is defined in the usual way on $\Lambda^p({\cal V}_\bullet)$ by
$$* \alpha := {1 \over p!}e_{i_1 \ldots i_n}\alpha_{j_i \ldots j_p}g^{i_1j_1}\ldots g^{i_pj_p},\eqno({\rm A}1.4)$$
where $\alpha_{a_1 \ldots a_p}$ denotes the $p$-form $\alpha$ as a tensor. The well known formula $\alpha \wedge *\beta = (\alpha,\beta)e_{a_1\ldots a_n}$ pertains, where $(\alpha,\beta)$ is the induced scalar product on $\Lambda^p({\cal V}_\bullet)$. Similarly, one defines the Hodge star operator on $\Lambda^p({\cal V})$ by
$$* U := {1 \over p!}V^{i_1 \ldots i_n}U^{j_i \ldots j_p}g_{i_1j_1}\ldots g_{i_pj_p},\eqno({\rm A}1.5))$$
and $U \wedge * W = (U,W)V^{a_1 \ldots a_n}$. These two Hodge star operators are related as follows:
$$\xi_g(*U) = (-1)^s*\xi_g(U).\eqno({\rm A}1.6)$$
When the context is clear, we may write either $\Lambda^p({\cal V})$ or $\Lambda^p({\cal V}_\bullet)$ simply as $\Lambda^p$. In the four-dimensional case, $\Lambda^2_\pm$ will denote the subspaces of self dual(SD)/anti-self dual (ASD) multivectors or forms.

Now let $(M,g,{\cal D})$ be a Walker four-manifold. We typically denote a set of Walker coordinates by $(u,v,x,y)$, but it is preferable to write coordinate expressions for geometrical objects in a form independent of the choice of letters used; to this end the Walker coordinates will be designated by the numerals 1, 2, 3, 4 respectively, when convenient; in particular, $\partial_1 := \partial_u$, $\partial_2 := \partial_v$, $\partial_3 := \partial_x$, and $\partial_4 := \partial_y$. The canonical form of the metric is given in (1.1) and (1.7).
\vskip 24pt
\noindent {\bf A1.1 The Christoffel Symbols}\hfil\break
The Christoffel symbols in a Walker coordinate system are:
$$\displaylines{\Gamma^{\bf i}_{11} = \Gamma^{\bf i}_{12} = \Gamma^{\bf i}_{22} = 0\cr
\noalign{\vskip 6pt}
\Gamma^{\bf i}_{13} = \cases{{1 \over 2}a_1,&${\bf i}=1$;\cr \cr {1 \over 2}c_1,&${\bf i}=2$;\cr \cr 0,& otherwise;\cr} \hskip 1.25in \Gamma^{\bf i}_{14} = \cases{{1 \over 2}c_1,&${\bf i}=1$;\cr \cr {1 \over 2}b_1,&${\bf i}=2$;\cr \cr 0,& otherwise;\cr}\cr
\noalign{\vskip 6pt}
\Gamma^{\bf i}_{23} = \cases{{1 \over 2}a_2,&${\bf i}=1$;\cr \cr {1 \over 2}c_2,&${\bf i}=2$;\cr \cr 0,& otherwise;\cr} \hskip 1.25in \Gamma^{\bf i}_{24} = \cases{{1 \over 2}c_2,&${\bf i}=1$;\cr \cr {1 \over 2}b_2,&${\bf i}=2$;\cr \cr 0,& otherwise;\cr}\cr
\noalign{\vskip 6pt}
\Gamma^{\bf i}_{33} = \cases{{1 \over 2}(aa_1+ca_2+a_3),&${\bf i}=1$;\cr \cr {1 \over 2}(2c_3+ca_1+ba_2-a_4),&${\bf i}=2$;\cr \cr -{1 \over 2}a_1,&${\bf i}=3$;\cr \cr -{1 \over 2}a_2,&${\bf i}=4$;\cr}\hskip 1.25in \Gamma^{\bf i}_{34} = \cases{{1 \over 2}(a_4+ac_1+cc_2),&${\bf i}=1$;\cr \cr {1 \over 2}(b_3+cc_1+bc_2),&${\bf i}=2$;\cr \cr -{1 \over 2}c_1,&${\bf i}=3$;\cr \cr -{1 \over 2}c_2,&${\bf i}=4$;\cr}\cr
\noalign{\vskip 6pt}
\Gamma^{\bf i}_{44} = \cases{{1 \over 2}(2c_4+ab_1+cb_2-b_3),&${\bf i}=1$;\cr \cr {1 \over 2}(b_4+cb_1+bb_2),&${\bf i}=2$\cr \cr -{1\over 2}b_1,&${\bf i}=3$;\cr \cr -{1 \over 2}b_2,&${\bf i}=4$.\cr}\cr}$$
\vskip 24pt
\noindent {\bf A1.2 The Geodesic Equations}\hfil\break
With $X^{\bf a} = (u,v,x,y)$, the Lagrangian is 
$${\cal L} := (1/2)g_{\bf ab}\dot X^{\bf a}\dot X^{\bf b} = \dot u\dot x + \dot v\dot y + {a \over 2}(\dot x)^2 + {b \over 2}(\dot y)^2 + c\dot x\dot y,$$
which is invariant under the interchange
$$u\ \leftrightarrow\ v \hskip .75in x\ \leftrightarrow\ y\hskip .75in a\ \leftrightarrow\ b.\eqno({\rm A}1.7)$$
Note that this interchange constitutes a fundamental symmetry of Walker coordinates which all Walker coordinate expressions must manifest. The Euler equations are:
$$\displaylines{\left({\partial{\cal L} \over \partial\dot u}\right)^\cdot - {\partial{\cal L} \over \partial u} = \ddot x - {a_u \over 2}(\dot x)^2 - {b_u \over 2}(\dot y)^2 - c_u\dot x\dot y = 0\cr
\noalign{\vskip 6pt}
\left({\partial{\cal L} \over \partial\dot v}\right)^\cdot - {\partial{\cal L} \over \partial v} = \ddot y - {a_v \over 2}(\dot x)^2 - {b_v \over 2}(\dot y)^2 - c_v\dot x\dot y = 0\cr}$$
$$\eqalign{\left({\partial{\cal L} \over \partial\dot x}\right)^\cdot - {\partial{\cal L} \over \partial x} &= \ddot u + a_u\dot u\dot x + a_v\dot v \dot x + c_u\dot u\dot y + c_v\dot v\dot y\cr
&+ {1 \over 2}(aa_u + ca_v+a_x)(\dot x)^2 + {1 \over 2}(ab_u+cb_v-b_x+2c_y)(\dot y)^2 + (a_y+ac_u+cc_v)\dot x\dot y = 0\cr}$$
$$\eqalign{\left({\partial{\cal L} \over \partial\dot y}\right)^\cdot - {\partial{\cal L} \over \partial y} &= \ddot v + c_u\dot u\dot x + c_v\dot v \dot x + b_u\dot u\dot y + b_v\dot v\dot y\cr
&+ {1 \over 2}(ca_u + ba_v+2c_x-a_y)(\dot x)^2 + {1 \over 2}(cb_u+bb_v+b_y)(\dot y)^2 + (b_x+cc_u+bc_v)\dot x\dot y = 0,\cr}$$
from which in fact the Christoffel symbols may be directly read off. The geodesic equations were previously published in [13].

In particular, putting $x = \hbox{constant}$, $y = \hbox{constant}$ reduces these four equations to the pair
$$\ddot u = 0 \hskip 1.5in \ddot v = 0,$$
i.e.,
$$(u,v,x,y) = (\alpha s + \beta,\gamma s + \delta,\mu,\nu),$$
for any constants $\alpha$, $\beta$, $\gamma$, $\delta$, $\mu$ and $\nu$, is a null geodesic lying in the integral surface of $\cal D$ through $(\beta,\delta,\mu,\nu)$. As the ratio $\alpha/\gamma$ varies, one obtains a one-parameter family of null geodesics lying in, and sweeping out, this $\alpha$-surface.

Writing $\nabla_{\bf i} := \nabla_{\partial_{\bf i}}$, and noting that of course $\nabla_{\bf i}\partial_{\bf j} = \Gamma^{\bf k}_{\bf ij}\partial_{\bf k} = \nabla_{\bf j}\partial_{\bf i}$, one computes the covariant derivatives of the coordinate basis:
$$\displaylines{\nabla_1\partial_1 = 0 \hskip .75in \nabla_2\partial_2 = 0 \hskip .75in \nabla_2\partial_1 = 0 = \nabla_1\partial_2\cr
\nabla_3\partial_3 = {1 \over 2}(a_3+ca_2+aa_1)\partial_1 + {1 \over 2}(2c_3-a_4+ba_2+ca_1)\partial_2 - {1 \over 2}(a_1\partial_3+a_2\partial_4)\cr
\nabla_4\partial_4 = {1 \over 2}(ab_1 + cb_2 - b_3 + 2c_4)\partial_1 + {1\over 2}(b_4 + cb_1 + bb_2)\partial_2 -{1 \over 2}(b_1\partial_3 + b_2\partial_4)\cr
\noalign{\vskip -6pt}
\hfill\llap({\rm A}1.8)\cr
\noalign{\vskip -6pt}
\nabla_1\partial_3 = {1 \over 2}(a_1\partial_1 + c_1\partial_2) = \nabla_3\partial_1 \hskip 1in
\nabla_1\partial_4 = {1 \over 2}(c_1\partial_1 + b_1\partial_2) = \nabla_4\partial_1\cr
\nabla_2\partial_3 = {1 \over 2}(a_2\partial_1+c_2\partial_2) = \nabla_3\partial_2\hskip 1in
\nabla_2\partial_4 = {1 \over 2}(c_2\partial_1+b_2\partial_2) = \nabla_4\partial_2\cr
\nabla_3\partial_4 = {1 \over 2}(a_4 + ac_1 + cc_2)\partial_1 + {1 \over 2}(b_3 + cc_1 + bc_2)\partial_2 -{ 1\over 2}(c_1\partial_3 + c_2\partial_4) = \nabla_4\partial_3.\cr}$$

One confirms that $D = \langle \partial_u,\partial_v \rangle_{\bf R}$ is indeed parallel. If $\partial_u$ and $\partial_v$ are actually parallel, then $a$, $b$ \& $c$ depend only on $x$ and $y$. In fact, Walker [42], [44] asserted that in this case, one can actually choose the coordinates $(u,v,x,y)$ so that $a=c=0$ and $b$ is a function of $x$ \& $y$ only.
\vskip 24pt
\noindent {\bf A1.3 Riemann Curvature}\hfil\break
We use the curvature conventions of [27], which for the Riemann curvature are
$$\displaylines{\hfill R(X,Y)Z = \nabla_{[X,Y]}Z - [\nabla_X,\nabla_Y]Z\hfill\llap({\rm A}1.9)\cr
\hfill R_{\bf ijkl} := R(X_{\bf i},X_{\bf j},X_{\bf k},X_{\bf l}) := g\bigl(R(X_{\bf k},X_{\bf l})X_{\bf j},X_{\bf i}\bigr)\hfill\llap({\rm A}1.10)\cr}$$
agreeing with those for the Riemann curvature in [28]. The coordinate expression for $R^{\bf i}{}_{\bf jkl}$ is then:
$$R^{\bf i}{}_{\bf jkl} = \Gamma^{\bf i}_{{\bf kj},{\bf l}} - \Gamma^{\bf i}_{{\bf lj},{\bf k}} + \Gamma^{\bf m}_{\bf kj}\,\Gamma^{\bf i}_{\bf lm} - \Gamma^{\bf m}_{\bf lj}\,\Gamma^{\bf i}_{\bf km}.\eqno({\rm A}1.11)$$
Direct computation yields:
$$\displaylines{R^{\bf i}{}_{{\bf j}12} = 0\cr
\noalign{\vskip 12pt} 
R^{\bf i}{}_{113} = \cases{-{1 \over 2}a_{11},&${\bf i}=1$;\cr \cr -{1 \over 2}c_{11},&${\bf i}=2$;\cr \cr 0,& otherwise;\cr} \hskip 1.25in R^{\bf i}{}_{213} = \cases{-{1 \over 2}a_{12},&${\bf i}=1$;\cr \cr -{1 \over 2}c_{12},&${\bf i}=2$;\cr \cr 0,& otherwise;\cr}\cr
\noalign{\vskip 12pt}
R^{\bf i}{}_{313} = \cases{-{1 \over 2}(aa_{11} + ca_{12}),&${\bf i}=1$;\cr \cr -{1 \over 4}(2c_{13} + 2ca_{11} + 2ba_{12} - 2a_{14} + b_1a_2 - c_1c_2),&${\bf i}=2$;\cr \cr {1 \over 2}a_{11},&${\bf i}=3$;\cr \cr {1 \over 2}a_{12},&${\bf i}=4$;\cr}\cr
\noalign{\vskip 12pt}
R^{\bf i}{}_{413} = \cases{{1 \over 4}(2c_{13} - a_{14} - ac_{11} - 2cc_{12} + a_2b_1 - c_1c_2),&${\bf i}=1$;\cr \cr -{1 \over 2}(cc_{11} + bc_{12}),&${\bf i}=2$;\cr \cr -{1 \over 2}c_{11},&${\bf i}=3$;\cr \cr -{1 \over 2}c_{12},&${\bf i}=4$;\cr}\cr
\noalign{\vskip 12pt}
R^{\bf i}{}_{114} = \cases{-{1 \over 2}c_{11},&${\bf i}=1$;\cr \cr -{1 \over 2}b_{11},&${\bf i}=2$\cr \cr 0,& otherwise;\cr} \hskip 1.5in R^{\bf i}{}_{214} = \cases{-{1 \over 2}c_{12},&${\bf i}=1$;\cr \cr -{1 \over 2}b_{12},&${\bf i}=2$\cr \cr 0,& otherwise;\cr}\cr
\noalign{\vskip 12pt}
R^{\bf i}{}_{314} = \cases{-{1 \over 2}(ac_{11} + cc_{12}),&${\bf i}=1$;\cr \cr {1 \over 4}\left(2c_{14} - 2b_{13} - 2cc_{11} - 2bc_{12} + a_1b_1 - b_1c_2 + b_2c_1- (c_1)^2\right),&${\bf i}=2$;\cr \cr {1 \over 2}c_{11},&${\bf i}=3$;\cr \cr {1 \over 2}c_{12},&${\bf i}=4$;\cr}\cr
\noalign{\vskip 12pt}
R^{\bf i}{}_{414} = \cases{-{1 \over 4}\left(2c_{14} + 2ab_{11} + 2cb_{12} - 2b_{13} + a_1b_1 + b_2c_1- b_1c_2 - (c_1)^2\right),&${\bf i}=1$;\cr \cr -{1 \over 2}(cb_{11} + bb_{12}),&${\bf i}=2$;\cr \cr {1 \over 2}b_{11},&${\bf i}=3$\cr \cr {1 \over 2}b_{12},&${\bf i}=4$;\cr}\cr
\noalign{\vskip 12pt}
R^{\bf i}{}_{123} = \cases{-{1 \over 2}a_{12},&${\bf i}=1$;\cr \cr -{1 \over 2}c_{12},&${\bf i}=2$;\cr \cr 0,& otherwise;\cr} \hskip 1.5in R^{\bf i}{}_{223} = \cases{-{1 \over 2}a_{22},&${\bf i}=1$;\cr \cr -{1 \over 2}c_{22},&${\bf i}=2$;\cr \cr 0,& otherwise;\cr}\cr
\noalign{\vskip 12pt}
R^{\bf i}{}_{323} = \cases{-{1 \over 2}(aa_{12} + ca_{22}),&${\bf i}=1$;\cr \cr -{1 \over 4}\left(2c_{23} + 2ca_{12} + 2ba_{22} - 2a_{24} + a_1c_2 + a_2b_2 - a_2c_1 - (c_2)^2\right),&${\bf i}=2$;\cr \cr {1 \over 2}a_{12},&${\bf i}=3$;\cr \cr {1 \over 2}a_{22},&${\bf i}=4$;\cr}\cr
\noalign{\vskip 12pt}
R^{\bf i}{}_{423} = \cases{{1 \over 4}\left(2c_{23} - 2a_{24} - 2ac_{11} - 2cc_{22} - a_2c_1+ a_1c_2 + a_2b_2 - (c_2)^2\right),&${\bf i}=1$;\cr \cr -{1 \over 2}(cc_{12} + bc_{22}),&${\bf i}=2$;\cr \cr {1 \over 2}c_{12},&${\bf i}=3$;\cr \cr{1 \over 2}c_{22},&${\bf i}=4$;\cr}\cr
\noalign{\vskip 12pt}
R^{\bf i}{}_{124} = \cases{-{1 \over 2}c_{12},&${\bf i}=1$;\cr \cr -{1 \over 2}b_{12},&${\bf i}=2$;\cr \cr 0,& otherwise;\cr} \hskip 1.5in R^{\bf i}{}_{224} = \cases{-{1 \over 2}c_{22},&${\bf i}=1$;\cr \cr -{1 \over 2}b_{22},&${\bf i}=2$;\cr \cr 0,& otherwise;\cr}\cr
\noalign{\vskip 12pt}
R^{\bf i}{}_{324} = \cases{-{1 \over 2}(ac_{12} + cc_{22}),&${\bf i}=1$;\cr \cr {1 \over 4}(2c_{24} - 2b_{32} - 2cc_{12} - 2bc_{22} + a_2b_1 - c_1c_2),&${\bf i}=2$;\cr \cr {1 \over 2}c_{12},&${\bf i}=3$;\cr \cr {1 \over 2}c_{22},&${\bf i}=4$;\cr}\cr
\noalign{\vskip 12pt}
R^{\bf i}{}_{424} = \cases{-{1 \over 4}(2c_{24} + 2ab_{12} + 2cb_{22} - 2b_{23} + a_2b_1 - c_1c_2),&${\bf i}=1$;\cr \cr -{1 \over 2}(cb_{12} + bb_{22}),&${\bf i}=2$;\cr \cr {1 \over 2}b_{12},&${\bf i}=3$;\cr \cr {1 \over 2}b_{22},&${\bf i}=4$;\cr}\cr
\noalign{\vskip 12pt}
R^{\bf i}{}_{134} = \cases{{1 \over 4}(2a_{14} - 2c_{13} + c_1c_2 - a_2b_1),&${\bf i}=1$;\cr \cr {1 \over 4}\left(2c_{14} - 2b_{13} + a_1b_1 + b_2c_1 - b_1c_2 - (c_1)^2\right),&${\bf i}=2$;\cr \cr 0,& otherwise;\cr}\cr
\noalign{\vskip 12pt}
R^{\bf i}{}_{234} = \cases{{1 \over 4}\left(2a_{24} - 2c_{23} + a_2c_1 - a_1c_2 - b_2c_2 + (c_2)^2\right),&${\bf i}=1$;\cr \cr {1 \over 4}(2c_{24} - 2b_{23} + a_2b_1 - c_1c_2),&${\bf i}=2$;\cr \cr 0,& otherwise;\cr}\cr
\noalign{\vskip 12pt}
R^{\bf i}{}_{334} = \cases{{ 1\over 4}\left(2aa_{14} + 2ca_{24} - 2ac_{13} - 2cc_{23} + ca_2c_1 + ac_1c_2 - ca_1c_2 - aa_2b_1 - ca_2b_2 + c(c_2)^2\right),&${\bf i}=1$;\cr \cr {1 \over 2}\left(2c_{34} + a_1c_4 - b_3c_2 + b_2c_3 - a_4c_1 + ca_{14} + ba_{24} - cc_{13} - bc_{23} - a_{44} - b_{33}\right)&\cr &${\bf i}=2$;\cr\ \ + {1 \over 4}\left(a_2b_4 + a_3b_1 - a_4b_2 - a_1b_3 + aa_1b_1 + ca_1b_2 - cc_1c_2 + ba_2c_1 - ba_1c_2 - a(c_1)^2\right),&\cr \cr -{1 \over 4}(2a_{14} - 2c_{13} - a_2b_1 + c_1c_2),&${\bf i}=3$;\cr \cr -{1 \over 4}\left(2a_{24} - 2c_{23} - a_1c_2 +a_2c_1 - a_2b_2 + (c_2)^2\right),&${\bf i}=4$;\cr}\cr
\noalign{\vskip 12pt}
R^{\bf i}{}_{434} = \cases{-{1 \over 2}\left(2c_{34} + a_1c_4 - a_4c_1 + b_2c_3 - b_3c_2 - ac_{14} - cc_{24} + ab_{13} + cb_{23} - a_{44} - b_{33}\right)&\cr &${\bf i}=1$;\cr\ \ - {1 \over 4}\left(a_3b_1 - a_1b_3 + a_2b_4 - a_4b_2 - cc_1c_2 + ab_1c_2 + ca_1b_2 + ba_2b_2 -ac_1b_2 - b(c_2)^2\right),&\cr \cr -{1 \over 4}\left(2bb_{23} + 2cb_{13} - 2bc_{24} - 2cc_{14} + cb_1c_2 + bc_1c_2 - cb_2c_1 - bb_1a_2 - cb_1a_1 + c(c_1)^2\right),&${\bf i}=2$;\cr \cr {1 \over 4}\left(2b_{13} - 2c_{14} + b_1c_2 - b_2c_1 - a_1b_1 + (c_1)^2\right),&${\bf i}=3$;\cr \cr {1 \over 4}(2b_{23} - 2c_{24} - a_2b_1 + c_1c_2),&${\bf i}=4$.\cr}\cr}$$
\vskip 24pt
\noindent {\bf A1.4 Fully Covariant Riemann Curvature}\hfil\break
Straightforward computation from A1.3 yields:
$$\displaylines{R_{12{\bf jk}} = 0\cr
R_{1313} = {1 \over 2}a_{11} \hskip .5in R_{1314} = {1 \over 2}c_{11} \hskip .5in R_{1323} = {1 \over 2}a_{12} \hskip .5in R_{1324} = {1 \over 2}c_{12}\cr
\noalign{\vskip 12pt}
R_{1334} = -{1 \over 4}(2a_{14} - 2c_{13} - a_2b_1 + c_1c_2)\cr
\noalign{\vskip 12pt}
R_{1414} = {1 \over 2}b_{11} \hskip .75in R_{1423} = {1 \over 2}c_{12} \hskip .75in R_{1424} = {1 \over 2}b_{12}\cr
\noalign{\vskip 12pt}
R_{1434} = {1 \over 4}\left(2b_{13} - 2c_{14} - a_1b_1 + b_1c_2 - b_2c_1 + (c_1)^2\right)\cr
\noalign{\vskip 12pt}
R_{2323} = {1 \over 2}a_{22} \hskip 1.25in R_{2324} = {1 \over 2}c_{22}\cr
\noalign{\vskip 12pt}
R_{2334} = -{1 \over 4}\left(2a_{24} - 2c_{23} - a_2b_2 + a_2c_1 - a_1c_2 + (c_2)^2\right)\cr
\noalign{\vskip 12pt}
R_{2424} = {1 \over 2}b_{22} \hskip 1in R_{2434} = {1 \over 4}(2b_{23} - 2c_{24} -a_2b_1 + c_1c_2)\cr}$$
$$\eqalign{R_{3434} &= -{1 \over 2}\left(2c_{34} + a_1c_4 - a_4c_1 + b_2c_3 - b_3c_2 - cc_1c_2 - a_{44} - b_{33}\right)\cr 
&-{1 \over 4}\left(a_3b_1 - a_1b_3 + a_2b_4 - a_4b_2 + aa_1b_1 + ba_2b_2 + ca_1b_2 + ca_2b_1 - a(c_1)^2 - b(c_2)^2\right)\cr}$$
\vskip 24pt
\noindent {\bf A1.5 The Ricci Tensor}\hfil\break
The definition of the Ricci tensor in [27] is:
$$R_{ab} := R^c{}_{bac} = R_{ac}{}^c_b = -R^c_{acb} = R^c{}_{abc}.\eqno({\rm A}1.12)$$
The Ricci curvature is defined in [28] as the negative of (A1.12) so we must modify the equations [28](4.6.20--23) and [28](4.6.25) by removing a minus sign so as to preserve the definitions of $\Phi_{ABA'B'}$ and $\Lambda$, see Appendix Two.
$$\displaylines{R_{11} = 0 \hskip .75in R_{12} = 0 \hskip .75in R_{13} = {1 \over 2}(a_{11} + c_{12}) \hskip .75in R_{14} = {1 \over 2}(b_{12} + c_{11})\cr
\noalign{\vskip 6pt}
R_{22} = 0 \hskip 1in R_{23} = {1 \over 2}(a_{12} + c_{22}) \hskip 1in R_{24} = {1 \over 2}(b_{22} + c_{12})\cr}$$
\vskip -12pt
$$\eqalign{R_{33} &= {1 \over 2}\left(2ca_{12} - 2a_{24} + 2c_{23} + aa_{11} + ba_{22} + a_2b_2 + a_1c_2 - a_2c_1 - (c_2)^2\right)\cr
R_{34} &= {1 \over 2}(2cc_{12} + ac_{11} + bc_{22} + a_{14} + b_{23} - c_{13} - c_{24} - a_2b_1 + c_1c_2)\cr
R_{44} &= {1 \over 2}\left(2c_{14} - 2b_{13} + 2cb_{12} + ab_{11} + bb_{22} + a_1b_1 + b_2c_1 - b_1c_2 - (c_1)^2\right)\cr}$$
\vskip 24pt
\noindent {\bf A1.6 The Scalar Curvature}\hfil\break
$$S = a_{11} + b_{22} + 2c_{12}.$$
\vskip 24pt
\noindent{\bf A1.7 The Ricci Endomorphism}\hfil\break
$$R^{\bf i}{}_1 = \cases{{1 \over 2}(a_{11} + c_{12}),&${\bf i}=1$;\cr \cr
{1 \over 2}(b_{12} + c_{11}),&${\bf i}=2$;\cr \cr 0,& otherwise;\cr}$$
$$R^{\bf i}{}_2 = \cases{{1 \over 2}(a_{12} + c_{22}),&${\bf i}=1$;\cr \cr
{1 \over 2}(b_{22} + c_{12}),&${\bf i}=2$;\cr \cr 0,& otherwise;\cr}$$
$$R^{\bf i}{}_3 = \cases{{1 \over 2}\left(2c_{23} - 2a_{24} + ba_{22} + ca_{12} - ac_{12} - cc_{22} + a_2b_2 + a_1c_2 - a_2c_1 - (c_2)^2\right) =: \zeta,&${\bf i}=1$;\cr \cr
{1 \over 2}(ac_{11} + cc_{12} - ca_{11} - ba_{12} + a_{14} + b_{23} - c_{13} - c_{24} - a_2b_1 + c_1c_2) =: \eta,&${\bf i}=2$;\cr \cr
{1 \over 2}(a_{11} + c_{12}),&${\bf i}=3$;\cr \cr
{1 \over 2}(a_{12} + c_{22}),&${\bf i}=4$;\cr}$$
$$R^{\bf i}{}_4 = \cases{{1 \over 2}(bc_{22} + cc_{12} - ab_{12} - cb_{22} + a_{14} + b_{23} - c_{13} - c_{24} - a_2b_1 + c_1c_2) =: \Xi,&${\bf i}=1$;\cr \cr
{1 \over 2}\left(2c_{14} - 2b_{13} + ab_{11} + cb_{12} - cc_{11} - bc_{12} + a_1b_1 + b_2c_1 - b_1c_2 - (c_1)^2\right) =: \Upsilon,&${\bf i}=2$;\cr \cr
{1 \over 2}(b_{12} + c_{11}),&${\bf i}=3$;\cr \cr
{1 \over 2}(b_{22} + c_{12}),&${\bf i}=4$;\cr}$$
\vskip 24pt
\noindent {\bf A1.8 The Einstein Endomorphism}\hfil\break
The diagonal elements of $E^a{}_b := R^a{}_b - (S/4)\delta^a{}_b$ are:
$$E^1{}_1 = E^3{}_3 = -E^2{}_2 = -E^4{}_4 = {a_{11} - b_{22} \over 4}.$$
The off-diagonal elements are as in A1.6. 
Putting
$$\theta := {a_{11} - b_{22} \over 4} \hskip .75in \mu := {b_{12} + c_{11} \over 2} \hskip .75in \nu := {a_{12} + c_{22} \over 2}$$
one computes
$$\eqalign{E(\partial_1) &= \theta\partial_1 + \mu\partial_2\cr
E(\partial_2) &= \nu\partial_1 - \theta\partial_2\cr
E(\partial_3) &= \zeta\partial_1 + \eta\partial_2 + \theta\partial_3 + \nu\partial_4\cr
E(\partial_4) &= \Xi\partial_1 + \Upsilon\partial_2 + \mu\partial_3 - \theta\partial_4\cr}$$
where $\zeta$, $\eta$, $\Xi$ and $\Upsilon$ are defined in A1.7. A1.4 was reported by Ghanam \& Thompson [13], though with a typographical error. Matsushita [23] reported A1.14--6 and the covariant form of the Einstein endomorphism. Chaichi et al. [5] computed curvature properties under the perhaps ad hoc assumption $c=0$. D\'{\i}az-Ramos et al. [8] also reported A1.4--6. Generally, these authors employed distinct curvature conventions to us; we have been motivated by choices which maintain a close correspondence with the conventions of [28] so as to facilitate the employment of spinors.
\vskip 24pt
\noindent {\bf A1.9 Conformal Curvature}\hfil\break
With the conventions of [27], the Weyl conformal curvature is given by:
$$\eqalignno{R_{abcd} &= C_{abcd} - {S \over 6}(g_{ad}g_{bc} - g_{ac}g_{bd}) + {1 \over 2}(g_{ad}R_{bc} - g_{ac}R_{bd} + g_{bc}R_{ad} - g_{bd}R_{ac})\cr
&= C_{abcd} + {S \over 12}(g_{ad}g_{bc} - g_{ac}g_{bd}) + {1 \over 2}(g_{ad}E_{bc} - g_{ac}E_{bd} + g_{bc}E_{ad} - g_{bd}E_{ac}),&({\rm A}1.13)\cr}$$
and the Weyl curvature components are:
$$\displaylines{C_{12{\bf jk}} = 0\qquad\hbox{for}\qquad ({\bf j},{\bf k}) \not= (3,4); \hskip 1in C_{1234} = {S \over 12}\cr
\noalign{\vskip 12pt}
C_{1313} = {1 \over 6}(a_{11} - c_{12} + b_{22}) \qquad C_{1314} = {1 \over 4}(c_{11} - b_{12}) \qquad C_{1323} = {1 \over 4}(a_{12} - c_{22}) \qquad C_{1324} = {1 \over 2}c_{12}\cr
\noalign{\vskip 12pt}
C_{1334} = -{1 \over 12}(3a_{14} - 3c_{13} - 5cc_{12} - 3bc_{22} + 3c_{24} - 3b_{23} - ca_{11} + 3ab_{12} + 2cb_{22})\cr
\noalign{\vskip 12pt}
C_{1414} = {1 \over 2}b_{11} \qquad C_{1423} = -{1 \over 12}(a_{11} - 4c_{12} + b_{22}) \qquad C_{1424} = {1 \over 4}(b_{12} - c_{11})\cr
\noalign{\vskip 12pt}
C_{1434} = {1 \over 12}(b\,a_{11} + 3ab_{11} + 3cb_{12} + bb_{22} - bc_{12} - 3cc_{11})\cr
\noalign{\vskip 12pt}
C_{2323} = {1 \over 2}a_{22} \hskip 1.25in C_{2324} = -{1 \over 4}(a_{12} - c_{22})\cr       
\noalign{\vskip 12pt}
C_{2334} = -{1 \over 12}(aa_{11} + 3ca_{12} + 3ba_{22} - 3cc_{22} - ac_{12} + ab_{22})\cr
\noalign{\vskip 12pt}
C_{2424} = {1 \over 6}(a_{11} - c_{12} + b_{22})\cr
\noalign{\vskip 12pt}
C_{2434} = {1 \over 12}(2ca_{11} + 3ba_{12} - 3a_{14} + 3b_{23} - cb_{22} - 3c_{24} - 3ac_{11} - 5cc_{12} + 3c_{13})\cr}$$
$$\eqalign{C_{3434} &= {1 \over 12}\bigl(3ba_1c_2 + 3cb_1a_2 + 6bca_{12}  + baa_{11} - 3ca_1b_2  - 4abc_{12} - 6cac_{11} - 3bc_1a_2 + abb_{22}\cr
&\qquad  + 6acb_{12} - 6cbc_{22} + 3ac_1b_2 + 3a^2b_{11} + 6c_1a_4 - 3ab_1c_2 - 6a_1c_4 + 3a_1b_3 + 6c_2b_3\cr
& \qquad - 3a_2b_4 - 3b_1a_3 - 6b_2c_3 + 3b_2a_4 - 12c_{34} + 6a_{44} + 6b_{33} + 6ac_{14} - 6ab_{13} - 8c^2c_{12}\cr
& \qquad - 6ca_{14} + 6cc_{13} + 6cc_{24} - 6cb_{23} + 3b^2a_{22} - 6ba_{24} + 6bc_{23} + 2c^2a_{11} + 2c^2b_{22}\bigr)\cr}$$
\vskip 24pt
\noindent {\bf A1.10 The Curvature Endomorphism}
\vskip 12pt
Define the {\sl curvature endomorphism\/} of the space $\Lambda^2(T_pM)$ of bivectors as the tensor contraction:
$${\cal R}(F) := {1 \over 2}R^{ab}{}_{cd}F^{cd} = {1 \over 2}C^{ab}{}_{cd}F^{cd} + {S \over 12}F^{ab} + {1 \over 2}(E^a{}_cF^{bc} - E^b{}_cF^{ac}),\eqno({\rm A}1.14)$$
where the second expression follows from (A1.13). If $R^{ab}{}_{cd}$ has components $R^{\bf ij}{}_{\bf kl}$ with respect to a frame $\{e_1,\ldots,e_4\}$, then $\cal R$ has matrix $\left(R^{\bf ij}{}_{\bf kl}\right)$ with respect to the induced basis $\{\,e^a_{\bf i} \wedge e^b_{\bf j}: {\bf i} < {\bf j}\,\}$ for $\Lambda^2(T_pM)$. The trace of this endomorphism is
$$\tr({\cal R}) = \sum_{{\bf i} < {\bf j}}\,R^{\bf ij}{}_{\bf ij} = {1 \over 2}\sum_{\bf ij}\,R^{\bf ij}{}_{\bf ij} = {S \over 2}.\eqno({\rm A}1.15)$$
This definition ensures that $\cal R$ is the identity on ${\bf S}^4$. 

$\Lambda^2(T_pM)$ equipped with the induced scalar product is isomorphic to ${\bf R}^{2,4}$ and the Hodge star operator $*$ induces, via its eigenspaces, the orthogonal decomposition $\Lambda^2 = \Lambda^2_+ \oplus \Lambda^2_-$, under which
$${\cal R} = \pmatrix{{\cal W}^+&{\cal Z}\cr {^*\!{\cal Z}}&{\cal W}^-\cr} + {S \over 12}1_6,\eqno({\rm A}1.16)$$
with $\cal R$ self adjoint. ${\cal W}^+$ and ${\cal W}^-$ are induced by the self-dual and anti-self-dual parts of the Weyl conformal tensor respectively and, with $\Lambda^2_+(T_pM) \cong \Lambda^2_-(T_pM) \cong {\bf R}^{1,2}$, ${\cal W}^+$ and ${\cal W}^-$ are self adjoint with respect to the induced scalar products, while $\cal Z$ is induced by the Einstein endomorphism $E^a{}_b$, and $^*\!{\cal Z}$ is the adjoint of ${\cal Z} \in \Hom(\Lambda^2_-,\Lambda^2_+)$.

To compute the curvature endomorphism, one requires a basis of $\Lambda^2$ compatible with the decomposition $\Lambda^2 = \Lambda^2_+ \oplus \Lambda^2_-$. If $\{e_1, e_2, e_3, e_4\}$ is an oriented $\Psi$-ON frame, then the following is a $\Psi$-ON frame of $\Lambda^2 = \Lambda^2_+ \oplus \Lambda^2_-$:
$$\displaylines{\hskip 1in s^+_1 := {e_1 \wedge e_2 + e_3 \wedge e_4 \over \sqrt2} \hskip 1.25in s^-_1 := {e_1 \wedge e_2 - e_3 \wedge e_4 \over \sqrt2}\hfill\cr
\hskip 1in s^+_2 := {e_1 \wedge e_3 + e_2 \wedge e_4 \over \sqrt2} \hskip 1.25in s^-_2 := {e_1 \wedge e_3 - e_2 \wedge e_4 \over \sqrt2}\hfill\llap({\rm A}1.17)\cr
\hskip 1in s^+_3 := {e_1 \wedge e_4 - e_2 \wedge e_3 \over \sqrt2} \hskip 1.25in s^-_3 := {e_1 \wedge e_4 + e_2 \wedge e_3 \over \sqrt2}.\hfill\cr}$$
A simple choice of $\Psi$-ON frame constructed from Walker coordinates $(u,v,x,y)$ is provided by
$$\vcenter{\openup1\jot \halign{$\hfil#$&&${}#\hfil$&\qquad$\hfil#$\cr
e_1 &:= \displaystyle{1 \over 2}(1-a)\partial_1 + \partial_3 & e_3 &:= \displaystyle-{1 \over 2}(1+a)\partial_1 + \partial_3\cr
e_2 &:= \displaystyle-c\partial_1 + {1 \over 2}(1-b)\partial_2 + \partial_4 & e_4 &:= \displaystyle-c\partial_1 - {1\over2}(1+b)\partial_2 + \partial_4.\cr}}\eqno({\rm A}1.18)$$
The matrix relating these frames has determinant one so $[e_1,e_2,e_3,e_4] = [\partial_1,\partial_2,\partial_3,\partial_4]$, i.e., $\{e_1,e_2,e_3,e_4\}$ is a $\Psi$-ON frame with the canonical orientation (see \S 1) and thus suitable for employment in (A1.17).

If one takes the Walker coordinates $(v,u,y,x)$ derived from the symmetry (A1.7), the $\Psi$-ON frame obtained from them via (A1.18) also possesses the canonical orientation but opposite ${\bf SO^\bfplus}$-orientation (i.e., opposite `time' and `space' orientations).

With the choice (A1.18), one obtains
$$\displaylines{\hskip 1.2in s^+_1 = {1+ab \over 2}\,\partial_1 \wedge \partial_2 + 2c\,\partial_1 \wedge \partial_3 - a\,\partial_1 \wedge \partial_4 + b\,\partial_2 \wedge \partial_3 + 2\,\partial_3 \wedge \partial_4,\hfill\cr
\hskip 1.2in s^+_2 = c\,\partial_1 \wedge \partial_2 + \partial_1 \wedge \partial_3 + \partial_2 \wedge \partial_4,\hfill\llap({\rm A}1.19)\cr
\hskip 1.2in s^+_3 = {ab-1 \over 2}\,\partial_1 \wedge \partial_2 + 2c\,\partial_1 \wedge \partial_3 - a\,\partial_1 \wedge \partial_4 + b\,\partial_2 \wedge \partial_3 + 2\,\partial_3 \wedge \partial_4,\hfill\cr}$$
and
$$\displaylines{\hskip 2in s^-_1 = -{a+b \over 2}\,\partial_1 \wedge \partial_2 + \partial_1 \wedge \partial_4 - \partial_2 \wedge \partial_3,\hfill\cr
\hskip 2in s^-_2 = - c\,\partial_1 \wedge \partial_2 + \partial_1 \wedge \partial_3 - \partial_2 \wedge \partial_4,\hfill\llap({\rm A}1.20)\cr
\hskip 2in s^-_3 = {a-b \over 2}\,\partial_1 \wedge \partial_2 + \partial_1 \wedge \partial_4 + \partial_2 \wedge \partial_3.\hfill\cr}$$
The matrix representations of the Weyl curvature endomorphisms have been reported in [8] and [10]. Putting
$$\displaylines{{\cal P} := a_{11} + b_{22} - 4c_{12} \hskip 1in {\cal Q} := a_{22} + b_{11} \hskip 1in {\cal T} := a_{12} - c_{22}\cr
\hfill {\cal X} := b_{12} - c_{11} \hskip 1.25in {\cal Y} := a_{22} - b_{11}\hfill\llap({\rm A}1.21)\cr}$$
then the matrix representation of ${\cal W}^-$ is
$$^-{\bf W} := -{1 \over 12}\pmatrix{-({\cal P} + 3{\cal Q})&3({\cal T}+{\cal X})&3{\cal Y}\cr
-3({\cal T}+{\cal X})&2{\cal P}&3({\cal T}-{\cal X})\cr
-3{\cal Y}&3({\cal T}-{\cal X})&-({\cal P}-3{\cal Q})\cr}.\eqno({\rm A}1.22)$$

With
$$\eqalignno{A &= 6ab_{13} - 6bc_{23} - 12cc_{13} - 12b_{33} + 12c_{34}\cr
&\qquad - 6ac_{14} + 6ba_{24} + 12ca_{14} + 12c_{34} - 12a_{44}\cr
&\qquad - 3a(-bc_{12} - 2c_2b_1 + ab_{11}) - 3b(ba_{22} - 2a_2c_1 - ac_{12}) - 6c(ba_{12} + a_1b_2 + a_2b_1 - ac_{11})\cr
&\qquad + 6(-bc_{23} - c_2b_3 + a_3b_1 + ab_{13}) + 6(ba_{24} + a_2b_4 - a_4c_1 - ac_{14})\cr
&\qquad + \left[-6ac_1b_2 + 6acc_{11} - 6bc_2a_1 - 6bca_{12} + 12ca_1b_2 - 12c^2a_{11} + 12b_2c_3 - 12cc_{13} + 12a_1c_4 + 12ca_{14}\right]\cr
&\qquad - 6a_4c_1 - 6a_4b_2 - 6a_1b_3 - 6b_3c_2 - a_{11} - b_{22} - 2c_{12},&({\rm A}1.23)\cr}$$
where we note that for use in 3.7 the right-hand side is expressed intentionally in a redundant form in that we have not grouped all like terms together, in particular the second term in the third bracketed quantity of the third line ($-6ca_1b_2$) cancels half of the fifth term in the fifth line ($12ca_1b_2$) to yield the term $6ca_1b_2$, and
$$B :=  2(a_{14} - b_{23} - c_{13} + c_{24}) - 2ca_{11} - ba_{12} + ab_{12} + ac_{11} - bc_{22}- 2cc_{12} ,\eqno({\rm A}1.24)$$
the matrix representation of ${\cal W}^+$ is
$$^+{\bf W} := -{1 \over 12}\pmatrix{A&3B&A+S\cr -3B&2S&-3B\cr -(A+S)&-3B&-(A+2S)\cr}.\eqno({\rm A}1.25)$$
One computes
$$\det({\cal W}^+ - \lambda1_3) = (S/6 + \lambda)(\lambda - S/12)^2,\eqno({\rm A}1.26)$$
whence the eigenvalues of ${\cal W}^+$ are
$$-{S \over 6},\ {S \over 12},\ {S \over 12}.\eqno({\rm A}1.27)$$
D\'{\i}az-Ramos et al. [8] showed that ${\cal W}^+$ possesses only certain possible Jordan canonical forms, see 2.6, indicating that generic four-dimensional Walker geometry manifests certain self duality properties, as explained by our spinor analysis of Walker geometry.

The Einstein endomorphism (i.e., traceless Ricci tensor) determines
$${\cal E} := \pmatrix{{\bf 0}&{\cal Z}\cr {^*\!{\cal Z}}&{\bf 0}\cr}.\eqno({\rm A}1.28)$$
From (A1.14)
$${\cal E}(F) = {1 \over 2}(E^a{}_cF^{bc} - E^b{}_cF^{ac}) = -{1 \over 2}(E^a{}_cF^{cb} - E^b{}_cF^{ca}).\eqno({\rm A}1.29)$$
In particular,
$${\cal E}(X \wedge Y) = -{1 \over 2}(X \wedge E(Y) + E(X) \wedge Y).\eqno({\rm A}1.30)$$
Using A1.8 and (A1.30), Davidov \& Mu\v skarov [10] computed the matrix representation of ${\cal Z}$ with respect to the bases  (A1.19--20):
$${\cal Z} = -{1 \over 2}\pmatrix{\Upsilon+\zeta+c(\nu-\mu)&\eta+\Xi-2\theta c&\Upsilon-\zeta-c(\nu+\mu)\cr
\mu-\nu&2\theta&\mu+\nu\cr
-\bigl(\Upsilon+\zeta+c(\nu-\mu)\bigr)&-\bigl(\eta+\Xi-2\theta c\bigr)&-\bigl(\Upsilon-\zeta-c(\nu+\mu)\bigr)\cr},\eqno({\rm A}1.31)$$
which completes the description of the curvature endomorphism.

Both D\'{\i}az-Ramos et al. [8] and Davidov \& Mu\v skarov [10] used these results to characterize the vanishing of the ASD Weyl curvature of a four-dimensional Walker geometry; the latter authors also obtained characterizations of some other curvature conditions while the former authors studied the Osserman condition on the Jacobi operator.
\vskip 24pt
\noindent {\section Appendix Two: Spinors for Four-Dimensional Neutral Metrics}
\vskip 12pt
The two-component spinor formalism for ${\bf R}^{2,2}$ is more directly analogous to that for ${\bf C}^4$ equipped with the standard {\bf C}-bilinear scalar product than that for ${\bf R}^{1,3}$. That said, it is mostly straightforward to adapt the results in [28] to the context of neutral signature, though one must be aware of a few features peculiar to neutral geometry.

The two-component spinor formalism is based on an isomorphism, via Clifford algebras, of ${\bf R}^4$ with $S \otimes S'$, where $S \cong S' \cong {\bf R}^2$, but $S$ and $S'$ are independent spaces. Each of $S$ and $S'$ is more appropriately viewed as isomorphic to the symplectic plane ${\bf R}^2_{\rm sp}$, with symplectic forms $\epsilon$ and $\epsilon'$. Objects constructed from the tensor algebra of $S'$ are indicated by abstract indices bearing a prime, in which case primes on the symbol denoting the object itself are dropped, whence $\epsilon_{A'B'}$ for $\epsilon'$. The actual isomorphism of interest is then ${\bf R}^{2,2} = ({\bf R}^4,\eta) \cong (S \otimes S',\epsilon \otimes \epsilon')$. See [19] for a brief sketch of this isomorphism and some basic spinor algebra and geometry. In particular, the isomorphism ${\bf R}^{2,2} \cong S \otimes S'$ may be taken to be
$$v = (v^1,v^2,v^3,v^4)\ \leftrightarrow\ {1 \over \sqrt2}\pmatrix{v^1+v^3&v^4-v^2\cr v^4+v^2&v^1-v^3\cr} =: \left(v^{\bf AA'}\right),\eqno({\rm A}2.1)$$
whence $s(v,v) = 2\det\left(v^{\bf AA'}\right)$; in particular, $v$ is null iff $\left(v^{\bf AA'}\right)$ is singular, equivalently $v^{AA'}$ is decomposable as an element of $S \otimes S'$.

When employing results from [28], the main fact to bear in mind is that there is no natural identification between $S$ and $S'$ (which for neutral signature are real linear spaces whereas in the case of Lorentz signature they are complex conjugate (linear) spaces of each other). In this appendix we will merely record a few results which we require in the main body of the paper.

For ${\bf R}^{2,2}$, (A1.3) indicates that $V_{abcd} = e_{abcd}$ and, by (A1.6) $\xi_g(*U) = *\xi_g(U)$ for any multivector $U$, i.e., the Hodge star operators on multivectors and forms coincide under the identification of multivectors and forms via the metric. The volume form of the standard orientation of ${\bf R}^{2,2}$ is
$$e_{abcd} = \epsilon_{AC}\epsilon_{BD}\epsilon_{A'D'}\epsilon_{B'C'} - \epsilon_{AD}\epsilon_{BC}\epsilon_{A'C'}\epsilon_{B'D'}.\eqno({\rm A}2.2)$$
Representing an element of $\Lambda^2\left({\bf R}^{2,2}\right)$ as a skew tensor $F^{ab}$, the decomposition into SD and ASD summands is
$$F^{ab} = \epsilon^{AB}\psi^{A'B'} + \phi^{AB}\epsilon^{A'B'},\eqno({\rm A}2.3)$$
where $\psi^{A'B'} \in S' \odot S'$ and $\phi^{AB} \in S \odot S$. Note that $*^2 = 1$.

The spinorial representation of the curvature may be obtained exactly as in [28], \S 4.6, the only difference being that all curvature spinors are real objects whence the SD and ASD Weyl spinors are independent objects: $\epsilon_{AB}\epsilon_{CD}\tilde\Psi_{A'B'C'D'}$ is the spinorial representation of the SD Weyl tensor and $\Psi_{ABCD}\epsilon_{A'B'}\epsilon_{C'D'}$ that of the ASD Weyl curvature tensor. The fully covariant Riemann tensor is given by
$$\eqalignno{R_{abcd} &= \epsilon_{AB}\epsilon_{CD}\tilde\Psi_{A'B'C'D'} + \Psi_{ABCD}\epsilon_{A'B'}\epsilon_{C'D'}\cr
&\qquad + \Phi_{ABC'D'}\epsilon_{A'B'}\epsilon_{CD} + \Phi_{CDA'B'}\epsilon_{AB}\epsilon_{C'D'}&({\rm A}2.4)\cr
&\qquad + 2\Lambda(\epsilon_{AC}\epsilon_{A'C'}\epsilon_{BD}\epsilon_{B'D'} - \epsilon_{AD}\epsilon_{A'D'}\epsilon_{BC}\epsilon_{B'C'})\cr}$$
where the Weyl spinors $\tilde\Psi_{A'B'C'D'}$ and $\Psi_{ABCD}$ are fully symmetric while the Ricci spinor satisfies $\Phi_{ABA'B'} = \Phi_{(AB)(A'B')}$.

Because our definition A1.5 of the Ricci tensor is the negative of that employed in [28], we obtain:
$$R_{ab} = 2\Phi_{ABA'B'} - 6\Lambda\epsilon_{AB}\epsilon_{A'B'},\eqno({\rm A}2.5)$$
whence
$$S = -24\Lambda \hskip 1.25in 2\Phi_{ab} := 2\Phi_{ABA'B'} = R_{ab} - {S \over 4}g_{ab} = E_{ab},\eqno({\rm A}2.6)$$
where $E_{ab}$ is the fully covariant version of the Einstein endomorphism, i.e., the trace-free Ricci tensor (not to be confused with the `Einstein tensor' $G_{ab}$ of [28]).

As in [28], \S 4.9,
$$\triangle_{ab} := 2\nabla_{[a}\nabla_{b]} = \epsilon_{A'B'}\,\Square_{AB} + \epsilon_{AB}\,\Square_{A'B'},\eqno({\rm A}2.7)$$
where 
$$\Square_{AB} := \nabla_{X'(A}\nabla_{B)}{}^{X'} \hskip 1in \Square_{A'B'} := \nabla_{X(A'}\nabla_{B')}{}^X.\eqno({\rm A}2.8)$$
The {\sl spinor Ricci identities\/} for arbitrary spinors $\kappa_A$ and $\tau_{A'}$ are:
$$\displaylines{\Square_{AB}\kappa_C = \Psi_{ABCE}\kappa^E - \Lambda(\kappa_A\epsilon_{BC} + \epsilon_{AC}\kappa_B) \hskip 1in \Square_{AB}\tau_{C'} = \Phi_{ABC'E'}\tau^{E'}\cr
\hfill\llap({\rm A}2.9)\cr
\Square_{A'B'}\tau_{C'} = \tilde\Psi_{A'B'C'E'}\tau^{E'} - \Lambda(\tau_{A'}\epsilon_{B'C'} + \epsilon_{A'C'}\tau_{B'}) \hskip .7in \Square_{A'B'}\kappa_C = \Phi_{A'B'CE}\kappa^E\cr}$$
The Bianchi equation may be written, see [28], \S 4.10,
$$\displaylines{\nabla^A_{B'}\Psi_{ABCD} = \nabla^{A'}_{(B}\Phi_{CD)A'B'} \hskip 1in \nabla^{A'}_B\tilde\Psi_{A'B'C'D'} = \nabla^A_{(B'}\Phi_{C'D')AB}\cr
\hfill\llap({\rm A}2.10)\cr
\nabla^{CA'}\Phi_{CDA'B'} = -3\nabla_{DB'}\Lambda.\cr}$$
\vskip 24pt
\baselineskip=12pt
\noindent {\bf BIBLIOGRAPHY}
\vskip 12pt
\newcount\q \q=0
\def\nref {\global\advance\q by1 \item{[\the\q]}}
\frenchspacing
\baselineskip=12pt
\nref B\'erard Bergery, L., Ikemakhen, A.: Sur l'holonomie des vari\'et\'es pseudo-riemanniennes de signature $(n,n)$. Bull. Soc. Math. France {\bf 125}, 93--114 (1997).
\vskip 1pt
\nref Bla\v zi\'c, N., Bokan N., Raki\'c, Z.: Osserman Pseudo-Riemannian Manifolds of Signature $(2,2)$. J. Austral. Math. Soc. {\bf 71}, 367--395 (2001).
\vskip 1pt
\nref Boyer, C.P., Finley III, J.D., Pleba\~nski, J.F.: Complex General Relativity, $\cal H$ and ${\cal HH}$ Spaces-A Survey of One Approach. In: Held, A. (ed.) General Relativity and Gravitation: One Hundred Years After the Birth of Albert Einstein, Vol. 2. New York: Plenum Press, 1980, pp. 241--281.
\vskip 1pt
\nref Brozos-V\'azquez, M., Garc\'{\i}a-R\'{\i}o, E., V\'azquez-Lorenzo, R.: Conformally Osserman four-dimensional manifolds whose conformal Jacobi operators have complex eigenvalues. Proc. R. Soc. London A {\bf 462}, 1425--1441 (2006).
\vskip 1pt
\nref Chaichi,  M., Garc\'\i a-R\'\i o,  E., Matsushita, Y.: Curvature properties of four-dimensional Walker metrics. Classical Quantum Gravity {\bf 22}, 559--577 (2005).
\vskip 1pt
\nref Cruceanu,  V., Gadea, P.M., Mu\~noz Masqu\'e, J.: Para-Hermitian and Para-K\"ahler Manifolds. Quaderni Inst. Mat. Univ. Messina {\bf 1}, 1--72 (1995).
\vskip 1pt
\nref Cruceanu, V., Fortuny, P., Gadea, P.M.: A survey on paracomplex geometry. Rocky Mountain J. Math. {\bf 26}, 83--115 (1996).
\vskip 1pt
\nref D\'{\i}az-Ramos,  J.C., Garc\'\i a-R\'\i o, E., V\'azquez-Lorenzo, R.: Four-dimensional Osserman metrics with nondiagonalizable Jacobi operators. J. Geom. Anal. {\bf 16}, 39--52 (2006).
\vskip 1pt
\nref Davidov,  J., D\'{\i}az-Ramos, J.C., Garc\'{\i}a-R\'{\i}o, E., Matsushita, Y., Mu\v skarov, O., V\'azquez-Lorenzo, R.: Almost K\"ahler Walker four-manifolds. J. Geom. Phys. {\bf 57}, 1075--1088 (2007).
\vskip 1pt
\nref Davidov, J., Mu\v skarov, O.: Self-dual Walker metrics with two-step nilpotent Ricci operator. J. Geom. Phys. {\bf 57}, 157--165 (2006).
\vskip 1pt
\nref Dunajski, M.: Anti-self-dual four-manifolds with a parallel spinor. Proc. R. Soc. London A {\bf 458}, 1205--1222 (2002).
\vskip 1pt
\nref Dunajski, M.,West, S.: Anti-self-dual conformal structures in neutral signature.\hfil\break www.arxiv.org.math.DG/0610280
\vskip 1pt
\nref Ghanam, R., Thompson, G.: The holonomy Lie algebras of neutral metrics in dimension four. J. Math. Phys. {\bf 42}, 2266--2284 (2001).
\vskip 1pt
\nref Jensen, G.R., Rigoli, M.: Neutral Surfaces in Neutral Four-Spaces. Le Matematiche {\bf XLV}, 407--443 (1990).
\vskip 1pt
\nref Kamada, H.: Self-duality of neutral metrics on four-dimensional manifolds. In: The Third Pacific Rim Geometry Conference (Seoul, 1996), Monographs Geom. Topology {\bf 25}. Cambridge, MA: International Press, 1998, pp. 79--98.
\vskip 1pt
\nref Kobayashi, S., Nomizu, K.: Foundations of Differential Geometry, Volume I. New York:\hfil\break Wiley-Interscience, 1963.
\vskip 1pt
\nref Law, P.R.: Neutral Einstein metrics in four dimensions. J. Math. Phys. {\bf 32}, 3039--3042 (1991).
\vskip 1pt
\nref Law, P.R.: Neutral Geometry and the Gauss-Bonnet Theorem for Two-Dimensional Pseudo-Riemannian Manifolds. Rocky Mountain J. Math. {\bf 22}, 1365--1383 (1992).
\vskip 1pt
\nref Law, P.R.: Classification of the Weyl curvature spinors of neutral metrics in four dimensions. J. Geom. Phys. {\bf 56}, 2093--2108 (2006).
\vskip 1pt
\nref Libermann, P.: Sur les structures presque paracomplexes. C. R. Acad. Sci. Paris S\'er. I Math. {\bf 234}, 2517--2519 (1952).
\vskip 1pt
\nref Libermann, P.: Sur le probl\`eme d'\'equivalence de certaines structures infinit\'esimales. Ann. Mat. Pura Appl. (4) {\bf 36}, 27--120 (1954).
\vskip 1pt
\nref Matsushita, Y.: Four-Dimensional Walker metrics and symplectic structures. J. Geom. Phys. {\bf 52}, 89--99 (2004).
\vskip 1pt
\nref Matsushita, Y.: Walker 4-manifolds with proper almost complex structures. J. Geom. Phys. {\bf 55}, 385--398 (2005).
\vskip 1pt
\nref Matsushita, Y.: The existence of indefinite metrics of signature $(++--)$ and two kinds of almost complex structures in dimension four, in: Contemporary Aspects of Complex Analysis, Differential Geometry and Mathematical Physics (Proceedings of the 7th International Workshop on Complex Structures and Vector Fields, 2004, Plovdiv, Bulgaria), 210--226. Hackensack, NJ: World Scientific Publishing Co., 2005.
\vskip 1pt
\nref Matsushita, Y., Haze, S., Law, P.R.: Almost K\"ahler-Einstein Structures on 8-Dimensional Walker Manifolds. Monatsh. Math. {\bf 150}, 41--48 (2007).
\vskip 1pt
\nref Matsushita, Y., Law, P.R.: Hitchin-Thorpe-Type Inequalities for Pseudo-Riemannian 4-Manifolds of Metric Signature $(++--)$. Geo. Ded. {\bf 87}, 65--89 (2001).
\vskip 1pt
\nref O'Neill, B. Semi-Riemannian Geometry: With Applications to Relativity. Orlando: Academic Press, 1983.
\vskip 1pt
\nref Penrose R., Rindler, W.: Spinors and Space-time, Volume 1: Two-Spinor Calculus and Relativistic Fields. Cambridge: Cambridge University Press, 1984.
\vskip 1pt
\nref Penrose, R., Rindler, W.: Spinors and Space-time, Volume 2: Spinor and Twistor Methods in Space-Time Geometry. Cambridge: Cambridge University Press, 1986.
\vskip 1pt
\nref Pleba\~nski, J.F.: Some Solutions of Complex Einstein Equations, J. Math. Phys. {\bf 16}, 2395--2402 (1975).
\vskip 1pt
\nref Pleba\~nski, J.F., Hacyan, S.: Null geodesic surfaces and Goldberg-Sachs theorem in complex Riemannian spaces. J. Math. Phys. {\bf 16}, 2403--2407 (1975).
\vskip 1pt
\nref Porteous, I.R.: Topological Geometry. Second Edition. Cambridge: Cambridge University Press, 1981.
\vskip 1pt
\nref Rashevskij, P.K.: The scalar field in a stratified space. Trudy Sem. Vektor. Tenzor. Anal. {\bf 6}, 225--248 (1948).
\vskip 1pt
\nref Rozenfeld, B.A.: On unitary and stratified spaces. Trudy Sem. Vektor. Tenzor. Anal. {\bf 7}, 260--275 (1949).
\vskip 1pt
\nref Sibata, T., Morinaga, K.: Complete and Simpler Treatment of Wave Geometry. Journal of Science of The Hiroshima University, Series A.6, 173--189 (1936).
\vskip 1pt
\nref Spivak, M.: A Comprehensive Introduction to Differential Geometry, Volume One. Second Edition. Berkeley: Publish or Perish Inc., 1979.
\vskip 1pt
\nref Spivak, M.: A Comprehensive Introduction to Differential Geometry, Volume Two. Second Edition. Berkeley: Publish or Perish Inc., 1979.
\vskip 1pt
\nref Sternberg, S.: Lectures on Differential Geometry. Second Edition. New York: Chelsea Publishing Company, 1983.
\vskip 1pt
\nref Thompson, G.: Normal Form of a metric admitting a parallel field of planes. J. Math. Phys. {\bf 33}, 4008--4010 (1992).
\vskip 1pt
\nref Warner, F.W.: Foundations of Differentiable Manifolds and Lie Groups. Graduate Texts in Mathematics 94. New York: Springer-Verlag, 1983.
\vskip 1pt
\nref Walker, A.G.: On parallel fields of partially null vector spaces. Quart. Journ. of Math. (Oxford) {\bf 20}, 135--145 (1949).
\vskip 1pt
\nref Walker, A.G.: Canonical form for a Riemannian space with a parallel field of null planes. Quart. J. Math. Oxford (2) {\bf 1}, 69--79 (1950).
\vskip 1pt
\nref Walker, A.G.: Canonical forms (II): Parallel partially null planes. Quart. J. Math. Oxford (2) {\bf 1}, 147--152 (1950).
\vskip 1pt
\nref Walker, A.G.: On Ruse's spaces of recurrent curvature. Proc. London Math. Soc. (2) {\bf 52}, 36--64 (1950).
\vskip 1pt
\nref Walker, A.G.: Connexions for parallel distributions in the large. Quart. J. Math. Oxford (2) {\bf 6}, 301--308 (1955).
\vskip 1pt
\nref Walker, A.G.: Connexions for parallel distributions in the large II. Quart. J. Math. Oxford (2) {\bf 9}, 221--231 (1958).
\vskip 1pt
\nref Yano, K.: Differential Geometry on Complex and Almost Complex Spaces. New York: MacMillan, 1965.

\bye